\documentclass[11pt]{article}
\usepackage{color}
\usepackage[colorlinks=true]{hyperref}
\usepackage[margin=1in]{geometry}
\usepackage{amssymb}
\usepackage{amsmath}
\usepackage{amsthm}
\usepackage{mathtools}
\usepackage{epsfig, graphicx}
\usepackage[boxed, noline, ruled, linesnumbered]{algorithm2e}
\newtheorem{theorem}{Theorem}
\newtheorem{lemma}{Lemma}

\newtheorem{remark}{Remark}
\numberwithin{equation}{section}

\allowdisplaybreaks % For long formula
\begin{document}
\title{Computing Multi-Eigenpairs of High-Dimensional Eigenvalue Problems Using Tensor
Neural Networks\footnote{This work was
supported in part by the National Key Research and Development Program of China
(2019YFA0709601), the National Center for Mathematics and Interdisciplinary Science, CAS.}}
\author{
Yifan Wang\footnote{LSEC, NCMIS, Institute
of Computational Mathematics, Academy of Mathematics and Systems
Science, Chinese Academy of Sciences, Beijing 100190,
China,  and School of Mathematical Sciences, University
of Chinese Academy of Sciences, Beijing 100049, China (wangyifan@lsec.cc.ac.cn).}
\ \  and \ \
Hehu Xie\footnote{LSEC, NCMIS, Institute
of Computational Mathematics, Academy of Mathematics and Systems
Science, Chinese Academy of Sciences, Beijing 100190,
China,  and School of Mathematical Sciences, University
of Chinese Academy of Sciences, Beijing 100049, China (hhxie@lsec.cc.ac.cn).}}

\date{}
\maketitle

\date{}
\maketitle

\begin{abstract}
In this paper, we propose a type of tensor-neural-network-based machine learning method to compute
multi-eigenpairs of high dimensional eigenvalue problems without Monte-Carlo procedure.
Solving multi-eigenvalues and their corresponding eigenfunctions is
one of the basic tasks in mathematical and computational physics.
With the help of tensor neural network and deep Ritz method, the high dimensional integrations
included in the loss functions of the machine learning process can be computed with high accuracy.
The high accuracy of high dimensional integrations can improve the accuracy
of the machine learning method for computing multi-eigenpairs of high dimensional
eigenvalue problems. Here, we introduce the tensor neural network and design the
machine learning method for computing multi-eigenpairs of the high dimensional eigenvalue problems.
The proposed numerical method is validated with plenty of numerical examples.

%The eigensolver adheres to the computational techniques developed in machine learning and is vastly different from
%traditional numerical methods.

%-----------------------------------------------------------------------------------------------------
\end{abstract}

\section{Introduction}
In modern sciences and engineers, there exist many high-dimensional problems,
which arise from quantum mechanics, statistical mechanics, and financial engineer.
There have appeared more and more high dimensional eigenvalue problems along with the developments
of science and engineer. Computing the eigenvalues and eigenfunctions of high dimensional
operators becomes more and more important in modern scientific computing.
The most famous and important high dimensional eigenvalue problem is the Schr\"{o}dinger equation from
quantum mechanics.
The classical numerical methods, such as finite difference, finite element, spectral,
can only solve the low dimensional Schr\"{o}dinger equations for some simple systems.
Using these classical numerical methods to solve high dimensional eigenvalue problems will
suffer from the so-called curse of dimensionality since the number of degrees of freedom and
computational complexity grows exponentially as the dimension increases.

For the high dimensional problems, the Monte-Carlo based methods attract more and more
attentions since it provides the possibility to solve high dimensional problems in the stochastic sense.
For example, there are two well known variational Monte Carlo (VMC) \cite{CeperleyChesterKalos} and
diffusion Monte Carlo (DMC)  \cite{NeedsTowlerDrummondRios} methods
for high-dimensional eigenvalue problems  in quantum mechanics.
Recently, many numerical methods based on machine learning have been proposed to solve the high-dimensional PDEs (\cite{BaymaniEffati,WeinanE, EYu,HanJentzenE, HanLuZhou,LagarisLikasPapageorgiou, LiYing,LiZhaiChen,LitsarevOseledets, RaissiPerdikarisKarniadakis, DGM, WAN, ZhangLiSchutte}). Among these machine learning methods,
neural network-based methods attract more and more attention since it can be used to build
approximations to the exact solutions of PDEs by machine learning methods.
The essential reason is that neural networks can approximate any function given enough parameters.
This type of method also provides a possible way to solve many useful high-dimensional PDEs from physics,
chemistry, biology, engineering, and so on.

Due to its universal approximation property, the fully-connected neural network (FNN) is the most widely
used architecture to build the functions for solving high-dimensional PDEs. There are several types of
FNN-based methods such as the well-known deep Ritz \cite{EYu}, deep Galerkin method \cite{DGM},
PINN \cite{RaissiPerdikarisKarniadakis}, and weak adversarial networks \cite{WAN}
for solving high-dimensional PDEs by designing different loss functions.
Among these methods, the loss functions always include computing high-dimensional
integration for the functions defined by FNN. For example,  the loss functions of the deep
Ritz method require computing the integrations on the high-dimensional domain for the functions
which is constructed by FNN. Direct numerical integration for the high-dimensional functions also meets the
``curse of dimensionality''.
Always, the Monte-Carlo method is adopted to do the high-dimensional integration
with some types of sampling methods \cite{EYu,HanZhangE}.  Due to the low convergence rate of the Monte-Carlo method,
the solutions obtained by the FNN-based numerical methods are difficult to obtain high accuracy
and stable convergence process. In other words, the Monte-Carlo method decreases computational
work in each forward propagation
by decreasing the simulation efficiency and stability of the FNN-based numerical methods
for solving high-dimensional PDEs.

In \cite{WangJinXie}, based on the deep Ritz method, a type of tensor neural network (TNN)
is proposed to build the trial
functions for solving high-dimensional PDEs.
The TNN is a function being designed by the tensor product operations on the neural networks
or by correlated low-rank CANDECOMP/PARAFAC (CP) tensor decomposition \cite{HongKoldaDuersch, KoldaBader}
approximations of FNNs.
An important advantage is that we do not need to use Monte-Carlo method to
do the integration for the functions which is constructed by TNN.   The high dimensional integration of the
TNN functions can be decomposed into one-dimensional integrations which can be computed by the
highly accurate Gauss type quadrature scheme.
The computational work for the integration of the functions by TNN is only a polynomial scale of the dimension,
which means the TNN overcomes the ``curse of dimensionality'' in some sense for
solving high-dimensional PDEs.  Furthermore, the high accuracy of high dimensional integration can
improve the TNN-based machine learning methods for solving high-dimensional problems.
The TNN-based machine learning methods have already been used to solve the smallest
eigenvalue and its corresponding eigenfunction of high-dimensional eigenvalue problems \cite{WangJinXie}
and Schr\"{o}dinger equations in \cite{WangLiaoXie}.

In this paper, we investigate the applications of TNN for computing multi-eigenpairs of
high dimensional eigenvalue problems.  The aim of this paper is to propose a universal
machine learning solver for computing multi-eigenpairs of high dimensional eigenvalue problems
based on the TNN and deep Ritz method without the a priori knowledge of the actual or guessed solutions.
Furthermore, the proposed eigensolver can compute multi-eigenpairs of the high dimensional
eigenvalue problem with obviously better accuracy than the Montor-Carlo based machine learning methods.
In this paper, we will also find that the applications of TNN and the high accuracy of high dimensional
integrations bring more choices to define the loss functions for the machine learning methods.
For example, the loss function in this paper for computing multi-eigenpairs is more like the classical
way for computing eigenvalue problems.
The proposed eigensolver provides a new way to solve high dimensional eigenvalue problems from
physics, chemistry, biology, engineering, and so on. We will provide plenty of numerical results for the
examples from the physics and material sciences.

An outline of the paper goes as follows. In Section \ref{Section_TNN}, we introduce the TNN architecture
and its approximation property. In Section \ref{Section_Method}, the TNN-based
machine learning method and corresponding numerical integration schemes
are introduced for solving multi-eigenpairs of the high dimensional
eigenvalue problems.  Section \ref{Section_Numerical} is devoted to providing numerical examples to validate the
accuracy and efficiency of the proposed numerical methods.
Finally, some concluding remarks are given in the last section.

%-----------------------------------------------------------------------------------------------------
\section{Tensor neural network architecture}\label{Section_TNN}
TNN structure, its approximation property and the computational complexity of related integration
have been discussed and investigated in \cite{WangJinXie}. In this section, we introduce the architecture of TNN.
In order to express clearly and facilitate the construction of the TNN method for solving multi-eigenpairs of
high dimensional eigenvalue problems, here, we will also elaborate on some important definitions and properties.

%\subsection{Tensor neural network architecture}

\begin{figure}[htb]
\centering
\includegraphics[width=16cm,height=12cm]{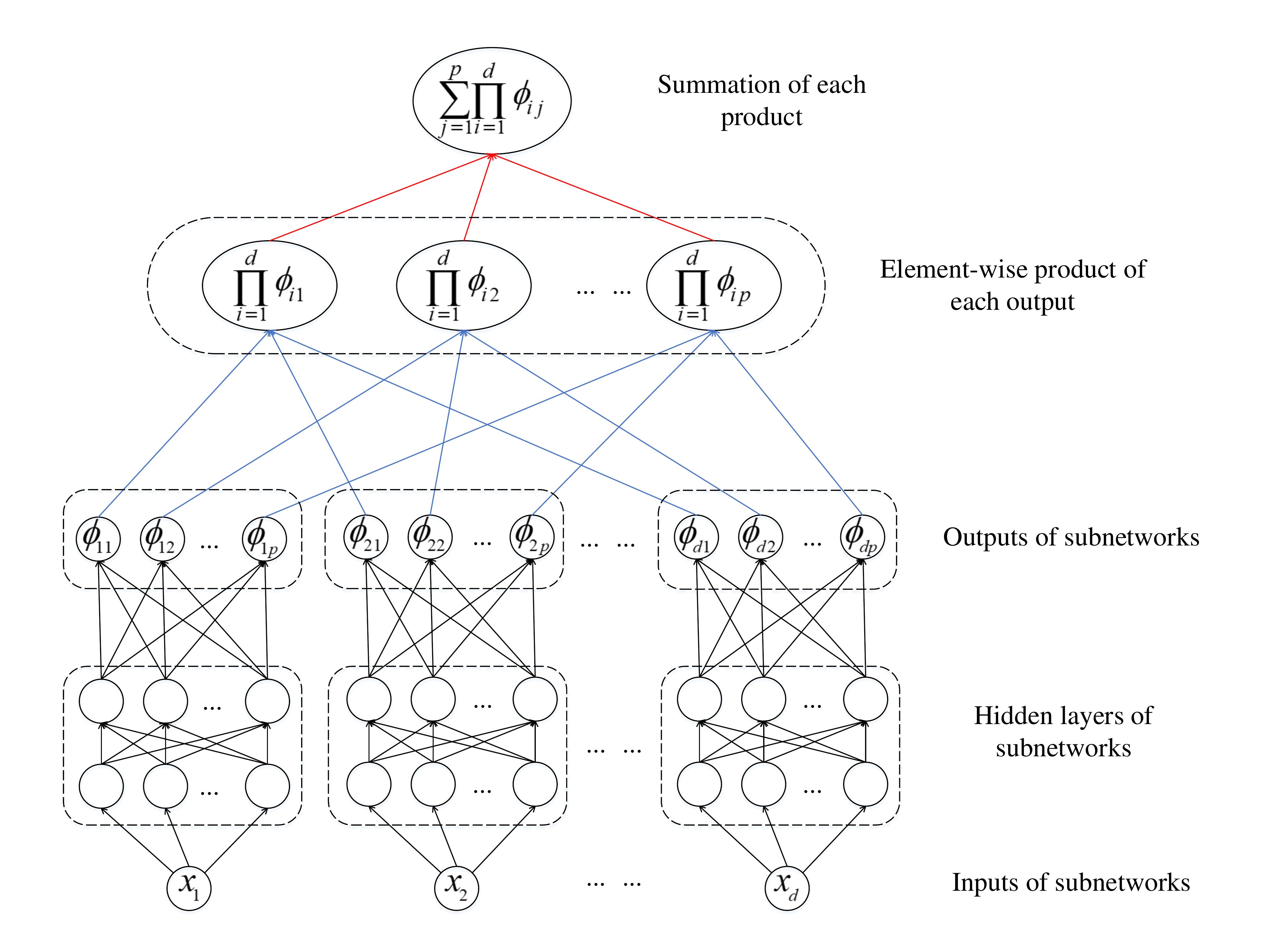}
\caption{Architecture of TNN. Black arrows mean linear transformation
(or affine transformation). Each ending node of blue arrows is obtained by taking the
scalar multiplication of all starting nodes of blue arrows that end in this ending node.
The final output of TNN is derived from the summation of all starting nodes of red arrows.}\label{TNNstructure}
\end{figure}

TNN is a neural network of low-rank structure, which is built by the tensor product of of several one-dimensional
input and multidimensional output subnetworks.
Due to the low-rank structure of TNN, an efficient and accurate quadrature scheme can be designed for the
TNN-related high dimensional integrations such as the inner product of two TNNs.
In \cite{WangJinXie}, we introduce TNN in detail and propose its numerical integration scheme with
the polynomial scale computational complexity of the dimension.
For each $i=1,2,\cdots,d$, we use $\Phi_i(x_i;\theta_i)=(\phi_{i,1}(x_i;\theta_i),\phi_{i,2}(x_i;\theta_i),\cdots,\phi_{i,p}(x_i;\theta_i))$
to denote a subnetwork that maps a set $\Omega_i\subset\mathbb R$ to $\mathbb R^p$, 
where $\Omega_i,i=1,\cdots,d,$ can be a bounded interval $(a_i,b_i)$, the whole line $(-\infty,+\infty)$
or the half line $(a_i,+\infty)$.
The number of layers and neurons in each layer, the selections of activation functions and other
hyperparameters can be different in different subnetworks. 
In this paper, in order to improve the numerical stability further, the TNN is defined as follows:
\begin{eqnarray}\label{def_TNN}
\Psi(x;\Theta)=\sum_{j=1}^pc_j\widehat\phi_{1,j}(x_1;\theta_1)\widehat\phi_{2,j}(x_2;\theta_2)
\cdots\widehat\phi_{d,j}(x_d;\theta_d)
=\sum_{j=1}^pc_j\prod_{i=1}^d\widehat\phi_{i,j}(x_i;\theta_i),
\end{eqnarray}
where $c=\{c_j\}_{j=1}^{p}$ is a set of trainable parameters, $\Theta=\{c,\theta_1,\cdots,\theta_d\}$
denotes all parameters of the whole architecture.
For $i=1,\cdots,d,j=1,\cdots,p$, $\widehat\phi_{i,j}(x_i,\theta_i)$ is a normalized functions as follows:
\begin{eqnarray}\label{eq_phi_normed}
\widehat\phi_{i,j}(x_i,\theta_i)=\frac{\phi_{i,j}(x_i,\theta_i)}{\|\phi_{i,j}(x_i,\theta_i)\|_{L^2(\Omega_i)}}.
\end{eqnarray}
In Section \ref{sec_quadrature}, we will discuss the architectures for each subnetwork $\Phi_i(x_i;\theta_i)$
in detail according to different types of $\Omega_i$.
The TNN architecture (\ref{def_TNN}) and the one defined in \cite{WangJinXie} are mathematically equivalent, but (\ref{def_TNN}) has better numerical stability during the training process. 
Figure \ref{TNNstructure} shows the corresponding architecture of TNN.
From Figure \ref{TNNstructure} and numerical tests, we can find the parameters for each rank of TNN are
correlated by the FNN, which guarantee the stability of the TNN-based machine learning methods. 
This is also an important difference from the tensor finite element methods.

In order to show the reasonableness of TNN, we now introduce the approximation property from \cite{WangJinXie}.
Since there exists the isomorphism relation between $H^m(\Omega_1\times\cdots\times\Omega_d)$
and the tensor product space $H^m(\Omega_1)\otimes\cdots\otimes H^m(\Omega_d)$, 
the process of approximating the function $f(x)\in H^m(\Omega_1\times\cdots\times\Omega_d)$
by the TNN defined as (\ref{def_TNN}) can be regarded as searching for a correlated CP decomposition structure
to approximate $f(x)$ in the space $H^m(\Omega_1)\otimes\cdots\otimes H^m(\Omega_d)$
with the rank being not greater than $p$. 
In \cite{WangJinXie} we introduce and prove the following approximation result to the functions in the space $H^m(\Omega_1\times\cdots\times\Omega_d)$ under the sense of $H^m$-norm.

\begin{theorem}\label{theorem_approximation}
Assume that each $\Omega_i$ is an interval in $\mathbb R$ for $i=1, \cdots, d$, $\Omega=\Omega_1\times\cdots\times\Omega_d$,
and the function $f(x)\in H^m(\Omega)$. Then for any tolerance $\varepsilon>0$, there exist a
positive integer $p$ and the corresponding TNN defined by (\ref{def_TNN})
such that the following approximation property holds
\begin{equation}\label{eq:L2_app}
\|f(x)-\Psi(x;\theta)\|_{H^m(\Omega)}<\varepsilon.
\end{equation}
\end{theorem}

The motivation for employing the TNN architectures is to provide high accuracy and high efficiency
in calculating variational forms of high-dimensional problems in which high-dimensional
integrations are included.
The TNN itself can approximate functions in Sobolev space with respect to $H^m$-norm.
Therefore we naturally put forward an approach to solve high-dimensional PDEs and the
ground state eigenpair in \cite{WangJinXie, WangLiaoXie}.
The major contribution of this paper is to propose a TNN-based machine learning
method to compute the leading multi-eigenpairs of high-dimensional eigenvalue problems.
We will find that the designing of the machine learning method in this paper for
computing multi-eigenpairs also shows the advantages of TNN.

\section{Machine learning method for computing multi-eigenpairs}\label{Section_Method}
This section is devoted to introducing the TNN-based machine learning method \cite{GoodfellowBengioCourville}
to compute the multi-eigenpairs of
high dimensional eigenvalue problems. We will introduce the way to build the eigenfunction approximations by
TNN functions and the quadrature schemes to compute the inner products of the TNN functions.

\subsection{Approximate eigen-subspace by TNNs}
In this subsection, we present the TNN-based discretization of eigenvalue problems for solving multi-eigenpairs.
Briefly speaking, analogous to the subspace projection method such as the finite element method,
instead of using the finite element basis, our approach uses several TNNs as the basis to
span a subspace of solution space and restrict the problem into this finite dimensional subspace.
Through a machine learning process, an optimal approximation of the eigen-subspace
represented by TNNs will be found. Then the final multi-eigenpair approximations are obtained by
solving a finite-dimensional matrix eigenvalue problem, which is similar to the Rayleigh-Ritz step in
the classical eigensolvers for the matrix eigenvalue problems \cite{Saad}.

For generality, we describe the eigenvalue and the TNN-based machine learning method in the abstract way.
More specifically, assume Sobolev spaces $\mathcal V$ and $\mathcal W$ are two Hilbert spaces and
satisfy $\mathcal V\subset \mathcal W$.
And let $a(\cdot,\cdot)$ and $b(\cdot,\cdot)$ denote two positive definite symmetric bilinear forms
on $\mathcal V\times \mathcal V$ and $\mathcal W\times \mathcal W$, respectively.
Furthermore, based on these bilinear forms, we can define the norms on the space $\mathcal V$ and $\mathcal W$ as follows
\begin{eqnarray}
\|v\|_a&=&\sqrt{a(v,v)},\ \ \ \forall v\in \mathcal V,\\
\|w\|_b&=&\sqrt{b(w,w)},\ \ \ \forall w\in \mathcal W.
\end{eqnarray}
Assume the norm $\|\cdot\|_a$ is relatively compact with respect to the norm $\|\cdot\|_b$ \cite{Conway}.

To describe our method and to build the loss function for computing leading $k$ eigenpairs briefly,
we focus on the following general eigenvalue problem: Find $(\lambda,u)\in\mathbb R\times \mathcal V$ such that $b(u,u)=1$ and
\begin{eqnarray}\label{Weak_Eigenvalue_Problem_0}
a(u,v)=\lambda b(u,v),\ \ \ \forall v\in \mathcal V.
\end{eqnarray}
It is well known that the eigenvalue problem (\ref{Weak_Eigenvalue_Problem_0})
has an eigenvalue sequence $\{\lambda_j \}$
(cf. \cite{BabuskaOsborn_Book}):
$$0<\lambda_1\leq \lambda_2\leq\cdots\leq\lambda_k\leq\cdots,\ \ \
\lim_{k\rightarrow\infty}\lambda_k=\infty,$$ and associated
eigenfunctions
$$u_1, u_2, \cdots, u_k, \cdots,$$
where $b(u_i,u_j)=\delta_{ij}$ ($\delta_{ij}$ denotes the Kronecker function).
In the sequence $\{\lambda_j\}$, the $\lambda_j$ are repeated according to their
geometric multiplicity.

The eigenvalues satisfy the minimum-maximum principle \cite{BabuskaOsborn_Book}
\begin{eqnarray}\label{min_max}
\lambda_k=\min_{\substack{\mathcal V_k\subset \mathcal V\\\dim \mathcal V_k=k}}\max_{f\in \mathcal V_k}\mathcal R(f)
=\max_{f\in{\rm span}\{u_1,\cdots,u_k\}}\mathcal R(f),
\end{eqnarray}
where $\mathcal R(f)=a(f,f)/b(f,f)$ denotes the Rayleigh quotient of the function $f$.
Let us define $\mathcal V_k={\rm span}\{v_1,\cdots,v_k\}$ and $v_i\in \mathcal V$ for $i=1,\cdots,k$ are
linearly independent of each other. Furthermore, we also define $\mathcal U_k={\rm span}\{u_1,\cdots,u_k\}$ which
denotes the eigensubspace with respect to the leading $k$ eigenvalues.
Use the terminology of the subspace approximation to the eigenvalue problems, we define the stiffness matrix
$\mathcal A(v_1,\cdots,v_k)$ and mass matrix $\mathcal B(v_1,\cdots,v_k)$ as follows
\begin{eqnarray}
&&\mathcal A(v_1,\cdots,v_k)=\Big(\mathcal A_{ij}(v_1,\cdots,v_k)\Big)_{1\leq i,j\leq k}
= \Big(a(v_i,v_j)\Big)_{1\leq i,j\leq k}  \in\mathbb R^{k\times k},\\
&&\mathcal B(v_1,\cdots,v_k)=\Big(\mathcal B_{ij}(v_1,\cdots,v_k)\Big)_{1\leq i,j\leq k} =
\Big(b(v_i,v_j)\Big)_{1\leq i,j\leq k}
\in\mathbb R^{k\times k}.
\end{eqnarray}
%where
%\begin{eqnarray}
%\mathcal A_{ij}(v_1,\cdots,v_k)=a(v_i,v_j),
%\end{eqnarray}
%and the mass matrix $\mathcal B(v_1,\cdots,v_k)$ as follows
%\begin{eqnarray}
%\mathcal B(v_1,\cdots,v_k)=\Big(\mathcal B_{ij}(v_1,\cdots,v_k)\Big)_{1\leq i,j\leq k}\in\mathbb R^{k\times k},
%\end{eqnarray}
%where
%\begin{eqnarray}
%\mathcal B_{ij}(v_1,\cdots,v_k)=b(v_i,v_j).
%\end{eqnarray}
Then due to the minimum-maximum principle, by simple derivation, the summation of the leading $k$ eigenvalues
satisfies the following optimization problem
\begin{eqnarray}
\sum_{i=1}^k\lambda_i=\min_{\substack{\mathcal V_k
={\rm span}\{v_1\cdots,v_k\}\\v_j\in \mathcal V,j=1,\cdots,k}}{\rm trace}\left(\mathcal B^{-1}(v_1,\cdots,v_k)
\mathcal A(v_1,\cdots,v_k)\right),
\end{eqnarray}
and the eigensubspace with respect to the leading $k$ eigenvalues consists with
\begin{eqnarray}
\mathcal U_k={\rm span}\{u_1,\cdots,u_k\}=\arg\min_{\mathcal V_k={\rm span}\{v_1\cdots,v_k\}
\subset \mathcal V}{\rm trace}
\Big(\mathcal B^{-1}\big(v_1,\cdots,v_k\big)\mathcal A\big(v_1,\cdots,v_k\big)\Big).
\end{eqnarray}

We approximate $\mathcal U_k$ by a space spanned by $k$ TNNs $\Psi_1(x;\Theta_1)$, $\cdots$,
$\Psi_k(x;\Theta_k)$, where each $\Psi_\ell(x;\Theta_\ell)$ is defined  by (\ref{def_TNN})
with $\Phi_{i,\ell}(x_i;\theta_{i,\ell})=(\phi_{i,1,\ell}(x_i;\theta_{i,\ell})$, $\phi_{i,2,\ell}(x_i;\theta_{i,\ell})$,
$\cdots$, $\phi_{i,p,\ell}(x_i;\theta_{i,\ell}))$ as the $i$-th subnetwork of the $\ell$-th TNN.
Then, for $\ell=1,\cdots,k$, the $\ell$-th TNN $\Psi_\ell(x;\Theta_\ell)$ is denoted as
\begin{eqnarray}\label{def_k_TNNs}
\Psi_\ell(x;\Theta_\ell)=\sum_{j=1}^{p_\ell}c_{j,\ell}\prod_{i=1}^d\widehat\phi_{i,j,\ell}(x_i;\theta_{i,\ell}),
\end{eqnarray}
where each $c_{\ell}=\{c_{j,\ell}\}_{j=1}^{p}$ is a set of trainable parameters, each $\Theta_\ell=\{c_{\ell},\theta_{1,\ell},\cdots,\theta_{d,\ell}\}$ denotes all parameters of the $\ell$-th TNN.
Similar to (\ref{eq_phi_normed}), each $\widehat\phi_{i,j,\ell}$ is normalized function. 
We can select the appropriate activation function such that all $\Psi_\ell(x;\Theta_\ell),\ell=1,\cdots,k$ belong to the
space $\mathcal V$.  Let us define $p=\max\{p_1,\cdots,p_k\}$.
These $k$ TNNs are trained using the following loss function
\begin{eqnarray}\label{loss}
&&{\rm Loss}\Big(\Psi_1(x;\Theta_1),\cdots,\Psi_k(x;\Theta_k)\Big)\nonumber\\
&&={\rm trace}\Big(\mathcal B^{-1}\big(\Psi_1(x;\Theta_1),\cdots,\Psi_k(x;\Theta_k)\big)
\mathcal A\big(\Psi_1(x;\Theta_1),\cdots,\Psi_k(x;\Theta_k)\big)\Big).
\end{eqnarray}
Since we do the inner-products of TNN in (\ref{loss}), the way to build the loss function is deep Ritz type.
In order to assemble the matrices $\mathcal A$ and $\mathcal B$ in (\ref{loss}), we need to compute the
high dimensional integrations included in the bilinear forms $a(\cdot, \cdot)$ and $b(\cdot,\cdot)$ in
(\ref{Weak_Eigenvalue_Problem_0}). 
The detailed method for assembling the matrices $\mathcal A$ and $\mathcal B$ in (\ref{loss}) will be introduced
in Section \ref{sec_quadrature}.
The loss function (\ref{loss}) is automatically differentiable thanks to packages that support
backpropagation such as TensorFlow and PyTorch.
In this paper, the gradient descent (GD) method is adopted to update all trainable parameters
\begin{eqnarray}\label{eq_update_par}
\Theta_\ell-\eta\nabla_{\Theta_\ell}{\rm Loss}\rightarrow\Theta_\ell,\ \ \ \ell=1,\cdots,k,
\end{eqnarray}
where $\eta$ is the learning rate and adjusted by the ADAM optimizer \cite{KingmaAdam}.

After sufficient training steps, we obtain a sequence of parameters $\{\Theta_1^*,\cdots,\Theta_k^*\}$
such that the loss function (\ref{loss}) arrives its minimum value under the required tolerance.
By solving the following finite-dimensional matrix eigenvalue problem
\begin{eqnarray}\label{eq_matrix_eig}
\mathcal A\big(\Psi_1(x;\Theta_1^*),\cdots,\Psi_k(x;\Theta_k^*)\big)\mathbf y
=\lambda\mathcal B\big(\Psi_1(x;\Theta_1^*),\cdots,\Psi_k(x;\Theta_k^*)\big)\mathbf y,
\end{eqnarray}
we obtain a sequence of eigenvalues
\begin{eqnarray}
0<\widehat\lambda_1\leq\widehat\lambda_2\leq\cdots\leq\widehat\lambda_k\,
\end{eqnarray}
and corresponding eigenvectors
\begin{eqnarray}
\mathbf y_1,\mathbf y_2,\cdots,\mathbf y_k,
\end{eqnarray}
where $\mathbf y_j=[y_{j,1},\cdots,y_{j,k}]^\top$ for $j=1,\cdots,k$.
Then $\widehat\lambda_1$, $\cdots$, $\widehat\lambda_k$ can be chosen as  the approximations to  the first $k$
eigenvalues $\lambda_1$, $\cdots$, $\lambda_k$ and
the corresponding approximations $\widehat u_1$, $\cdots$, $\widehat u_k$
 to the first $k$ eigenfunctions of problem (\ref{Weak_Eigenvalue_Problem_0}) can be
obtained by the following linear combination procedure
\begin{eqnarray}\label{eq_eigenpairs_estimation}
\widehat u_j(x)&=&\sum_{\ell=1}^ky_{j,\ell}\Psi_\ell(x;\Theta_\ell^*).
%\lambda_j&\approx&\widehat\lambda_j,\ \ \ j=1,\cdots,k.
\end{eqnarray}

\begin{remark}
In \cite{ZhangLiSchutte}, the loss function is defined by the summation of $k$ Rayleigh quotients
and a penalty term which constrains the $k$ neural networks to be mutually orthogonal and normalized.
Different form \cite{ZhangLiSchutte}, in this paper, the loss function (\ref{loss}) is defined without the
penalty term and the orthogonalization condition of $k$ TNNs is not imposed directly.
The reason is that the use of TNN architectures can provide a high-accuracy and
high-efficiency quadrature scheme, as we will see in Section \ref{sec_quadrature}, for
assembling the matrices in the loss function (\ref{loss}) and eigenvalue problem (\ref{eq_matrix_eig}).
The penalty term makes optimization process much more difficult and thus affects the final accuracy.
%Based on our experience with numerical examples, it seems that ignoring the
The orthogonalization condition is implicitly guaranteed by solving the optimization problem
in the machine learning method with the loss function (\ref{loss}).
%process does not affect the validity of our method.
%\comm{It is worth remarking that the loss function (\ref{loss}) we use and the loss function used
%in \cite{ZhangLiSchutte} derive essentially from the same minimum-maximum principle.
%\mathcal We do not claim to propose this loss function.
The main innovation of the present paper is for the first time to apply the TNN architecture to compute
multi-eigenpairs of high dimensional eigenvalue problems with high accuracy.
In Section \ref{sec_ex2}, we will give an example from \cite{ZhangLiSchutte} to show the advantages of
proposed approach.
\end{remark}

\subsection{Quadrature scheme for inner product}\label{sec_quadrature}
In this subsection, we provide a detailed description to calculate inner products in the loss function (\ref{loss}).
For simplicity, we are concerned with the following model problem as an example:
Find $(\lambda,u)\in\mathbb R\times H_0^1(\Omega)$ such that
\begin{equation}\label{def_EVP_strong}
\left\{\begin{aligned}
-\Delta u +V(x)u&=\lambda u, & & \text {in}\ \ \Omega, \\
u &=0, & & \text {on}\ \partial \Omega,
\end{aligned}\right.
\end{equation}
where $\Omega=\Omega_1\times\cdots\times\Omega_d$, 
each $\Omega_i,i=1,\cdots,d,$ can be a bounded interval $(a_i,b_i)$, the whole line $(-\infty,+\infty)$
or the half line $(a_i,+\infty)$, $V(x)\in L^2(\Omega)$ is a potential function.
We assume that the potential $V(x)$ has separated representation in the
tensor product space $L^2(\Omega_1)\otimes\cdots\otimes L^2(\Omega_d)$ as follows
\begin{eqnarray}\label{def_V(x)}
V(x)=\sum_{j=1}^q\prod_{i=1}^dV_{i,j}(x_i),
\end{eqnarray}
where $V_{i,j}(x_i)\in L^2(\Omega_i)$.
Then the equivalent variational form of the eigenvalue problem (\ref{def_EVP_strong}) can be defined as follows:
Find $(\lambda,u)\in\mathbb R\times H_0^1(\Omega)$ such that
\begin{eqnarray}\label{Weak_Eigenvalue_Problem}
a(u,v) = \lambda b(u,v), \ \ \ \ \forall v\in H_0^1(\Omega),
\end{eqnarray}
where the bilinear forms $a(\cdot,\cdot)$ and $b(\cdot,\cdot)$ are defined as follows
\begin{eqnarray}\label{inner_product_a_b}
a(u,v)=\int_{\Omega}(\nabla u\cdot\nabla v +Vuv)d\Omega,
 \ \ \ \  \ \ b(u,v) = \int_{\Omega}uv d\Omega.
\end{eqnarray}
Due to the definition of $k$ TNNs (\ref{def_k_TNNs}) and (\ref{def_V(x)}), entries of the matrix $\mathcal A$
in the loss function (\ref{loss}) have following expansions
\begin{eqnarray}\label{eq_A_mn}
\mathcal A_{mn}&=&a\big(\Psi_m(x;\Theta_m),\Psi_n(x;\Theta_n)\big)\nonumber\\
&=&\sum_{s=1}^d\sum_{j=1}^{p_m}\sum_{\ell=1}^{p_n}c_{j,m}c_{\ell,n}\prod_{i\neq s}^d\int_{\Omega_i}\phi_{i,j,m}(x_i;\theta_{i,m})\phi_{i,\ell,n}(x_i;\theta_{i,n})dx_i\nonumber\\
&&\ \ \ \ \
\cdot\int_{\Omega_s}\frac{\partial \phi_{s,j,m}}{\partial x_s}(x_s;\theta_{s,m})\frac{\partial\phi_{s,\ell,n}}{\partial x_s}(x_s;\theta_{s,n})dx_s\nonumber\\
&&\ \ +\sum_{s=1}^q\sum_{j=1}^{p_m}\sum_{\ell=1}^{p_n}c_{j,m}c_{\ell,n}\prod_{i=1}^d
\int_{\Omega_i}V_{i,s}(x_i)\phi_{i,j,m}(x_i;\theta_{i,m})
\phi_{i,\ell,m}(x_i;\theta_{i,n})dx_i,
\end{eqnarray}
and entries of the matrix $\mathcal B$ can be expanded as follows
\begin{eqnarray}\label{eq_B_mn}
\mathcal B_{mn}&=&b\big(\Psi_m(x;\Theta_m),\Psi_n(x;\Theta_n)\big)\nonumber\\
&=&\sum_{j=1}^{p_m}\sum_{\ell=1}^{p_n}c_{j,m}c_{\ell,n}\prod_{i=1}^d\int_{\Omega_i}\phi_{i,j,m}(x_i;\theta_{i,m})
\phi_{i,\ell,n}(x_i;\theta_{i,n})dx_i.
\end{eqnarray}

Since only one-dimensional integrations are involved in (\ref{eq_A_mn}) and (\ref{eq_B_mn}),
there is no need to use the Monte Carlo procedure to do these high dimensional integrations. 
In order to guarantee the high accuracy of the high dimensional integrations included in the loss functions,
we should use the high order one dimensional quadrature scheme,
such as Gauss-type rules \cite{Gautschi}, for the integrations in (\ref{eq_A_mn}) and (\ref{eq_B_mn}).

\subsubsection{Legendre-Gauss quadrature scheme for bounded domain}
Without loss of generality,
we decompose each $\Omega_i$ into  $M_i$ equal subintervals with length $h_i = |\Omega_i|/M_i$ and choose $N_i$
Legendre-Gauss points in each subinterval.
Denote the total quadrature points and corresponding weights as follows
\begin{eqnarray}\label{quad_p_w}
\left\{x_i^{(n_i)}\right\}_{n_i=1}^{M_iN_i},\ \ \ \left\{w_i^{(n_i)}\right\}_{n_i=1}^{M_iN_i},\ \ \ i=1,2,\cdots,d,
\end{eqnarray}
and define $N=\max\{N_1,\cdots,N_d\}$, $\underline N = \min\{N_1,\cdots,N_d\}$, $M=\max\{M_1,\cdots,M_d\}$ and $\underline M = \min\{M_1,\cdots,M_d\}$.
Then using the quadrature scheme (\ref{quad_p_w}), entries of the matrix $\mathcal A$ defined by (\ref{eq_A_mn})
have the following numerical format
\begin{eqnarray}\label{eq_quad_A_mn}
&&\mathcal A_{mn}=a\big(\Psi_m(x;\Theta_m),\Psi_n(x;\Theta_n)\big)\nonumber\\
&&\approx\sum_{s=1}^d\sum_{j=1}^{p_m}\sum_{\ell=1}^{p_n}c_{j,m}c_{\ell,n}\prod_{i\neq s}^d\sum_{n_i=1}^{M_iN_i}w_i^{(n_i)}\phi_{i,j,m}(x_i^{(n_i)};\theta_{i,m})
\phi_{i,\ell,n}(x_i^{(n_i)};\theta_{i,n})\nonumber\\
&&\cdot\sum_{n_s=1}^{M_sN_s}w_s^{(n_s)}\frac{\partial \phi_{s,j,m}}{\partial x_s}(x_s^{(n_s)};\theta_{s,m})\frac{\partial\phi_{s,\ell,n}}{\partial x_s}(x_s^{(n_s)};\theta_{s,n})\nonumber\\
&&+\sum_{s=1}^q\sum_{j=1}^{p_m}\sum_{\ell=1}^{p_n}c_{j,m}c_{\ell,n}\prod_{i=1}^d\sum_{n_i=1}^{M_iN_i}
w_i^{(n_i)}V_{i,s}(x_i^{(n_i)})\phi_{i,j,m}(x_i^{(n_i)};\theta_{i,m})\phi_{i,\ell,m}(x_i^{(n_i)};\theta_{i,n}),
\end{eqnarray}
and entries of the matrix $\mathcal B$ defined by (\ref{eq_B_mn}) can be computed as follows
\begin{eqnarray}\label{eq_quad_B_mn}
\mathcal B_{mn}&=&b\big(\Psi_m(x;\Theta_m),\Psi_n(x;\Theta_n)\big)\nonumber\\
&\approx&\sum_{j=1}^{p_m}\sum_{\ell=1}^{p_n}c_{j,m}c_{\ell,n}\prod_{i=1}^d\sum_{n_i=1}^{M_iN_i}
w_i^{(n_i)}\phi_{i,j,m}(x_i^{(n_i)};\theta_{i,m})\phi_{i,\ell,n}(x_i^{(n_i)};\theta_{i,n}).
\end{eqnarray}
Then the complete TNN-based algorithm for computing the first $k$ eigenpairs can be summarized in Algorithm \ref{Algorithm_1}.

\begin{algorithm}[htb!]
\caption{TNN-based method for the first $k$ eigenpairs }\label{Algorithm_1}
\begin{enumerate}
\item Initialization step: Build $k$ initial TNNs $\Psi_1^{(0)}(x;\Theta_1^{(0)}),\cdots,\Psi_k^{(0)}(x;\Theta_k^{(0)})$ as in (\ref{def_k_TNNs}),
    maximum training steps $L$, learning rate $\eta$, quadrature points and weights (\ref{quad_p_w}).

\item Assemble matrices $\mathcal A\big(\Psi_1^{(\ell)}(x;\Theta_1^{(\ell)}),\cdots,\Psi_k^{(\ell)}(x;\Theta_k^{(\ell)})\big)$
and $\mathcal B\big(\Psi_1^{(\ell)}(x;\Theta_1^{(\ell)}),\cdots,\Psi_k^{(\ell)}(x;\Theta_k^{(\ell)})\big)$
according to quadrature schemes (\ref{eq_quad_A_mn}) and (\ref{eq_quad_B_mn}), respectively.

\item Compute the value of the loss function (\ref{loss}),  where the matrix $\mathcal B^{-1}\mathcal A=:\mathcal C$
is obtained by solving the matrix equation $\mathcal B\mathcal C=\mathcal A$.

\item Compute the gradient of the loss function with respect to parameters $\{\Theta_1,\cdots,\Theta_k\}$ by automatic differentiation and update parameters $\{\Theta_1,\cdots,\Theta_k\}$ by (\ref{eq_update_par}).

\item Set $\ell=\ell+1$ and go to Step 2 for the next step until $\ell=L$.

\item Post-processing step: Solve the matrix eigenvalue problrm (\ref{eq_matrix_eig}) and compute the first $k$ eigenpair approximations by (\ref{eq_eigenpairs_estimation}).
\end{enumerate}
\end{algorithm}

It is worth mentioning that the quadrature schemes (\ref{eq_quad_A_mn}) and (\ref{eq_quad_B_mn})
which use $1$-dimensional $M_iN_i$ quadrature points in each dimension are actually equivalent
to implementing the following $d$-dimensional full
tensor quadrature scheme
\begin{eqnarray}
\mathcal A_{mn}&\approx&\sum_{n\in\mathcal N}w^{(n)}\nabla\Psi_m(x^{(n)};\Theta_m)
\cdot\nabla\Psi_n(x^{(n)};\Theta_n)\nonumber\\
&&+\sum_{n\in\mathcal N}w^{(n)}V(x^{(n)})\Psi_m(x^{(n)};\Theta_m)\Psi_n(x^{(n)};\Theta_n)\nonumber,\\
\mathcal B_{mn}&\approx&\sum_{n\in\mathcal N}w^{(n)}\Psi_m(x^{(n)};\Theta_m)\Psi_n(x^{(n)};\Theta_n),
\end{eqnarray}
where $d$-dimensional quadrature points and weights are defined as follows
\begin{eqnarray}\label{def_tensor_quad_p_w}
\left.
\begin{array}{rcl}
\Big\{x^{(n)}\Big\}_{n\in\mathcal N}&=&\left\{\left\{x_1^{(n_1)}\right\}_{n_1=1}^{M_1N_1}\times\ \left\{x_2^{(n_2)}\right\}_{n_2=1}^{M_2N_2}\times\ \cdots\times\ \left\{x_d^{(n_d)}\right\}_{n_d=1}^{M_dN_d}\right\},\\
\Big\{w^{(n)}\Big\}_{n\in\mathcal N}&=&\left\{\left\{w_1^{(n_1)}\right\}_{n_1=1}^{M_1N_1}\times \left\{w_2^{(n_2)}\right\}_{n_2=1}^{M_2N_2}\times \cdots, \times  \left\{w_d^{(n_d)}\right\}_{n_d=1}^{M_dN_d}\right\}.
\end{array}
\right.
\end{eqnarray}
The quadrature points (\ref{def_tensor_quad_p_w}) for general $d$-dimensional integrand has accuracy
$\mathcal O(h^{2\underline N}/(2\underline N)!)$ but $\mathcal N$ has $\mathcal O(M^dN^d)$ elements.
This becomes too large for $d\gg 1$ and makes the approximation of integration intractable.
Thanks to TNN having the low-rank tensor type structure, using the full tensor quadrature scheme
for the inner product of two TNNs has the splitting schemes (\ref{eq_quad_A_mn}) and (\ref{eq_quad_B_mn}).
In quadrature schemes (\ref{eq_quad_A_mn}) and (\ref{eq_quad_B_mn}), assembling matrices $\mathcal A$
and $\mathcal B$ costs only $\mathcal O(d^2p^2MN+dp^2qMN)$ and $\mathcal O(dp^2MN)$ operations, respectively.
Furthermore, due to the equivalence of using full tensor quadrature points (\ref{def_tensor_quad_p_w}), quadrature schemes (\ref{eq_quad_A_mn}) and (\ref{eq_quad_B_mn}) also have accuracy $\mathcal O(h^{2\underline N}/(2\underline N)!)$.
This means the TNN-based method proposed in this paper overcomes the ``curse of dimensionality''
in the sense of numerical integration.

For generality, in the next two sections, we will also introduce the integration of TNN on unbounded domains. 
The corresponding complexity analysis will be omitted since they are similar to that in this section.

\subsubsection{Hermite-Gauss quadrature scheme for the whole line}
For case $\Omega_i=(-\infty,+\infty)$, we use Hermite-Gauss quadrature to assemble matrix $\mathcal A$ and $\mathcal B$.
Hermite-Gauss quadrature scheme satisfies following property.
\begin{lemma}\cite[Theorem 7.3]{ShenTangWang}
Let $\{x^{(k)}\}_{k=0}^N$ be the zeros of the $(N+1)$-th order Hermite polynomial $H_{N+1}(x)$,
and let $\{w^{(k)}\}_{k=1}^N$ be given by
\begin{eqnarray}
w^{(k)}=\frac{\sqrt{\pi}2^NN!}{(N+1)H_N^2(x^{(k)})},\ \ \ 0\leq k\leq N.
\end{eqnarray}
Then we have
\begin{eqnarray}\label{hermite_gauss_form}
\int_{-\infty}^{+\infty}p(x)e^{-x^2}dx=\sum_{k=0}^Np(x^{(k)})w^{(k)},\ \ \ \forall p\in P_{2N+1},
\end{eqnarray}
where $P_{2N+1}$ denotes the set of polynomial functions of degree less than $2N+1$.
\end{lemma}

For using Hermite-Gauss quadrature, the $i$-th subnetwork of $\ell$-th TNN is defined as follows
\begin{eqnarray*}
\Phi_{i,\ell}(x_i;\theta_{i,\ell})&=&\big(\phi_{i,1,\ell}(x_i;\theta_{i,\ell}),\phi_{i,2,\ell}(x_i;\theta_{i,\ell}),
\cdots,\phi_{i,p,\ell}(x_i;\theta_{i,\ell})\big)\\
&=&e^{-\frac{\beta_i^2x_i^2}{2}}\big(\varphi_{i,1,\ell}(\beta_ix_i;\theta_{i,\ell}),\varphi_{i,2,\ell}
(\beta_ix_i;\theta_{i,\ell}),\cdots,\varphi_{i,p,\ell}(\beta_ix_i;\theta_{i,\ell}) \big),
\end{eqnarray*}
where $\beta_i$ is trainable parameter and $\varphi_{i,\ell}=\big(\varphi_{i,1,\ell}(\beta_ix_i;\theta_{i,\ell}),\varphi_{i,2,\ell}
(\beta_ix_i;\theta_{i,\ell}),\cdots,\varphi_{i,p,\ell}(\beta_ix_i;\theta_{i,\ell})\big)$ is a fully-connected neural network which maps $\mathbb R$ to $\mathbb R^{p}$.
Then the $\ell$-th TNN can be written as
\begin{eqnarray}
\Psi_{\ell}(x;\Theta_{\ell})=\sum_{j=1}^pc_{j,\ell}\prod_{i=1}^de^{-\frac{\beta_i^2x_i^2}{2}}
\widehat\varphi_{i,j,\ell}(\beta_ix_i;\theta_{i,\ell}),
\end{eqnarray}
where $\widehat\varphi_{i,j,\ell}$ satisfies normalization property $\|e^{-\frac{\beta_i^2x_i^2}{2}}\widehat\varphi_{i,j,\ell}(\beta_ix_i;\theta_{i,\ell})\|_{L^2(\Omega_i)}=1$
and is defined as follows
\begin{eqnarray}
\widehat\varphi_{i,j,\ell}(\beta_ix_i;\theta_{i,\ell})
=\frac{\varphi_{i,j,\ell}(\beta_ix_i;\theta_{i,\ell})}{\Big\|e^{-\frac{\beta_i^2x_i^2}{2}}
\varphi_{i,j,\ell}(\beta_ix_i;\theta_{i,\ell})\Big\|_{L^2(\Omega_i)}}.
\end{eqnarray}
For $i=1,\cdots,d$, let $\{z_i^{(k_i)}\}_{k_i=1}^{N_i}$ and $\{w_i^{(k_i)}\}_{k_i=1}^{N_i}$
be the Hermite-Gauss quadrature points and weights, respectively.
Under coordinate transformation $z_i=\beta_ix_i$, using the quadrature scheme (\ref{hermite_gauss_form}), entries of the matrix $\mathcal A$ defined by (\ref{eq_A_mn}) have the following numerical format
\begin{eqnarray}\label{Scheme_11}
&&\mathcal A_{mn}=\sum_{i=1}^d\sum_{j=1}^{p_m}\sum_{\ell=1}^{p_n}c_{j,m}c_{\ell,n}\prod_{i\neq s}^d\frac{1}{\beta_i}\sum_{k_i=1}^{N_i}\widehat\varphi_{i,j,m}(z_i^{(k_i)};\theta_{i,m})
\widehat\varphi_{i,\ell,n}(z_i^{(k_i)};\theta_{i,n})w_i^{(k_i)}\nonumber\\
&&\cdot\frac{1}{\beta_s}\sum_{k_s}^{N_s}\big(\widehat\varphi_{s,j,m}^\prime(z_s^{(k_s)};\theta_{s,m})
-\beta_s\widehat\varphi_{s,j,m}(z_s^{(k_s)};\theta_{s,m})\big)\nonumber\\
&&\cdot\big(\widehat\varphi_{s,j,m}^\prime(z_s^{(k_s)};\theta_{s,m})
-\beta_s\widehat\varphi_{s,j,m}(z_s^{(k_s)};\theta_{s,m})\big)w_s^{(k_s)}\nonumber\\
&&+\sum_{s=1}^q\sum_{j=1}^{p_m}\sum_{\ell=1}^{p_n}c_{j,m}c_{\ell,n}
\prod_{i=1}^d\frac{1}{\beta_i}\sum_{k_i=1}^{N_i}V_{i,s}(\frac{z_i^{(k_i)}}{\beta_i})
\widehat\varphi_{i,j,m}(z_i^{(k_i)};\theta_{i,m})\widehat\varphi_{i,\ell,n}(z_i^{(k_i)};\theta_{i,n})w_i^{(k_i)}.
\end{eqnarray}
And entries of the matrix $\mathcal B$ defined by (\ref{eq_B_mn}) have the following numerical format
\begin{eqnarray}\label{Scheme_10}
\mathcal B_{mn}=\sum_{j=1}^{p_m}\sum_{\ell=1}^{p_n}c_{j,m}c_{\ell,n}\prod_{i=1}^d\frac{1}{\beta_i}
\sum_{k_i=1}^{N_i}\widehat\varphi_{i,j,m}(z_i^{(k_i)};\theta_{i,m})
\widehat\varphi_{i,\ell,n}(z_i^{(k_i)};\theta_{i,n})w_i^{(k_i)}.
\end{eqnarray}

\subsubsection{Laguerre-Gauss quadrature scheme for the half line}
For case $\Omega_i=(0,+\infty)$, we use Laguerre-Gauss quadrature to assemble matrices $\mathcal A$ and $\mathcal B$.
Laguerre-Gauss quadrature satisfies the following theorem
\begin{lemma}\cite[Theorem 7.1]{ShenTangWang}
Let $\{x^{(k)}\}_{k=0}^N$ be the zeros of the $(N+1)$-th order Laguerre polynomial
$\mathcal L_{N+1}(x)$, and let $\{w^{(k)}\}_{k=0}^N$ be given by
\begin{eqnarray}
w^{(k)}=\frac{\Gamma(N+1)}{(N+1)(N+1)!}\frac{x^{(k)}}{\big[\mathcal L_N(x^{(k)})\big]^2},\ \ \ 0\leq 0\leq N.
\end{eqnarray}
Then we have
\begin{eqnarray}\label{lagurerre_gauss_form}
\int_{-\infty}^{+\infty}p(x)e^{-x}dx=\sum_{k=0}^Np(x^{(k)})w^{(k)},\ \ \ \forall p\in P_{2N+1},
\end{eqnarray}
where $P_{2N+1}$ denotes the set of polynomial functions of degree less than $2N+1$.
\end{lemma}

For using Laguerre-Gauss quadrature, the $i$-th subnetwork of the $\ell$-th TNN is defined as follows
\begin{eqnarray*}
\Phi_{i,\ell}(x_i;\theta_{i,\ell})&=&\big(\phi_{i,1,\ell}(x_i;\theta_{i,\ell}),\phi_{i,2,\ell}(x_i;\theta_{i,\ell}),
\cdots,\phi_{i,p,\ell}(x_i;\theta_{i,\ell})\big)\\
&=&e^{-\frac{\beta_ix_i}{2}}\big(\varphi_{i,1,\ell}
(\beta_ix_i;\theta_{i,\ell}),\varphi_{i,2,\ell}(\beta_ix_i;\theta_{i,\ell}),
\cdots,\varphi_{i,p,\ell}(\beta_ix_i;\theta_{i,\ell}) \big),
\end{eqnarray*}
where $\beta_i$ is trainable parameter and $\varphi_{i,\ell}=\big(\varphi_{i,1,\ell}(\beta_ix_i;\theta_{i,\ell}),\varphi_{i,2,\ell}(\beta_ix_i;\theta_{i,\ell}),
\cdots,\varphi_{i,p,\ell}(\beta_ix_i;\theta_{i,\ell})\big)$ is a fully-connected neural network which maps $\mathbb R$
to $\mathbb R^{p}$. Then the $\ell$-th TNN can be written as
\begin{eqnarray}
\Psi_{\ell}(x;\Theta_{\ell})=\sum_{j=1}^pc_{j,\ell}\prod_{i=1}^de^{-\frac{\beta_ix_i}{2}}
\widehat\varphi_{i,j,\ell}(\beta_ix_i;\theta_{i,\ell}),
\end{eqnarray}
where $\widehat\varphi_{i,j,\ell}$ satisfies normalization property $\|e^{-\frac{\beta_ix_i}{2}}\widehat\varphi_{i,j,\ell}(\beta_ix_i;\theta_{i,\ell})\|_{L^2(\Omega_i)}=1$
and is defined as follows
\begin{eqnarray}
\widehat\varphi_{i,j,\ell}(\beta_ix_i;\theta_{i,\ell})=\frac{\varphi_{i,j,\ell}
(\beta_ix_i;\theta_{i,\ell})}{\Big\|e^{-\frac{\beta_ix_i}{2}}
\varphi_{i,j,\ell}(\beta_ix_i;\theta_{i,\ell})\Big\|_{L^2(\Omega_i)}}.
\end{eqnarray}
For $i=1,\cdots,d$, let $\{z_i^{(k_i)}\}_{k_i=1}^{N_i}$ and $\{w_i^{(k_i)}\}_{k_i=1}^{N_i}$
be Laguerre-Gauss quadrature points and weights, respectively.
Under coordinate transformation $z_i=\beta_ix_i$, using the quadrature scheme (\ref{hermite_gauss_form}),
entries of the matrix $\mathcal A$ defined by (\ref{eq_A_mn}) have the following numerical format
\begin{eqnarray}\label{Scheme_21}
&&\mathcal A_{mn}=\sum_{i=1}^d\sum_{j=1}^{p_m}\sum_{\ell=1}^{p_n}c_{j,m}c_{\ell,n}\prod_{i\neq s}^d\frac{1}{\beta_i}\sum_{k_i=1}^{N_i}\widehat\varphi_{i,j,m}(z_i^{(k_i)};\theta_{i,m})
\widehat\varphi_{i,\ell,n}(z_i^{(k_i)};\theta_{i,n})w_i^{(k_i)}\nonumber\\
&&\cdot\frac{1}{\beta_s}\sum_{k_s}^{N_s}\big(\widehat\varphi_{s,j,m}^\prime(z_s^{(k_s)};\theta_{s,m})
-\frac{1}{2}\widehat\varphi_{s,j,m}(z_s^{(k_s)};\theta_{s,m})\big)\nonumber\\
&&\cdot\big(\widehat\varphi_{s,j,m}^\prime(z_s^{(k_s)};\theta_{s,m})-\frac{1}{2}
\widehat\varphi_{s,j,m}(z_s^{(k_s)};\theta_{s,m})\big)w_s^{(k_s)}\nonumber\\
&&+\sum_{s=1}^q\sum_{j=1}^{p_m}\sum_{\ell=1}^{p_n}c_{j,m}c_{\ell,n}
\prod_{i=1}^d\frac{1}{\beta_i}\sum_{k_i=1}^{N_i}V_{i,s}(\frac{z_i^{(k_i)}}{\beta_i})
\widehat\varphi_{i,j,m}(z_i^{(k_i)};\theta_{i,m})\widehat\varphi_{i,\ell,n}(z_i^{(k_i)};\theta_{i,n})w_i^{(k_i)}.
\end{eqnarray}
And entries of the matrix $\mathcal B$ defined by (\ref{eq_B_mn}) have the following numerical format
\begin{eqnarray}\label{Scheme_20}
\mathcal B_{mn}=\sum_{j=1}^{p_m}\sum_{\ell=1}^{p_n}c_{j,m}c_{\ell,n}\prod_{i=1}^d\frac{1}{\beta_i}
\sum_{k_i=1}^{N_i}\widehat\varphi_{i,j,m}(z_i^{(k_i)};\theta_{i,m})
\widehat\varphi_{i,\ell,n}(z_i^{(k_i)};\theta_{i,n})w_i^{(k_i)}.
\end{eqnarray}

\section{Numerical examples}\label{Section_Numerical}
In this section, we provide several examples to investigate the performance of the TNN-based eigensolver
proposed in this paper.
To show the convergence behavior and accuracy of our method, we define the relative errors
for the approximated eigenvalues $\widehat\lambda_\ell$ and eigenfunctions $\mathbf y_\ell$ as follows
\begin{eqnarray}\label{relative_errors}
{\rm err}_{\lambda,\ell}:=\frac{|\widehat\lambda_\ell-\lambda_\ell|}{|\lambda_\ell|},\
{\rm err}_{L^2,\ell}:=\frac{\|u_\ell-\mathcal Q_\ell u_\ell\|_{L^2(\Omega)}}{\|u_\ell\|_{L^2(\Omega)}},
\  {\rm err}_{H^1,\ell}:=\frac{\left|u_\ell-\mathcal P_\ell u_\ell\right|_{H^1(\Omega)}}{\left|u_\ell\right|_{H^1(\Omega)}},\
\ell = 1, \cdots, k,
\end{eqnarray}
where $\lambda_\ell$ and $u_\ell$ are reference eigenvalues and eigenfunctions.
In the first two examples, the reference eigenpairs are obtained by high order finite element methods on the meshes
with a sufficiently small mesh size.
As for the harmonic oscillator problems, the exact eigenpairs are chosen as references.
In (\ref{relative_errors}), $\mathcal P_\ell:H_0^1(\Omega)\rightarrow M_\ell$ and
$\mathcal Q:H_0^1(\Omega)\rightarrow M_\ell$ denote $L^2(\Omega)$ projection operator and $H^1(\Omega)$ projection operator, respectively, to the approximate eigenspace corresponding to the eigenvalue $\lambda_\ell$.
These two projection operators are defined as follows
\begin{eqnarray}
\left\langle\mathcal P_\ell u,v\right\rangle_{L^2}&=&\left\langle u,v\right\rangle_{L^2}:=\int_\Omega uvdx,
\ \ \ \forall v\in M_\ell\ \ {\rm for}\ u\in H_0^1(\Omega),\\
\left\langle\mathcal Q_\ell u,v\right\rangle_{H^1}&=&\left\langle u,v\right\rangle_{H^1}
:=\int_\Omega\nabla u\cdot\nabla vdx,
\ \ \ \forall v\in M_\ell\ \ {\rm for}\ u\in H_0^1(\Omega).
\end{eqnarray}
In implementation, we use the quadrature scheme similar to that in (\ref{eq_quad_A_mn}) and (\ref{eq_quad_B_mn}),
(\ref{Scheme_11}) and (\ref{Scheme_10}), (\ref{Scheme_21}) and (\ref{Scheme_20}),
to compute ${\rm err}_{L^2,\ell}$ and ${\rm err}_{H^1,\ell}$ with the same tensor product quadrature
points and weights as computing
the loss functions if the reference solution $u_\ell(x)$ has a low-rank representation, otherwise
we only report ${\rm err}_{\lambda,\ell}$.
With the help of Theorem 2 in \cite{WangJinXie}, the high efficiency and accuracy for
computing ${\rm err}_{L^2,\ell}$ and ${\rm err}_{H^1,\ell}$ can be guaranteed.

% \subsection{An example of linear Schr\"{o}dinger equation}
% In this subsection, we consider the same eigenvalue problem as in \cite{HanLuZhou} and \cite{LiYing} associated
% with the following linear Schr\"{o}dinger operator
% \begin{eqnarray}\label{Linear_Schrodinger_Operator}
% \mathcal L=-\Delta+4\pi^2\sum_{i=1}^dc_i\cos(2\pi x_i),
% \end{eqnarray}
% with the periodic boundary condition.
% The periodical property of this problem ensures that we can solve it in the domain $[0,1]^d$.
% Since the operator $\mathcal L$ is essentially decoupled, \cite{HanLuZhou} describes a decoupled way to
% obtain the reference eigenvalues and eigenfunctions by using the spectral method to solve the decoupled
% one dimensional eigenvalue problem in each dimension.
% In this paper, instead of the spectral method, we use a similar process to obtain reference
% eigenvalues and eigenfunctions by using the high order finite element method to solve the decoupled one dimensional
% eigenvalue problem in each dimension.
% As same as \cite{HanLuZhou} and \cite{LiYing}, the coefficients $c_i$ in (\ref{Linear_Schrodinger_Operator})
% take values in $[0,0.2]$ randomly.

\subsection{Infinitesimal generators of metastable diffusion processes}\label{sec_ex2}
In this subsection, we study the eigenvalue problem from \cite{ZhangLiSchutte} associated with the operator
\begin{eqnarray}
\mathcal L_d=\nabla V_d\cdot\nabla-\Delta.
\end{eqnarray}
The potential functions $V_d:\mathbb R^d\rightarrow\mathbb R$ for $d=2,50,100$ are defined as follows
\begin{eqnarray}
V_d(x)=V(\theta)+2(r-1)^2+5e^{-5r^2}+5\sum_{i=3}^dx_i^2,\ \ \ \forall x=(x_1,x_2,\cdots,x_d)\in\mathbb R^d,
\end{eqnarray}
where $(r,\theta)\in[0,+\infty)\times[-\pi,\pi)$ denotes the polar coordinates which relate to
the first two dimensional Eucild space $(x_1,x_2)\in\mathbb R^2$ by
\begin{eqnarray}
x_1=r\cos\theta,\ \ \ x_2=r\sin\theta,
\end{eqnarray}
and $V:[-\pi,\pi)\rightarrow\mathbb R$ is a double-well potential function which is defined as follows
\begin{eqnarray}
V(\theta)=
\left\{
\begin{aligned}
&\left[1-\left(\frac{3\theta}{\pi}+1\right)^2\right]^2,&\ \ \ &\theta\in \left[-\pi,-\frac{\pi}{3}\right),\\
&\frac{1}{5}\left(3-2\cos(3\theta)\right),&\ \ \ &\theta\in\left[-\frac{\pi}{3},\frac{\pi}{3}\right),\\
&\left[1-\left(\frac{3\theta}{\pi}-1\right)^2\right]^2,&\ \ \ &\theta\in\left[\frac{\pi}{3},\pi\right).
\end{aligned}
\right.
\end{eqnarray}

\begin{figure}[htb!]
\centering
\includegraphics[width=3cm,height=3cm]{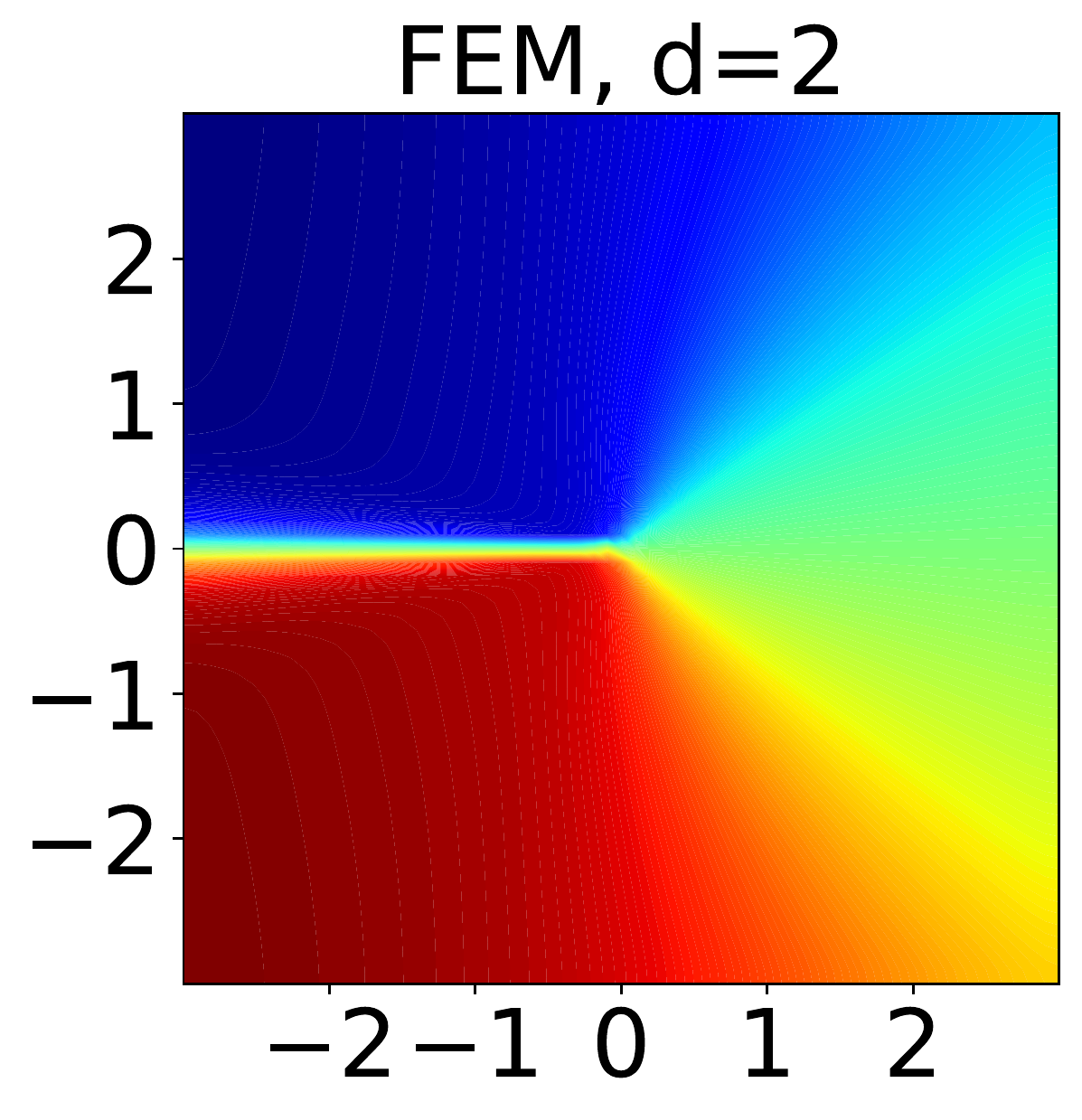}
\includegraphics[width=3cm,height=3cm]{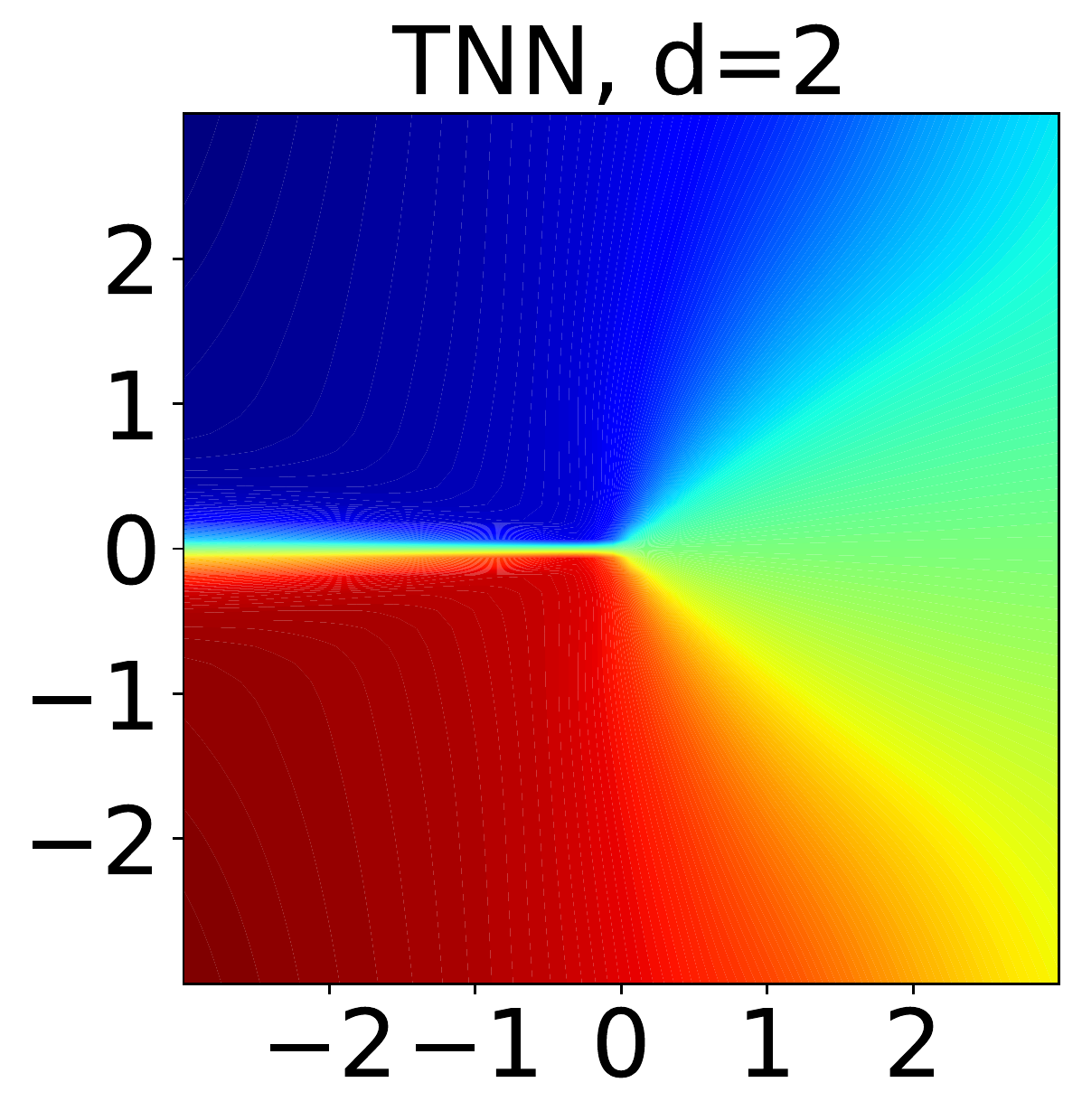}
\includegraphics[width=3cm,height=3cm]{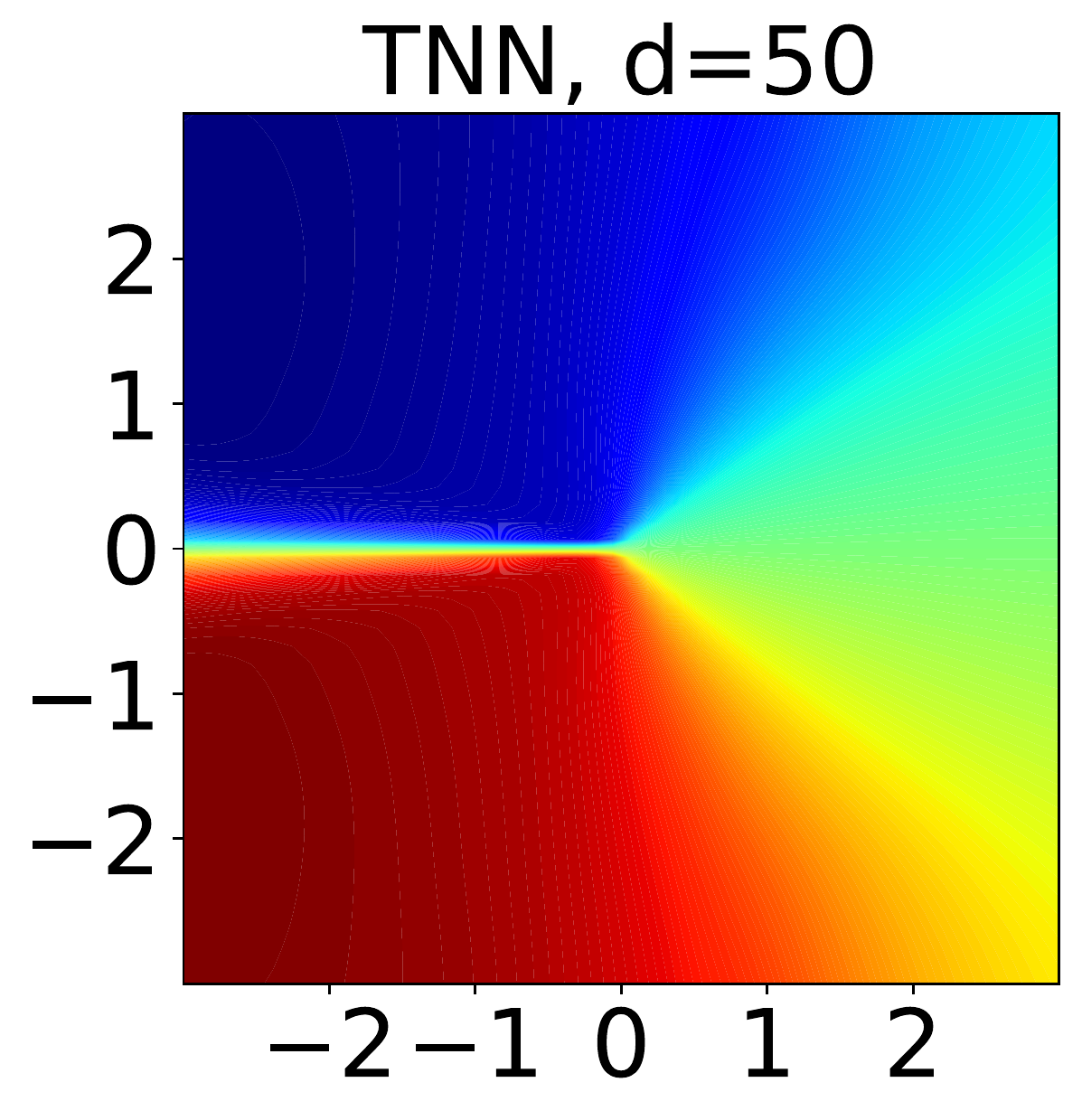}
\includegraphics[width=3cm,height=3cm]{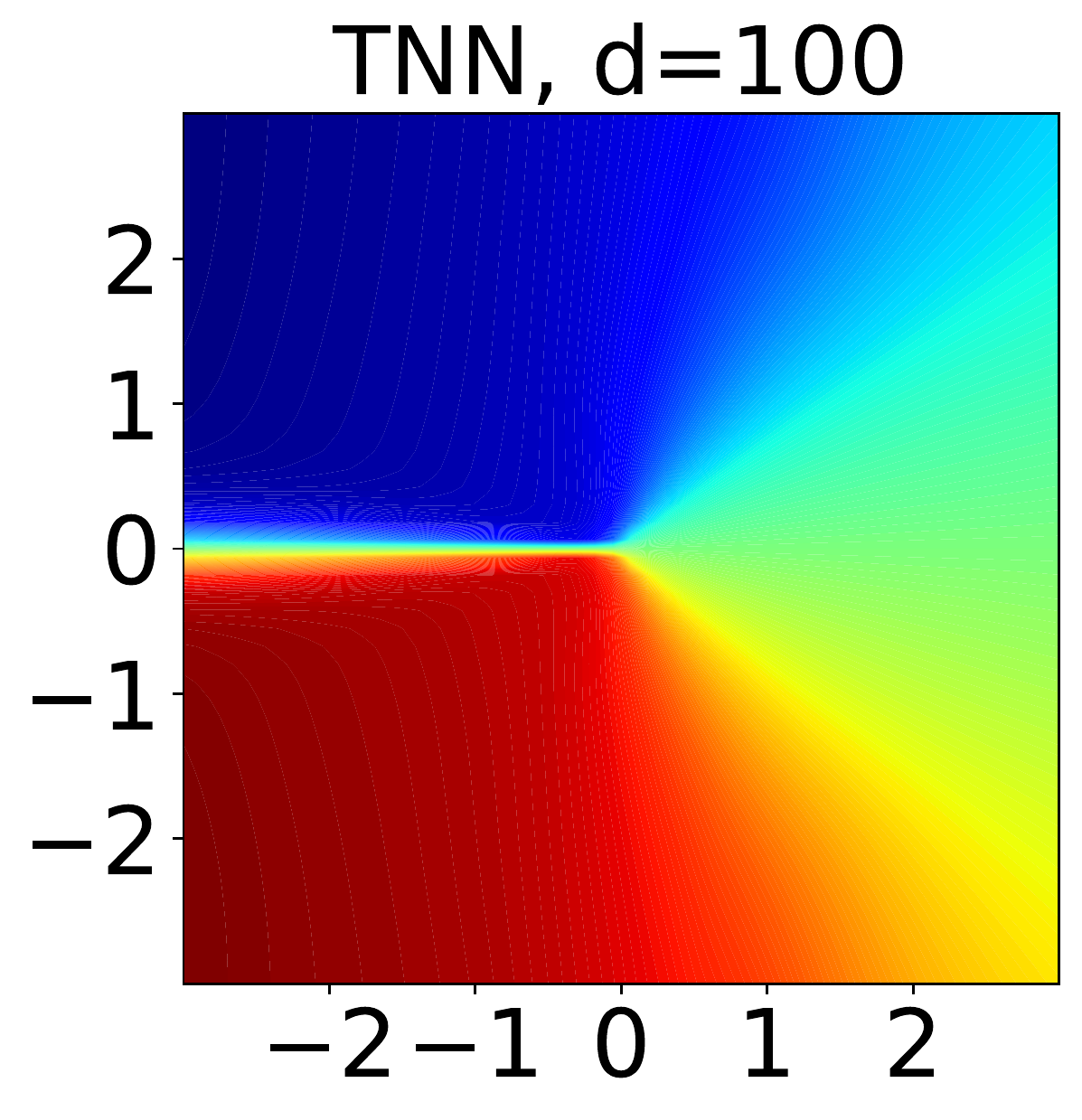}\\
\includegraphics[width=3cm,height=3cm]{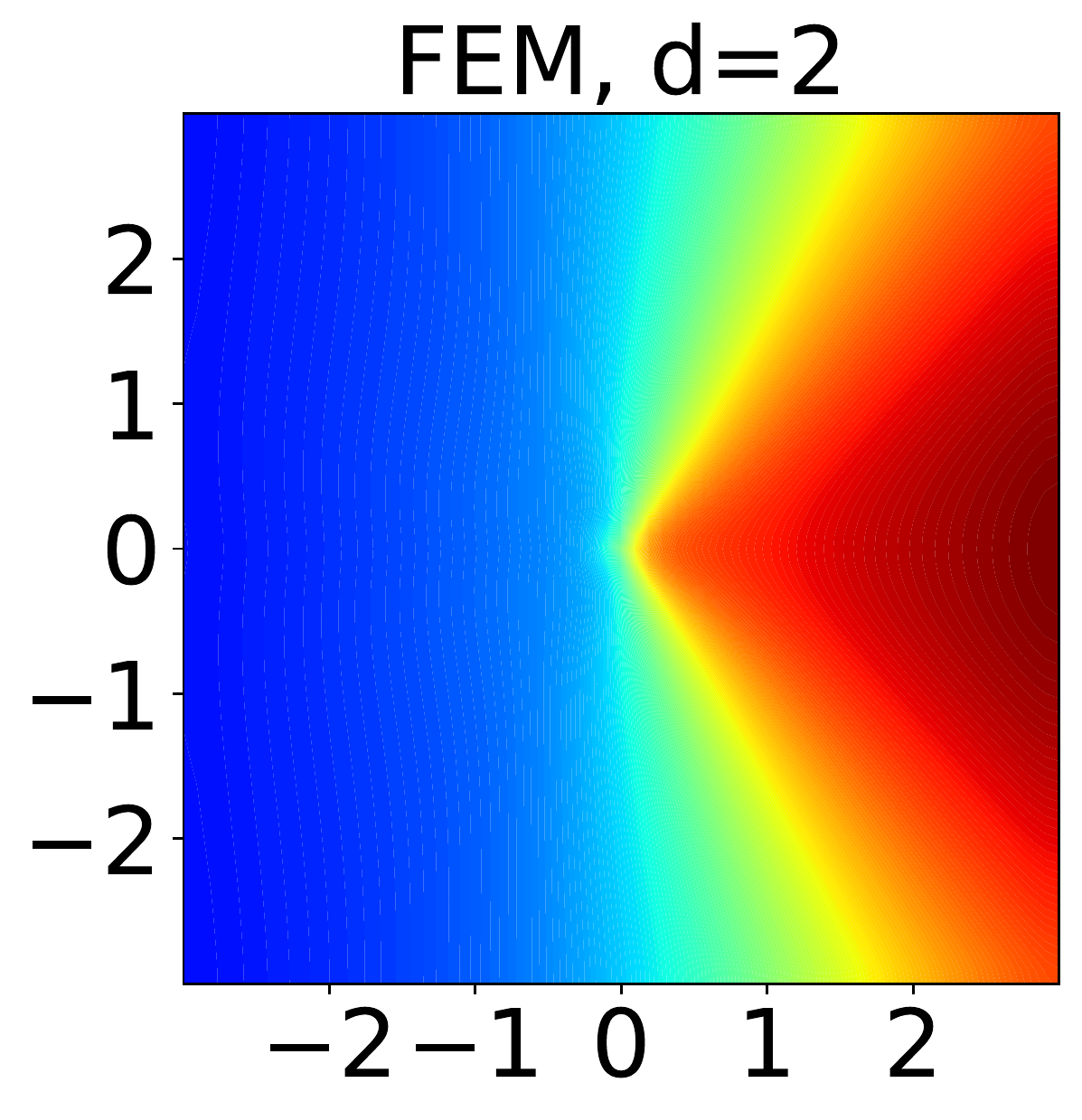}
\includegraphics[width=3cm,height=3cm]{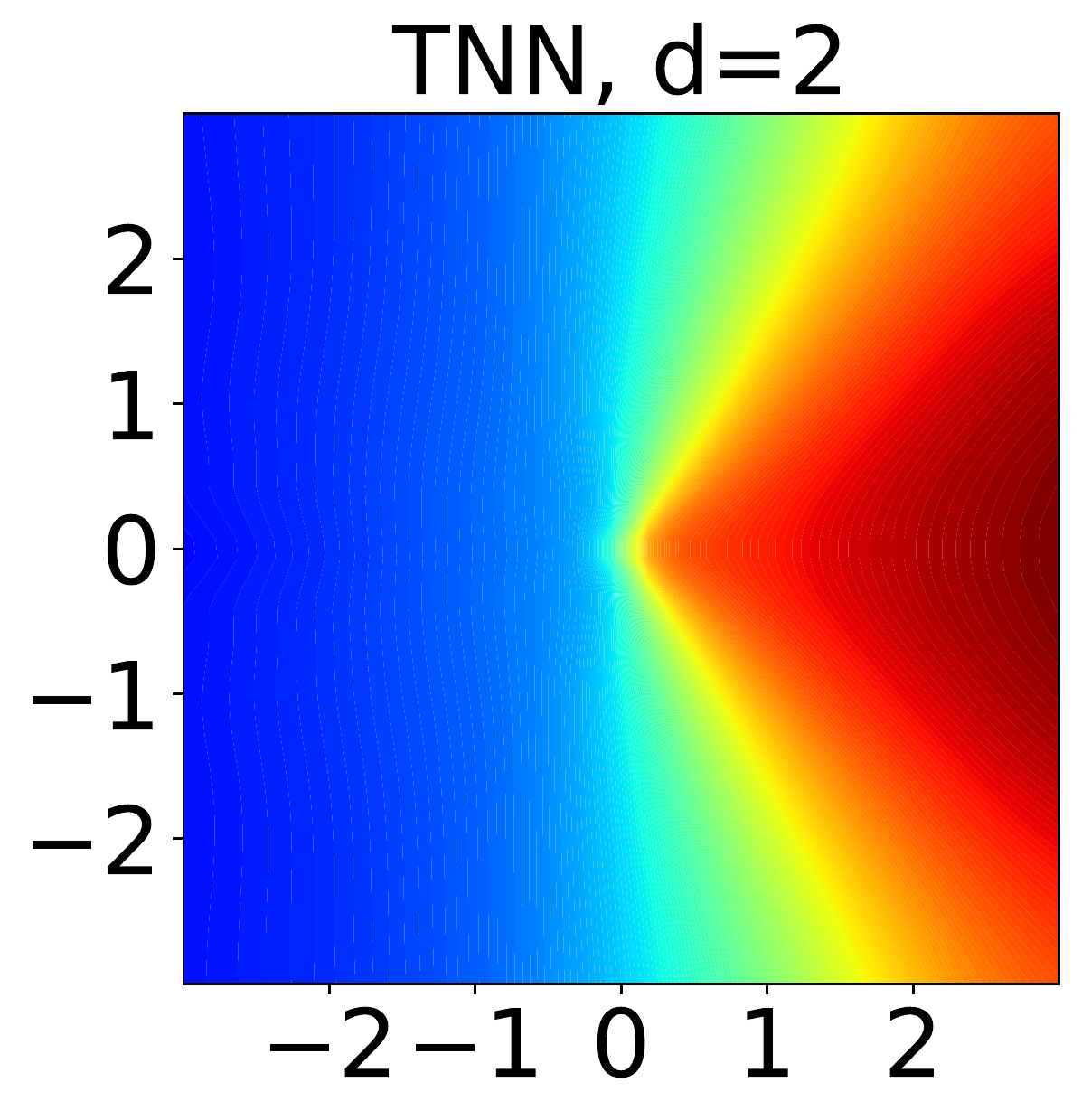}
\includegraphics[width=3cm,height=3cm]{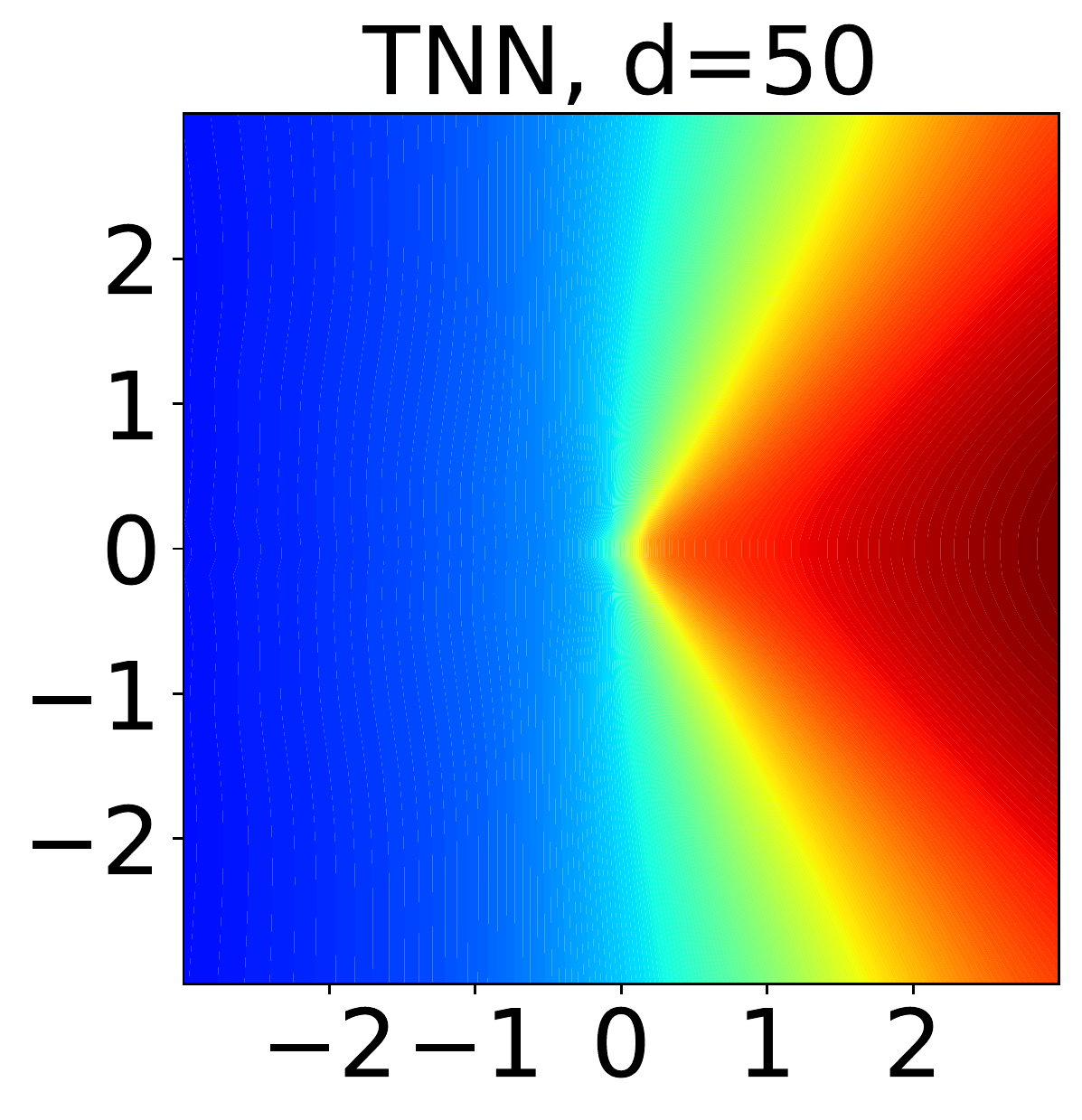}
\includegraphics[width=3cm,height=3cm]{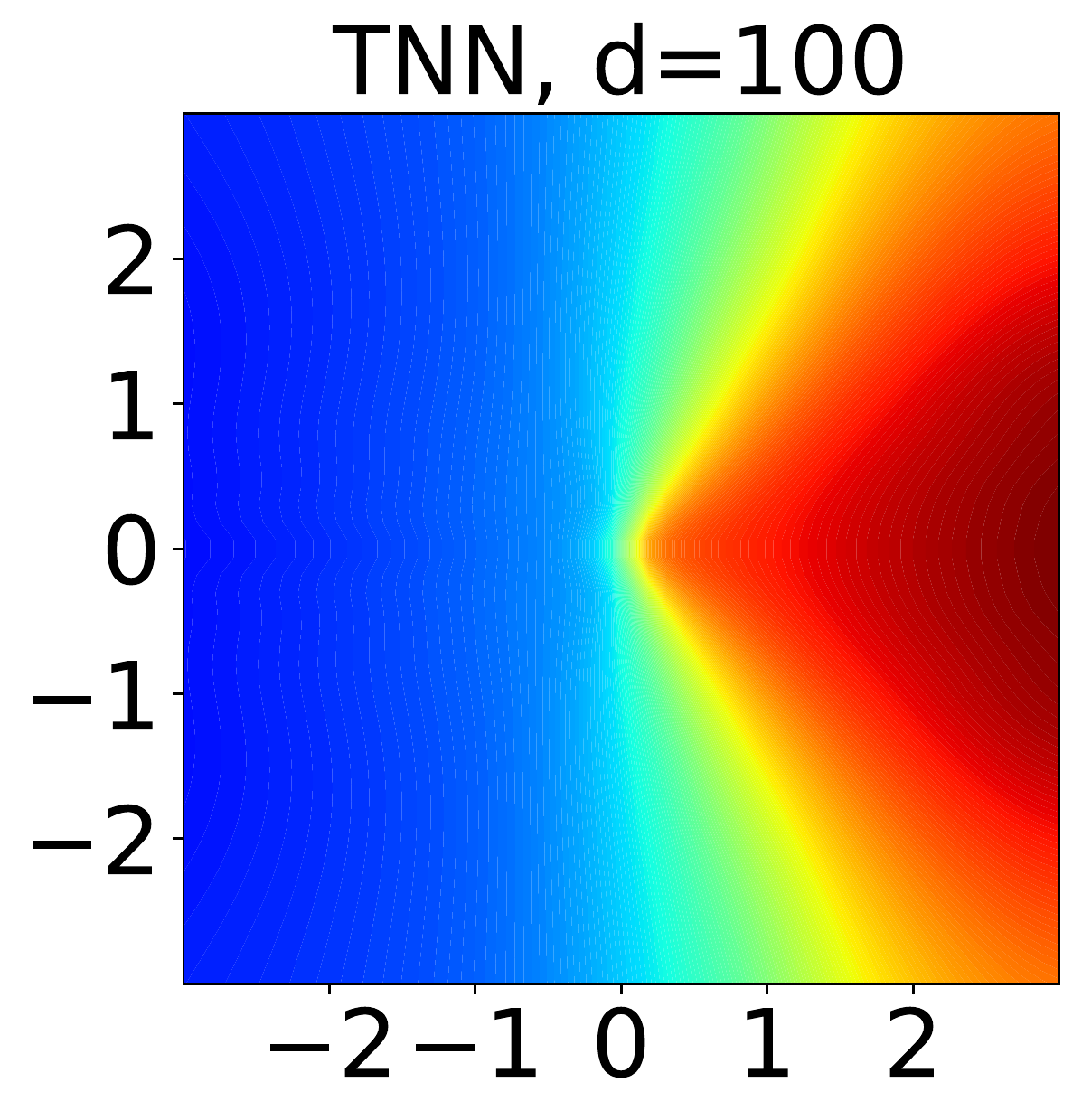}\\
\includegraphics[width=3cm,height=3cm]{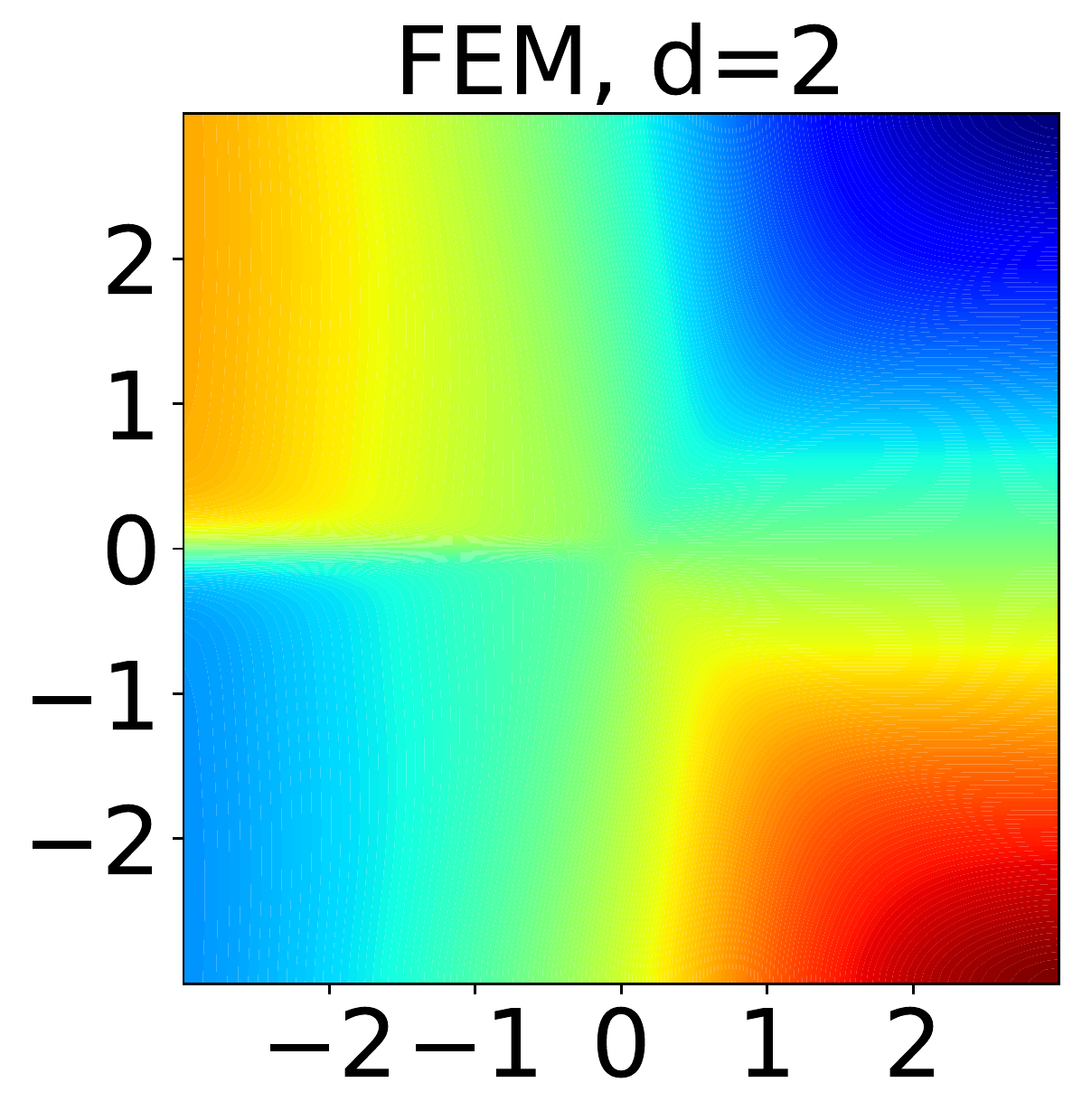}
\includegraphics[width=3cm,height=3cm]{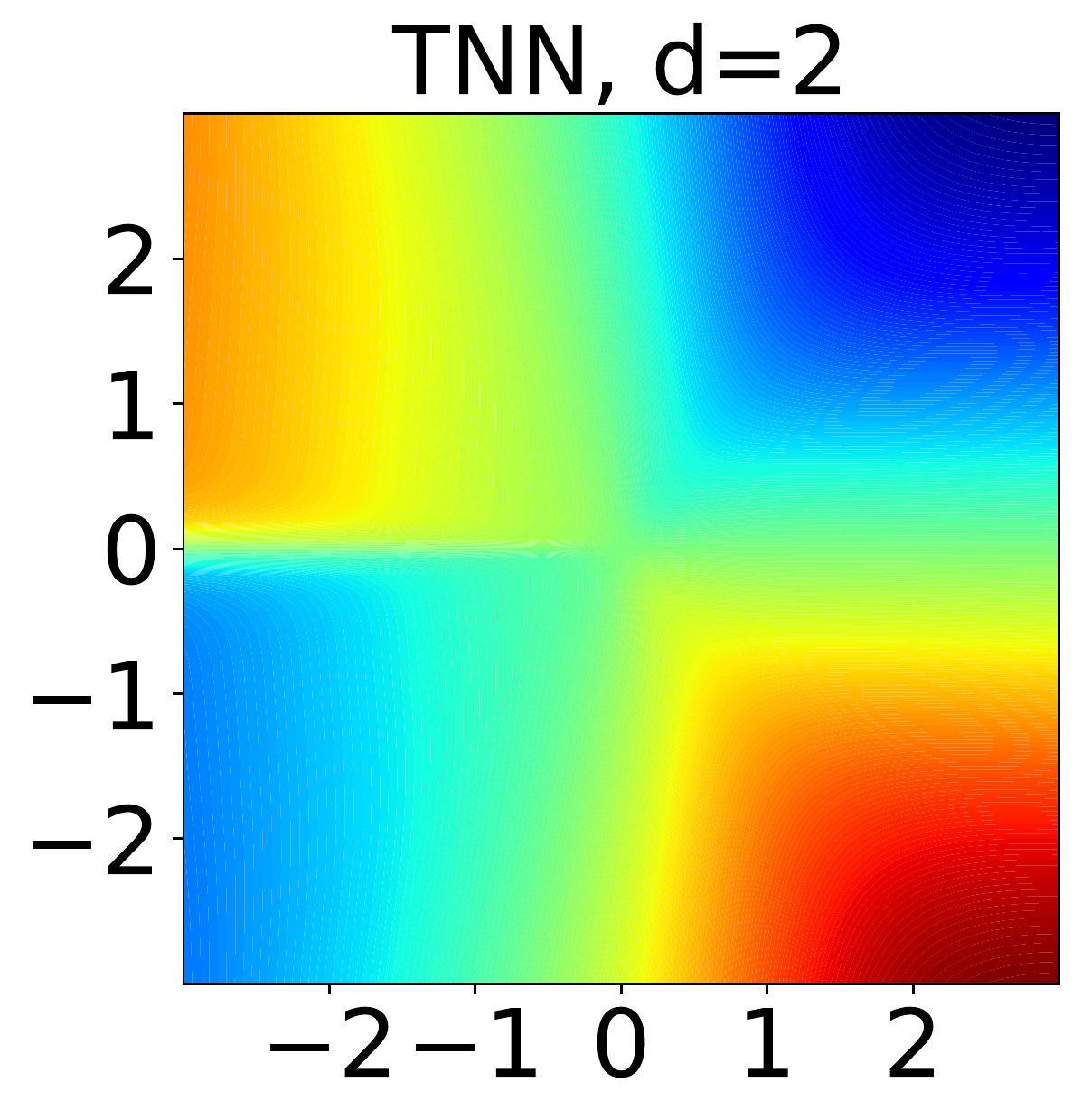}
\includegraphics[width=3cm,height=3cm]{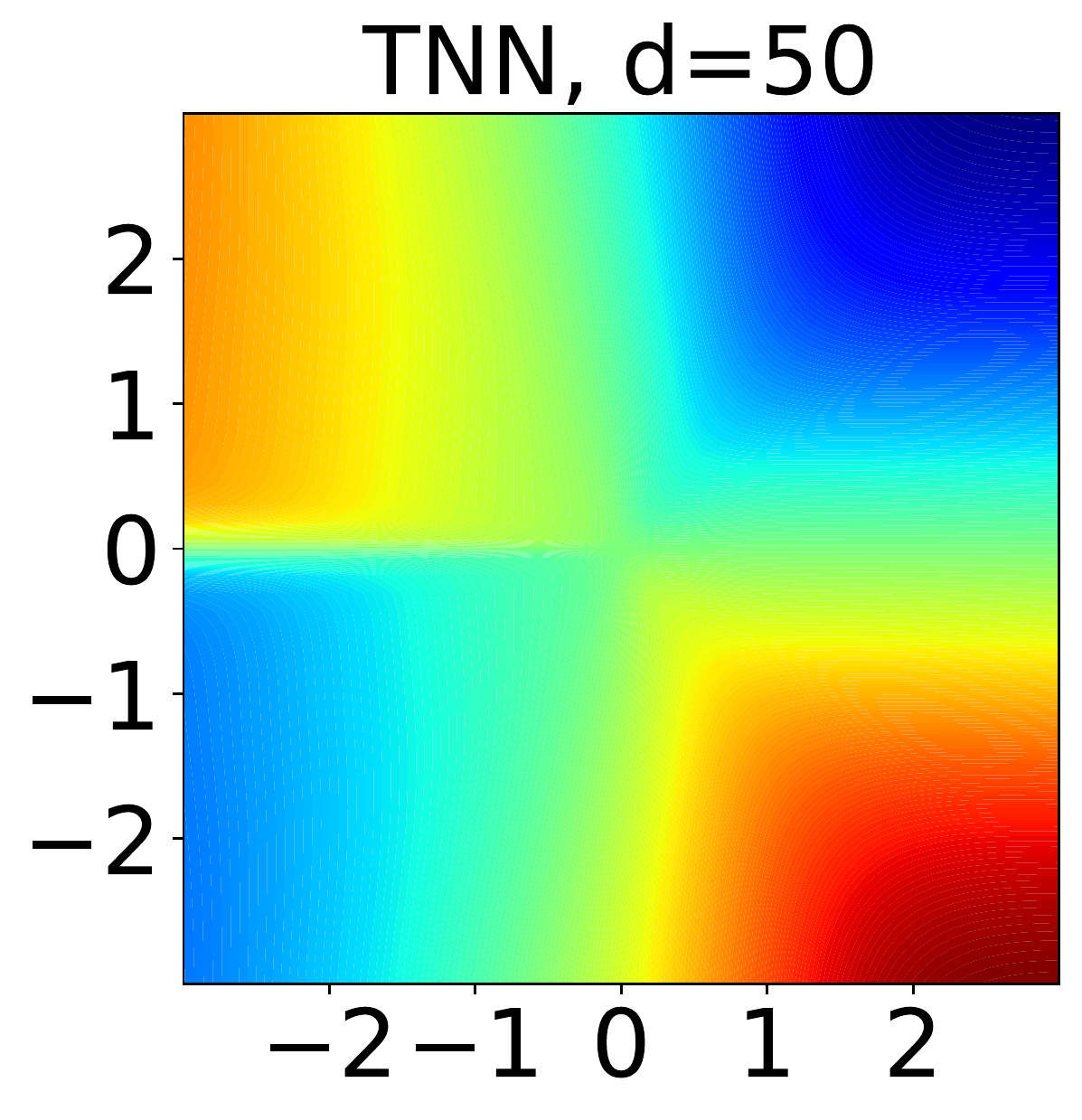}
\includegraphics[width=3cm,height=3cm]{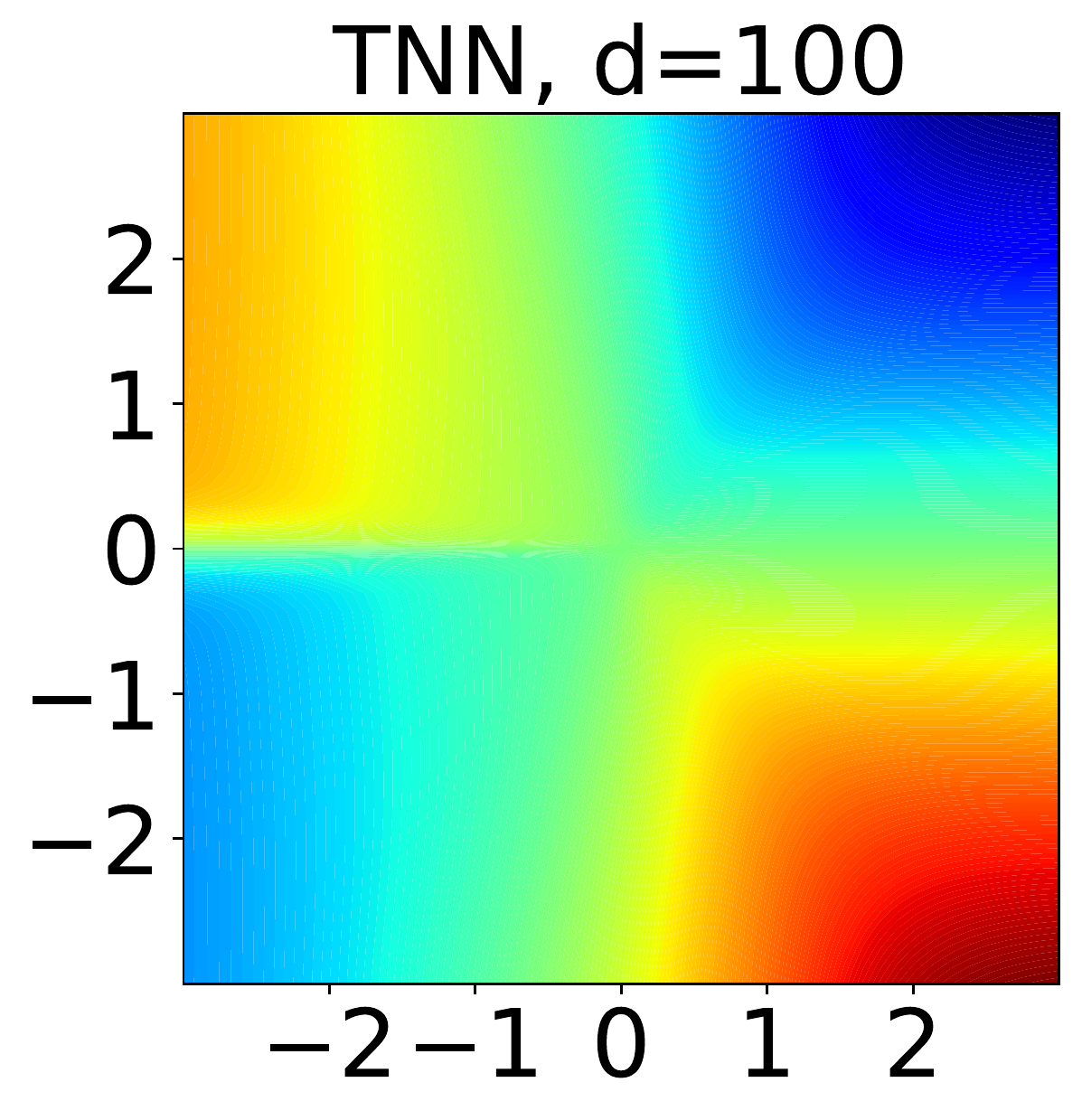}\\
\includegraphics[width=6cm,height=0.5cm]{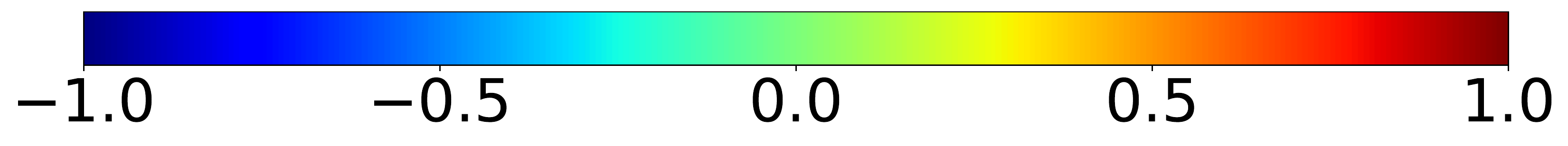}
\caption{The contour plots of the first three eigenfunctions in the first example, computed using the finite element method for $d=2$ (column ``FEM, $d=2$'') and by TNN-based
eigensolver for $d=2$ (column ``TNN, $d=2$''), $d=50$ (column ``TNN, $d=50$''), $d=100$ (column ``TNN, $d=100$''), respectively.}\label{fig_meta}
\end{figure}

The reference eigenvalue for $d=2$ is obtained by using the third order conforming finite element method
with the mesh size $h=\frac{3}{1024}\sqrt{2}$ on the domain $[-3,3]^2\subset \mathbb R^2$.
Here we use the open parallel
finite element package OpenPFEM \cite{OpenPFEM} to do the discretization and then using Krylovschur method from SLEPc \cite{SLEPc} to solve the corresponding algebraic eigenvalue problem \footnote{We express our thanks to Yangfei Liao for this computation}.
In this way, we obtain the first three  eigenvalues
$$\lambda_1 = 0.21881493133369, \ \ \ \lambda_2=0.76371970025476,\ \ \ \lambda_3=2.79019347384363, $$
as the reference values for our numerical investigation.
The corresponding eigenfunctions are shown in the first column of Figure \ref{fig_meta}.

In implementation, we use $3$ TNNs to learn the lowest $3$ eigenvalues. 
For $d=2,50,100$ cases, each TNN has depth 3 and width 20 and the rank is chosen to be $p=10$.
The Adam optimizer is employed with a learning rate 0.003.
We use the Adam optimizer in the first 100000 steps and then the LBFGS in the subsequent 10000 steps.
The final eigenvalue approximations are represented in Table \ref{table_meta} and the corresponding eigenfunction approximations are shown in Figure \ref{fig_meta}, where we can find the proposed method has obviously better accuracy than that in \cite{ZhangLiSchutte}.

\begin{table}[!htb]
\caption{Errors of infinitesimal generators of metastable diffusion processes
problem for the 3 lowest eigenvalues.}\label{table_meta}
\begin{center}
\begin{tabular}{ccccccc}
\hline
&  FEM, $d=2$&   TNN, $d=2$&  TNN, $d=50$& TNN, $d=100$\\
\hline
$\lambda_1$&   0.21881493133369&   0.21882643485835&   0.21883758216194&   0.21884057151253\\
$\lambda_2$&   0.76371970025476&   0.76372599038784&   0.76372707631296&   0.76378255169557\\
$\lambda_3$&   2.79019347384363&   2.79020711711705&   2.79022006384330&   2.79022583685384\\
\hline
\end{tabular}
\end{center}
\end{table}

\subsection{Harmonic oscillator problems}
In this subsection, we consider the Schr\"{o}dinger equation associated with the $d$-dimensional harmonic oscillator
and the corresponding Hamiltonian operator reads as
\begin{eqnarray}\label{eq_HO_x}
H=-\frac{1}{2}\sum_{i=1}^d\nabla_i^2+\frac{1}{2}x^TAx,
\end{eqnarray}
where $x=(x_1,x_2,\cdots,x_d)^T$, $A=(a_{ij})_{d\times d}\in\mathbb R^{d\times d}$ is a symmetric
positive definite matrix.

The exact wavefunctions and the corresponding energy (i.e. eigenvalue) of different states can be obtained with
the similar way in \cite{CHO}.
According to the property of the symmetric positive definite matrix, there exists an orthogonal matrix
$Q=(q_{ij})_{d\times d}\in\mathbb R^{d\times d}$ such that
\begin{eqnarray}\label{eq_SVD}
Q^TAQ={\rm diag}\big\{\mu_1,\mu_2,\cdots,\mu_d\big\},
\end{eqnarray}
where $\mu_j$ and $\mathbf q_j=(q_{1j},q_{2j},\cdots,q_{dj})^T$, $j=1,2,\cdots,d$ are eigenvalues and
corresponding normalized eigenvectors of the matrix $A$.
Then under the rotation transformation $y=Q^Tx$, the Hamiltonian operator (\ref{eq_HO_x})
in the coordinate system  $(y_1,y_2,\cdots,y_d)$ can be written as the following decoupled harmonic form
\begin{eqnarray}
H=-\frac{1}{2}\sum_{i=1}^d\nabla_i^2+\frac{1}{2}\sum_{i=1}^d\mu_iy_i^2.
\end{eqnarray}
Then the exact wavefunction of state $(n_1,n_2,\cdots,n_d)$ can be immediately obtained as follows
\begin{eqnarray}
\Psi_{n_1,n_2,\cdots,n_d}(y_1,y_2,\cdots,y_d)=\prod_{i=1}^d\mathcal H_{n_i}(\mu_i^{1/4}y_i)e^{-\mu_i^{1/2}y_i^2/2},
\end{eqnarray}
and the corresponding exact energy is
\begin{eqnarray}
E_{n_1,n_2,\cdots,n_d}=\sum_{i=1}^d\Big(\frac{1}{2}+n_i\Big)\mu_i^{1/2},
\end{eqnarray}
where $y_i=\sum_{j=1}^dq_{ij}x_j$ and $\mathcal H_n$ is physicists' Hermite polynomials \cite{AtakishievSuslov}.

In all three examples of this subsection, we adopt the same parameters $a_{ij}$ as in \cite{LiZhaiChen}
and calculate the lowest $16$ energy states as \cite{LiZhaiChen} done. The aim here is to
demonstrate the performance of the proposed approach through these comparisons.
%Figure \ref{Figure_LiZhaiChen} shows the corresponding numerical results.  From Figure \ref{Figure_LiZhaiChen},
%we can find the proposed numerical method in this paper can obtain better accuracy than that in \cite{LiZhaiChen}.

\subsubsection{Two-dimensional harmonic oscillator}
We first examine a simple case of the two-dimensional harmonic oscillator problem
\begin{eqnarray}
-\frac{1}{2}\Big(\nabla_1^2+\nabla_2^2\Big)\Psi+\frac{1}{2}(x_1^2+x_2^2)\Psi=E\Psi.
\end{eqnarray}
Since the problem is essentially decoupled, the exact eigenfunction has low-rank representation
in coordinate $(x_1,x_2)$ as follows
\begin{eqnarray}
\Psi_{n_1,n_2}=\mathcal H_{n_1}(x_1)e^{-x_1^2/2}\mathcal H_{n_2}(x_2)e^{-x_2^2/2}.
\end{eqnarray}
The corresponding exact energy is
\begin{eqnarray}
E_{n_1,n_2}=\Big(\frac{1}{2}+n_1\Big)+\Big(\frac{1}{2}+n_2\Big).
\end{eqnarray}
We use $16$ TNNs to learn the lowest $16$ energy states
and each subnetwork of a single TNN with depth $3$ and width $50$, and $p=20$.
The Adam optimizer is employed with a learning rate 0.001 in the first 500000 epochs 
and then the L-BFGS in the subsequent 10000 steps to produce the final result.
The Hermite-Gauss quadrature scheme with $99$ points  are adopted in each dimension. 

The corresponding numerical results are shown in Table \ref{table_2dHO}, where we can find the
proposed TNN-based machine learning method has obvious better accuracy than that in \cite{LiZhaiChen},
where the accuracy is $1.0${\rm e}-$2$ by the Monte-Carlo based machine learning methods.
\begin{table}[!htb]
\caption{Errors of two-dimensional harmonic oscillator problem for the 16 lowest energy states.}\label{table_2dHO}
\begin{center}
\begin{tabular}{ccccccc}
\hline
$n$&  $(n_1,n_2)$&   Exact $E_n$&  Approx $E_n$&  ${\rm err}_E$&  ${\rm err}_{L^2}$&   ${\rm err}_{H^1}$\\
\hline
0&   (0,0)&   1.0&   1.000000000000441&   4.414e-13&   2.935e-07&    7.314e-07\\
1&   (0,1)&   2.0&   2.000000000138887&   6.944e-11&   5.889e-06&    1.021e-05\\
2&   (1,0)&   2.0&   2.000000000369235&   1.846e-10&   9.371e-06&    1.550e-05\\
3&   (0,2)&   3.0&   3.000000000851601&   2.839e-10&   1.524e-05&    2.180e-05\\
4&   (1,1)&   3.0&   3.000000001492529&   4.975e-10&   1.970e-05&    2.958e-05\\
5&   (2,0)&   3.0&   3.000000005731970&   1.911e-09&   3.964e-05&    5.852e-05\\
6&   (0,3)&   4.0&   4.000000000287978&   7.200e-11&   8.690e-06&    1.312e-05\\
7&   (1,2)&   4.0&   4.000000000727938&   1.820e-10&   1.374e-05&    1.838e-05\\
8&   (2,1)&   4.0&   4.000000001748148&   4.370e-10&   2.272e-05&    3.100e-05\\
9&   (3,0)&   4.0&   4.000000005090556&   1.273e-09&   3.645e-05&    5.110e-05\\
10&  (1,3)&   5.0&   5.000000000117699&   2.354e-11&   5.199e-06&    2.009e-05\\
11&  (2,2)&   5.0&   5.000000000746078&   1.492e-10&   1.821e-05&    3.091e-05\\
12&  (3,1)&   5.0&   5.000000001093248&   2.186e-10&   1.973e-05&    3.467e-05\\
13&  (0,4)&   5.0&   5.000000001562438&   3.125e-10&   2.651e-05&    2.520e-05\\
14&  (4,0)&   5.0&   5.000000004861336&   9.723e-10&   4.059e-05&    4.161e-05\\
15&  (3,2)&   6.0&   6.000000043151862&   7.192e-09&   3.095e-05&    3.760e-05\\
\hline
\end{tabular}
\end{center}
\end{table}

\subsubsection{Two-dimensional coupled harmonic oscillator}
Since TNN-based methods carry out operations separately in each dimension, it is not
unexpected that our method has an impressive performance in the completely decoupled case.
To show the generalities of TNN-based machine learning method,
the next two examples are to compute the multi-states of operators with coupled oscillators.

First, we consider the following two-dimensional eigenvalue problem of the operator with
coupled harmonic oscillator
\begin{eqnarray}
-\frac{1}{2}\Big(\nabla_1^2+\nabla_2^2\Big)\Psi+\frac{1}{2}(a_{11}x_1^2+2a_{12}x_1x_2+a_{22}x_2^2)\Psi=E\Psi.
\end{eqnarray}
The coefficients $a_{11},a_{12},a_{22}$ are chosen as in \cite{LiZhaiChen}.
Since in \cite{LiZhaiChen} the coefficients are rounded to four decimal places,
we use rounded values as exact coefficients, that is, $a_{11}=0.8851$, $a_{12}=-0.1382$,
$a_{22}=1.1933$. The exact eigenfunction has the following representation in coordinate $(y_1,y_2)$
\begin{eqnarray}
\Psi_{n_1,n_2}(y_1,y_2)=\mathcal H_{n_1}(\mu_1^{1/4}y_1)e^{-\mu_1^{1/2}y_1^2/2}\cdot\mathcal H_{n_2}(\mu_2^{1/4}y_2)e^{-\mu_2^{1/2}y_2^2/2},
\end{eqnarray}
where $y_1=-0.9339352418x_1-0.3574422527x_2$ and $y_2=0.3574422527x_1-0.9339352418$, $\mu_1=0.8322071257$,
$\mu_2=1.2461928742$, the two quantum numbers $n_1,n_2$ take the values $0,1,2,\cdots$.
The exact energy for the state $(n_1,n_2)$ is
\begin{eqnarray}
E_{n_1,n_2}=\Big(\frac{1}{2}+n_1\Big)\mu_1^{1/2}+\Big(\frac{1}{2}+n_2\Big)\mu_2^{1/2}.
\end{eqnarray}

We use $16$ TNNs to learn the lowest $16$ energy states.
In each TNN, the rank is chosen to be $p=20$, the subnetwork is built with depth $3$ and width $50$.
The Adam optimizer is employed with a learning rate 0.001.
We use the Adam optimizer in the first 500000 steps and then the L-BFGS in the subsequent
10000 steps. In each direction, $99$ points Hermite-Gauss quadrature scheme are used to do the integration.

The corresponding numerical results are collected in Table \ref{table_2dCHO}, where we can find the
 proposed numerical method can also obtain the higher accuracy for the Schr\"{o}dinger equation with
 coupled  harmonic oscillator.
 \begin{table}[!htb]
\caption{Errors of two-dimensional coupled harmonic oscillator problem for the 16 lowest energy states.}\label{table_2dCHO}
\begin{center}
\begin{tabular}{ccccccc}
\hline
$n$&  $(n_1,n_2)$&   Exact $E_n$&  Approx $E_n$&  ${\rm err}_E$  &${\rm err}_{L^2}$&   ${\rm err}_{H^1}$\\
\hline
0&   (0,0)&   1.014291981649766&   1.014291988589516&   6.842e-09&    2.801e-05&  9.650e-05\\
1&   (1,0)&   1.926545852963290&   1.926545854461407&   7.776e-10&    1.264e-05&  3.167e-05\\
2&   (0,1)&   2.130622073635773&   2.130622076362180&   1.280e-09&    1.766e-05&  4.081e-05\\
3&   (2,0)&   2.838799724276814&   2.838799728095222&   1.345e-09&    2.099e-05&  4.402e-05\\
4&   (1,1)&   3.042875944949297&   3.042875947694890&   9.023e-10&    1.633e-05&  3.673e-05\\
5&   (0,2)&   3.246952165621781&   3.246952166999784&   4.244e-10&    1.171e-05&  2.485e-05\\
6&   (3,0)&   3.751053595590338&   3.751053597870394&   6.078e-10&    1.500e-05&  3.116e-05\\
7&   (2,1)&   3.955129816262821&   3.955129818022521&   4.449e-10&    1.289e-05&  2.784e-05\\
8&   (1,2)&   4.159206036935306&   4.159206038606588&   4.018e-10&    1.281e-05&  2.388e-05\\
9&   (0,3)&   4.363282257607788&   4.363282258514639&   2.078e-10&    9.008e-06&  2.002e-05\\
10&  (4,0)&   4.663307466903863&   4.663307470243584&   7.162e-10&    1.856e-05&  3.314e-05\\
11&  (3,1)&   4.867383687576346&   4.867383691191087&   7.426e-10&    1.987e-05&  3.338e-05\\
12&  (2,2)&   5.071459908248830&   5.071459911990555&   7.378e-10&    1.992e-05&  3.220e-05\\
13&  (1,3)&   5.275536128921312&   5.275536131659159&   5.190e-10&    1.741e-05&  2.794e-05\\
14&  (0,4)&   5.479612349593796&   5.479612351630867&   3.718e-10&    1.421e-05&  2.538e-05\\
15&  (5,0)&   5.575561338217387&   5.575561344662695&   1.156e-09&    2.769e-05&  3.627e-05\\
\hline
\end{tabular}
\end{center}
\end{table}

\begin{figure}[htb!]
\centering
\includegraphics[width=2.5cm,height=2.5cm]{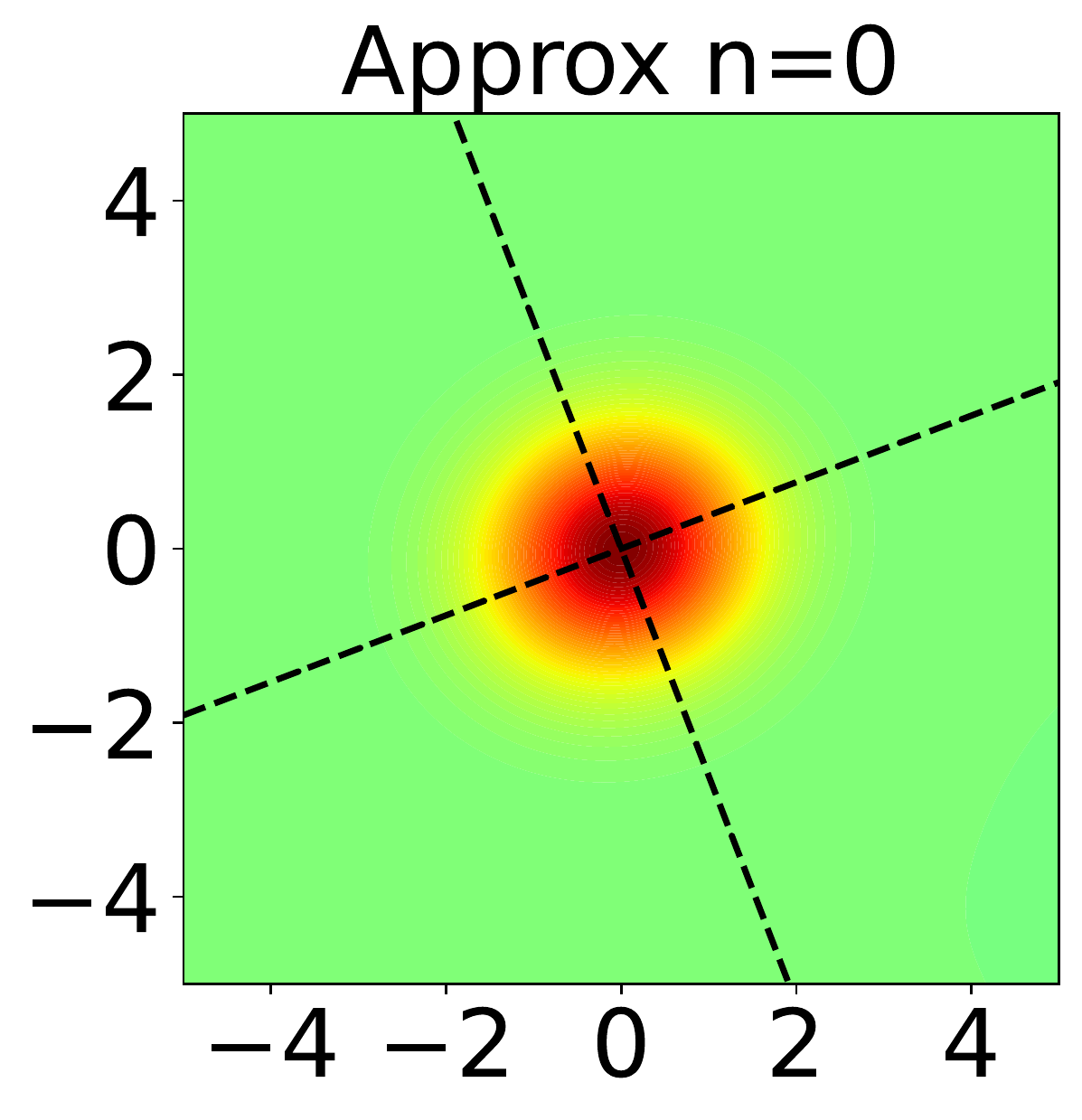}
\includegraphics[width=2.5cm,height=2.5cm]{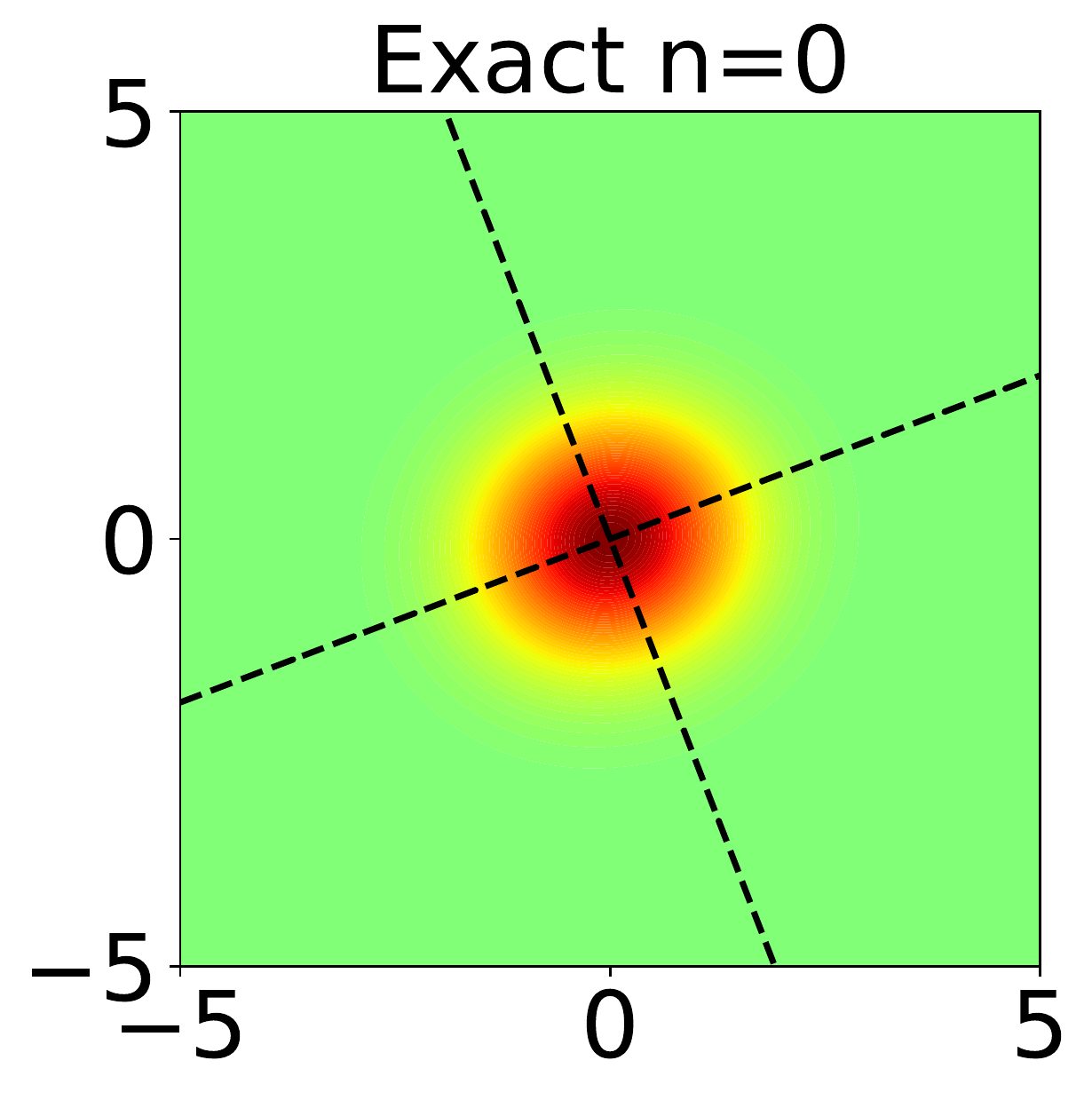}
\includegraphics[width=2.5cm,height=2.5cm]{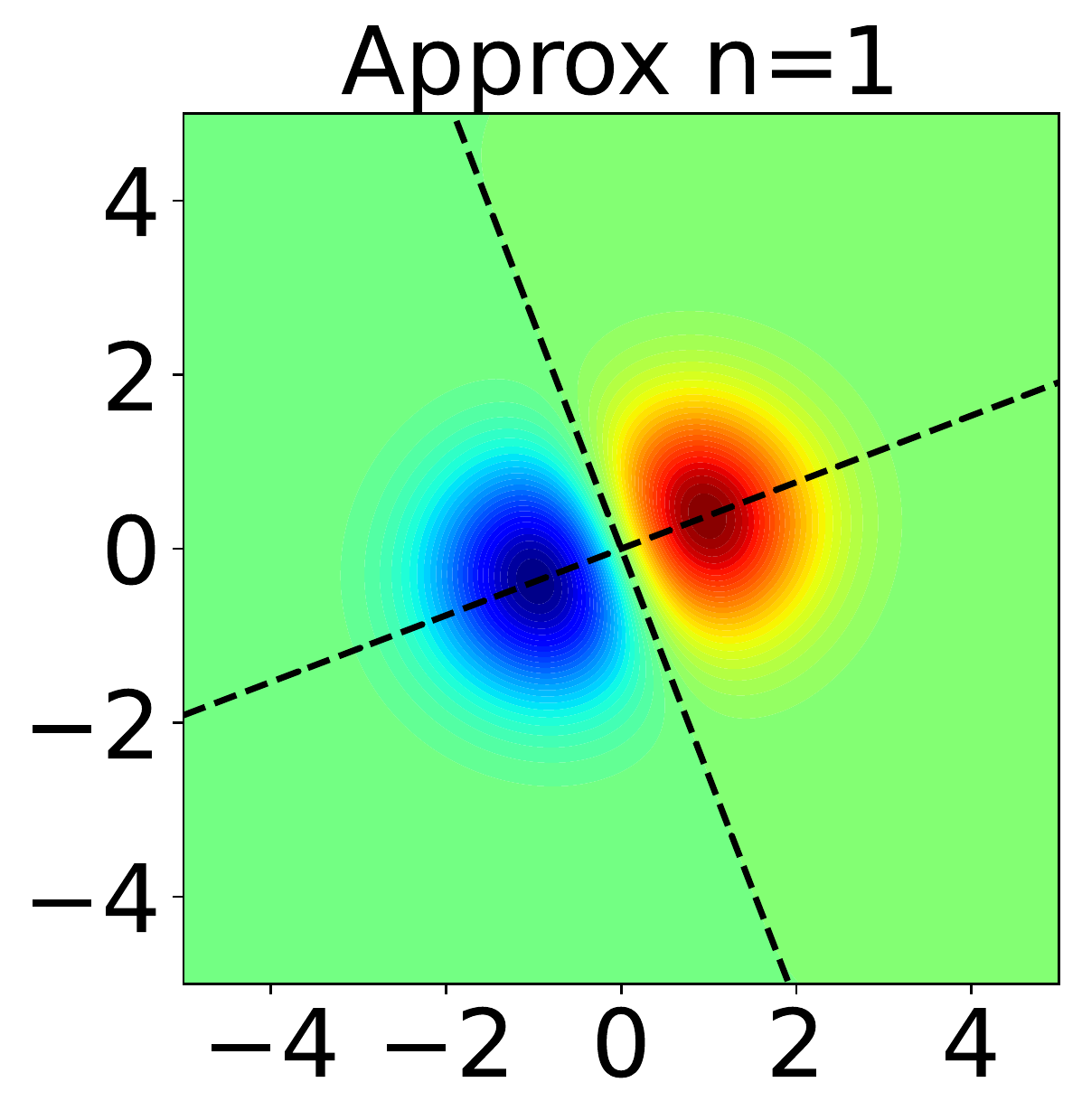}
\includegraphics[width=2.5cm,height=2.5cm]{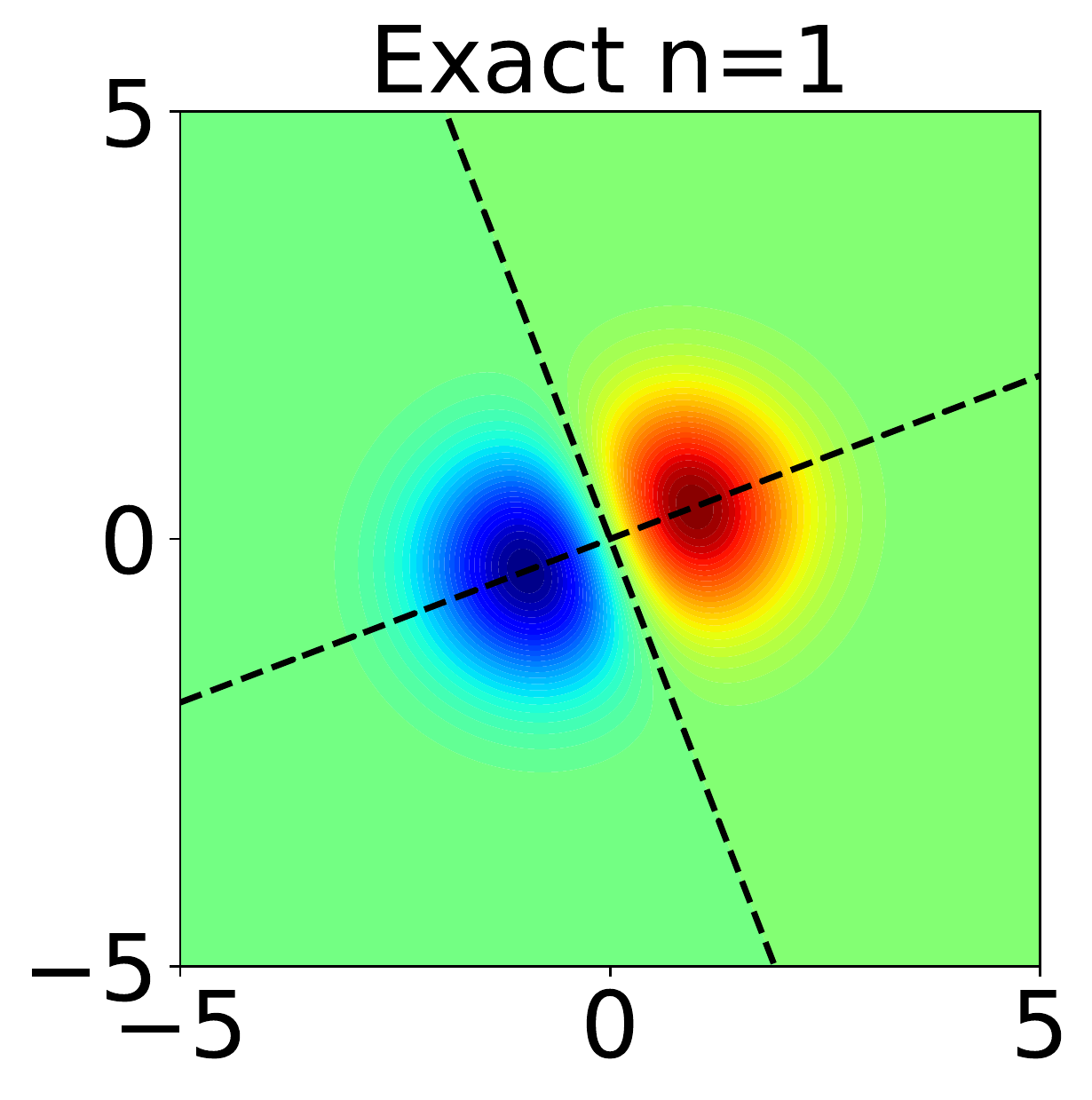}
\includegraphics[width=2.5cm,height=2.5cm]{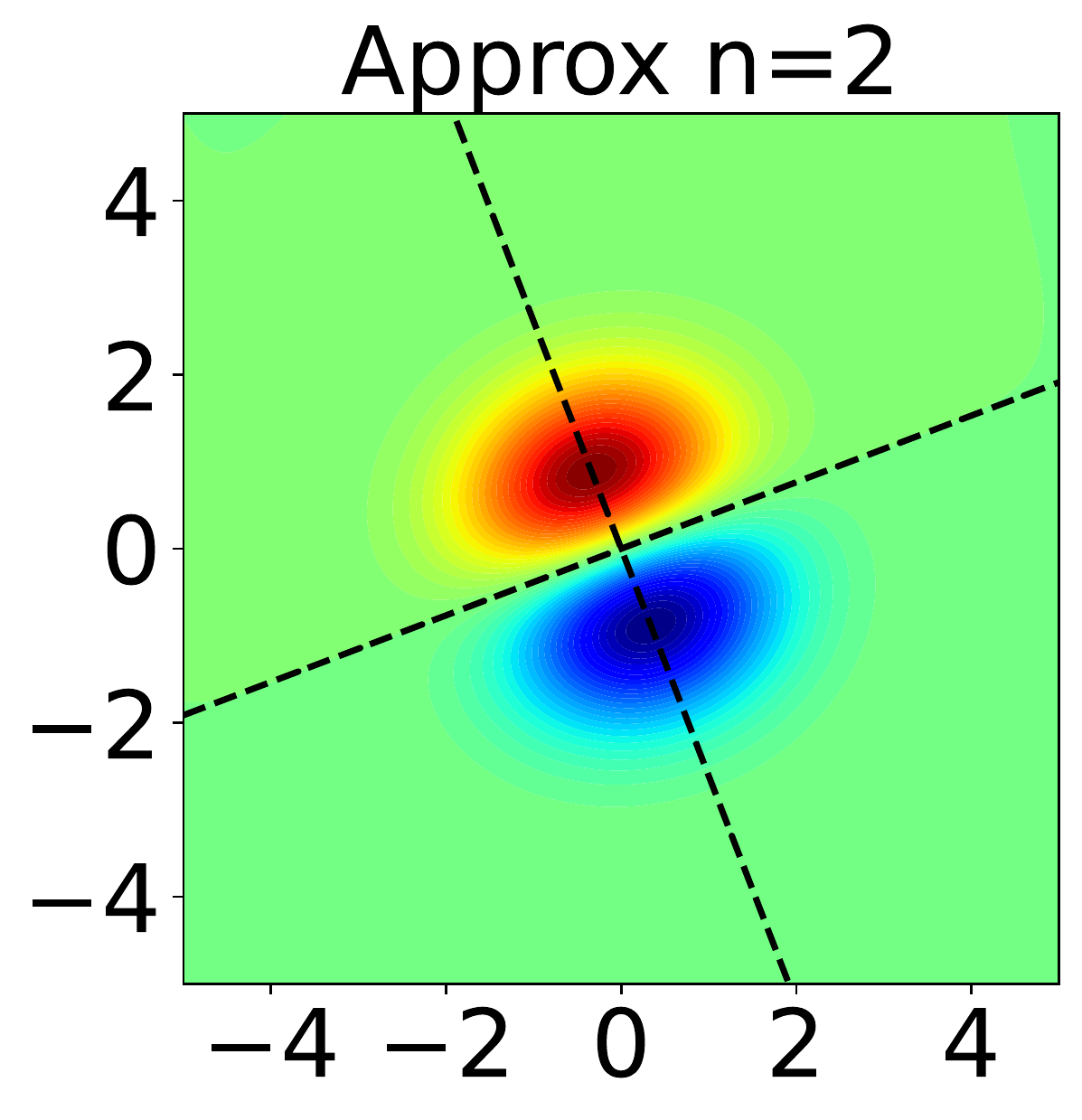}
\includegraphics[width=2.5cm,height=2.5cm]{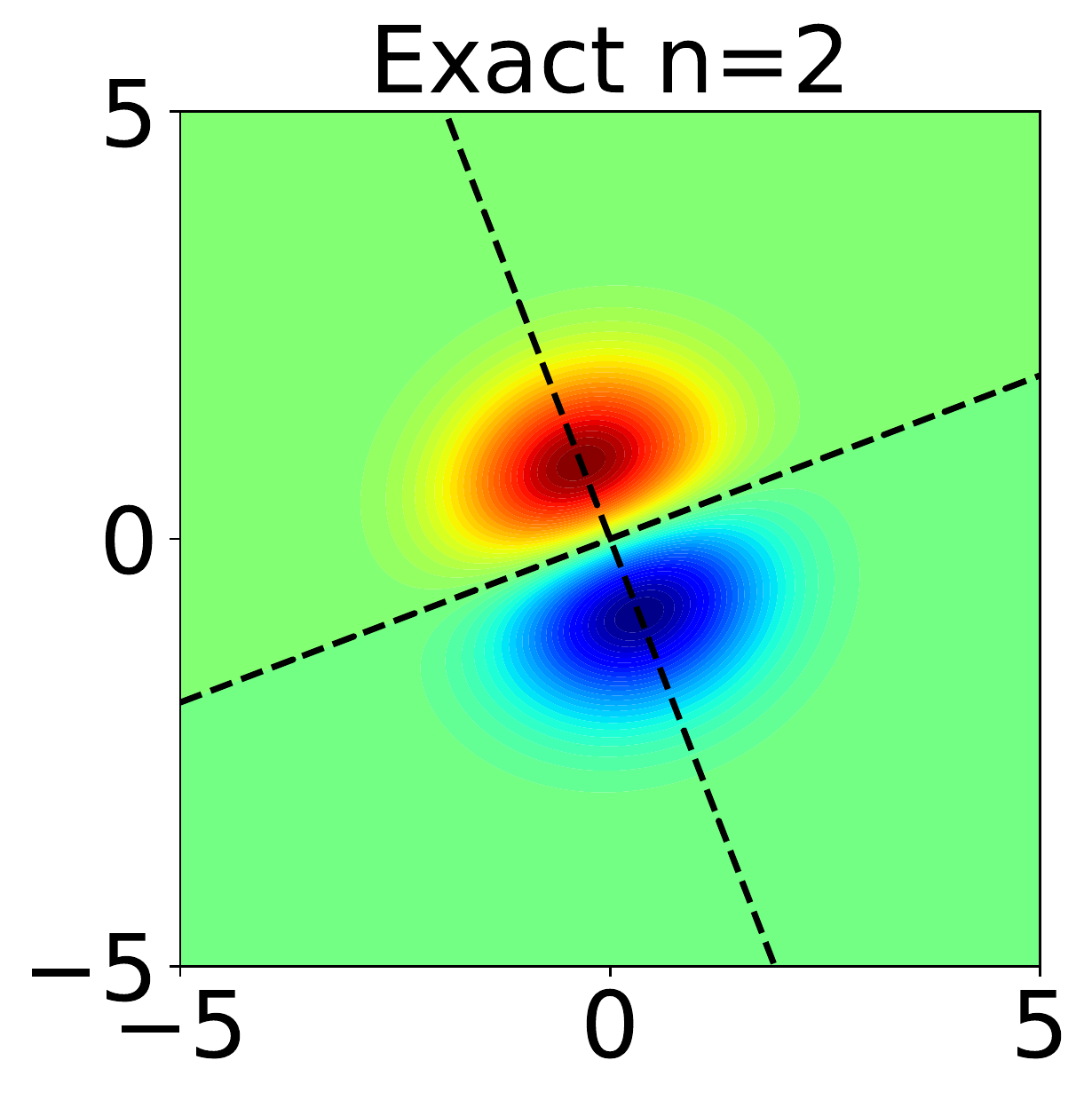}\\
\includegraphics[width=2.5cm,height=2.5cm]{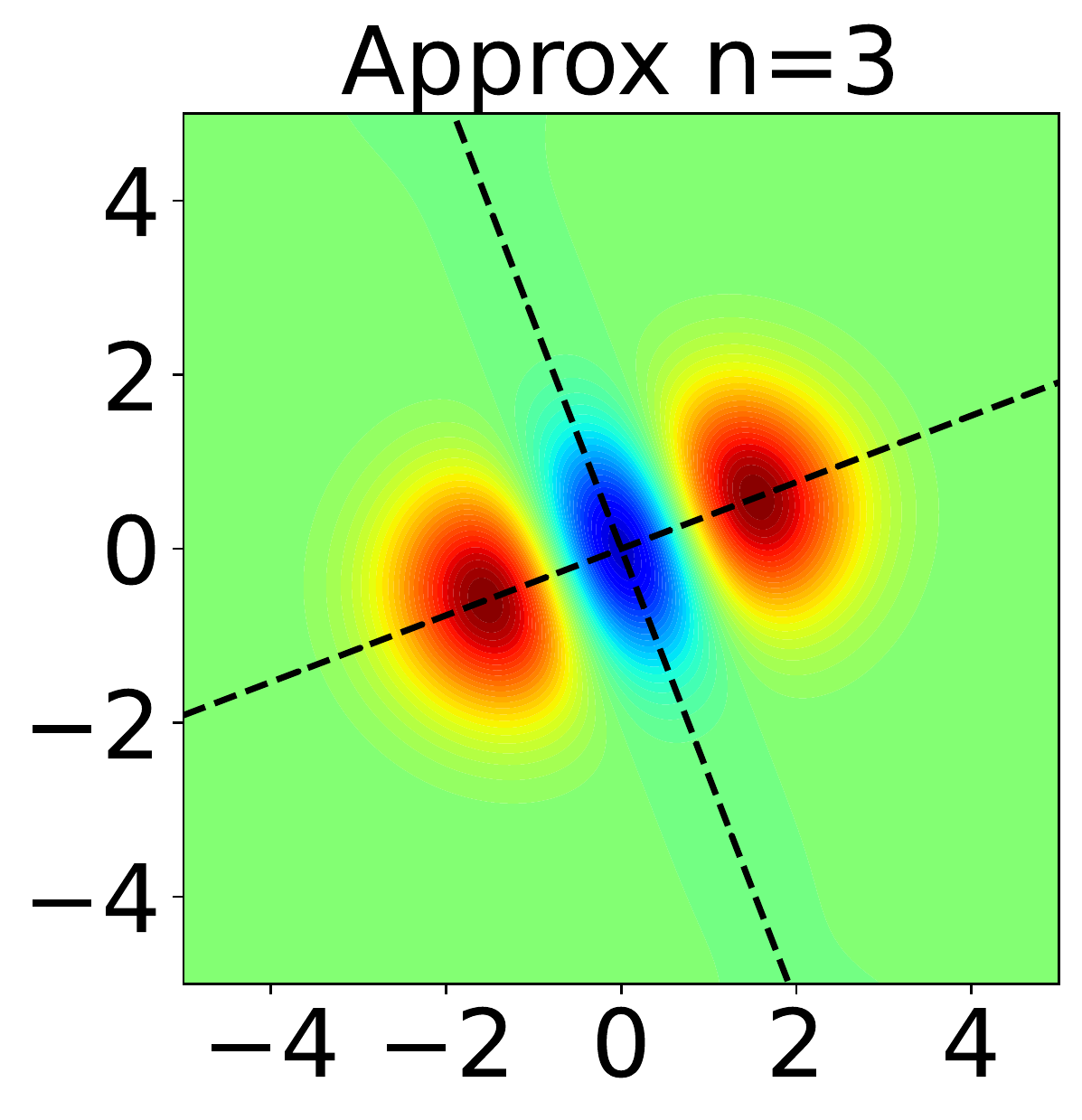}
\includegraphics[width=2.5cm,height=2.5cm]{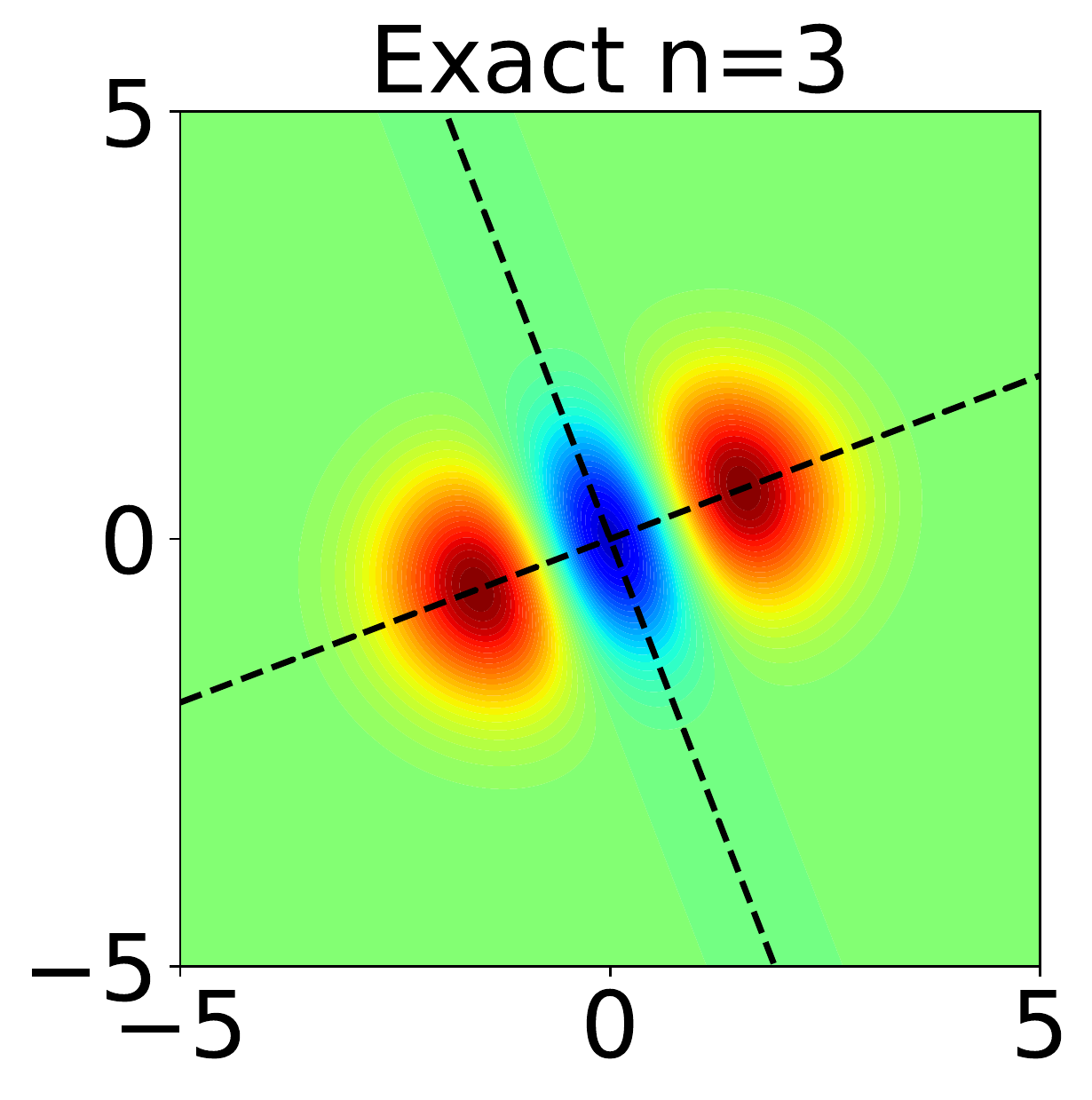}
\includegraphics[width=2.5cm,height=2.5cm]{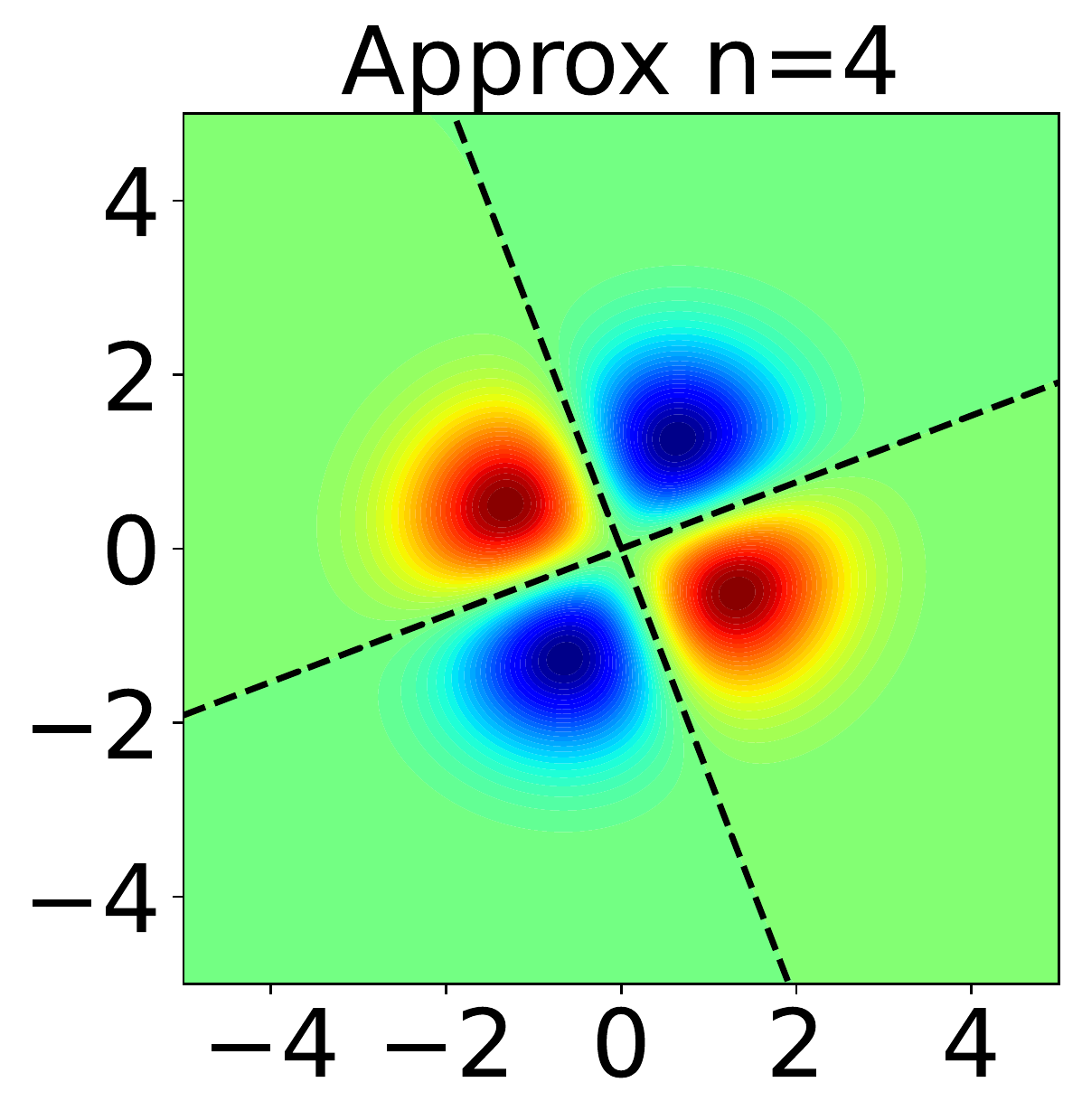}
\includegraphics[width=2.5cm,height=2.5cm]{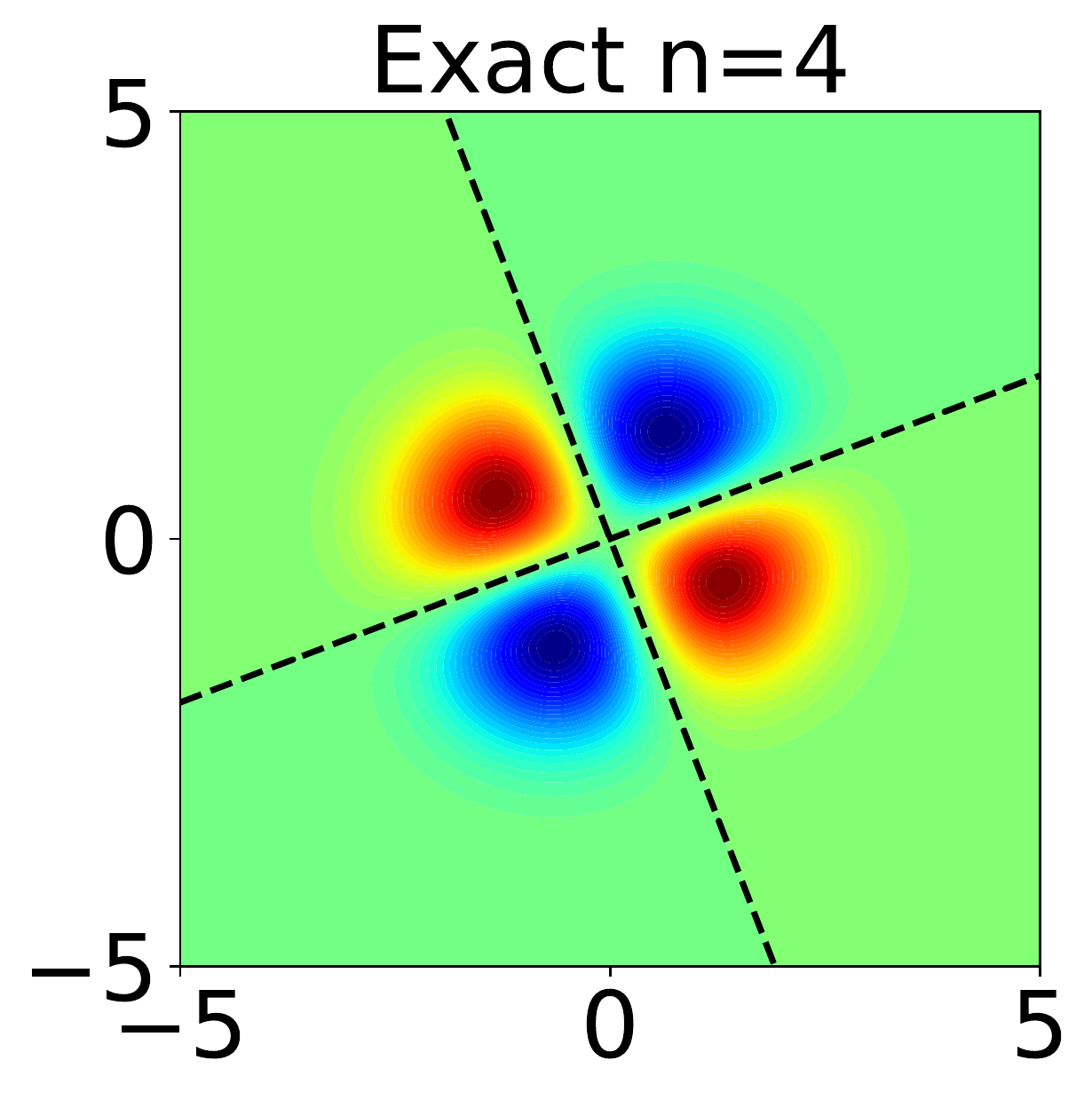}
\includegraphics[width=2.5cm,height=2.5cm]{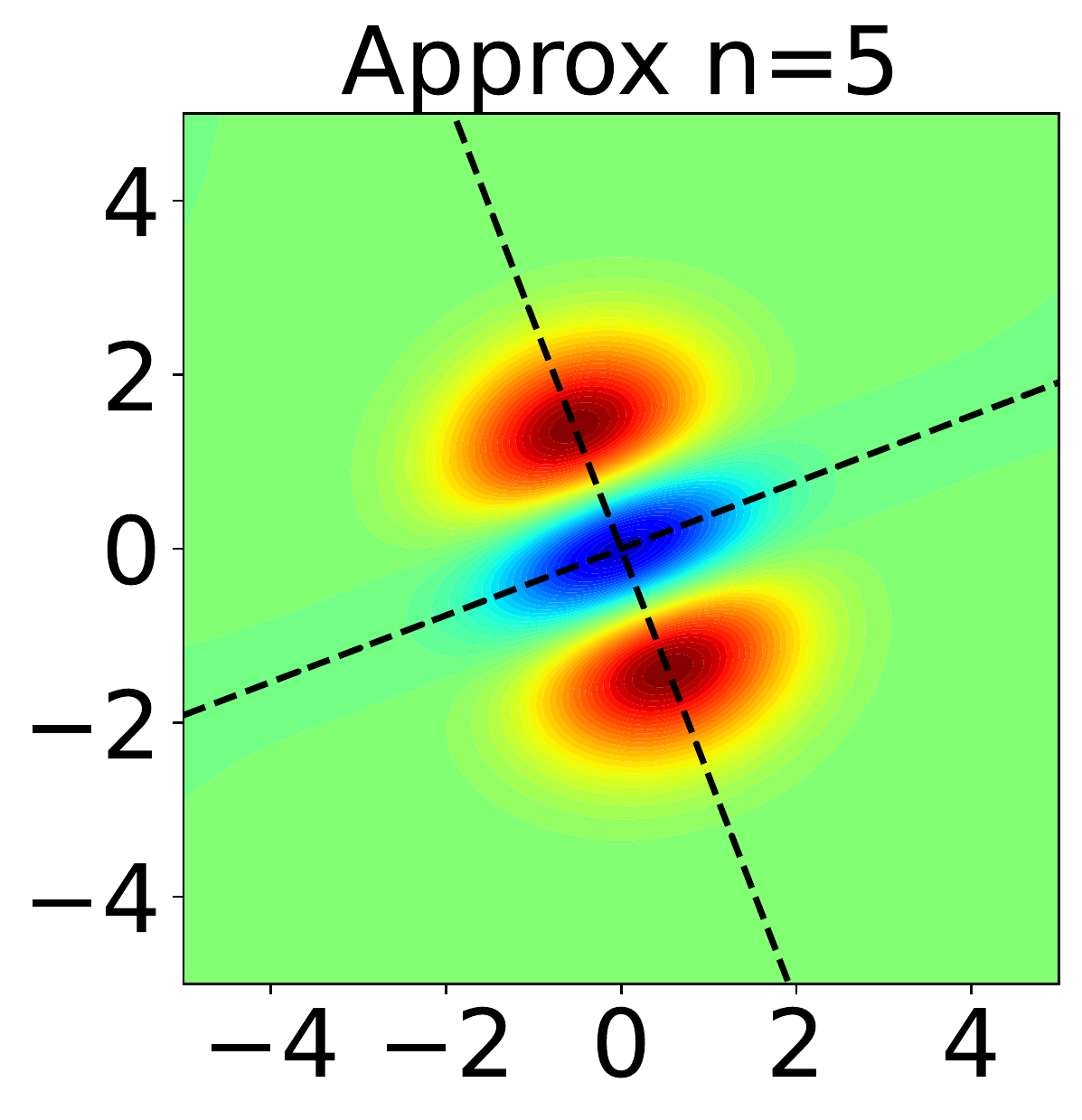}
\includegraphics[width=2.5cm,height=2.5cm]{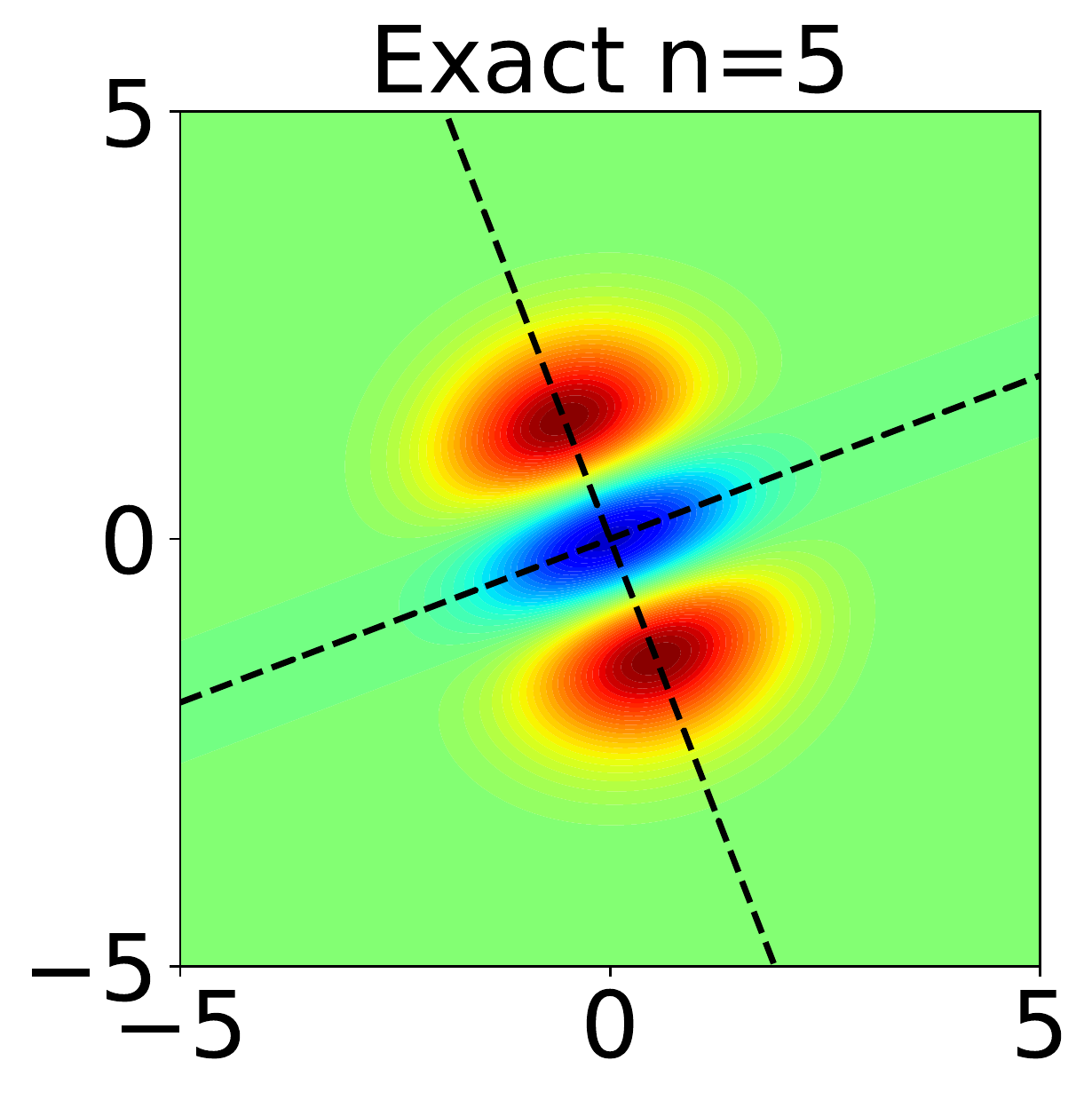}\\
\includegraphics[width=2.5cm,height=2.5cm]{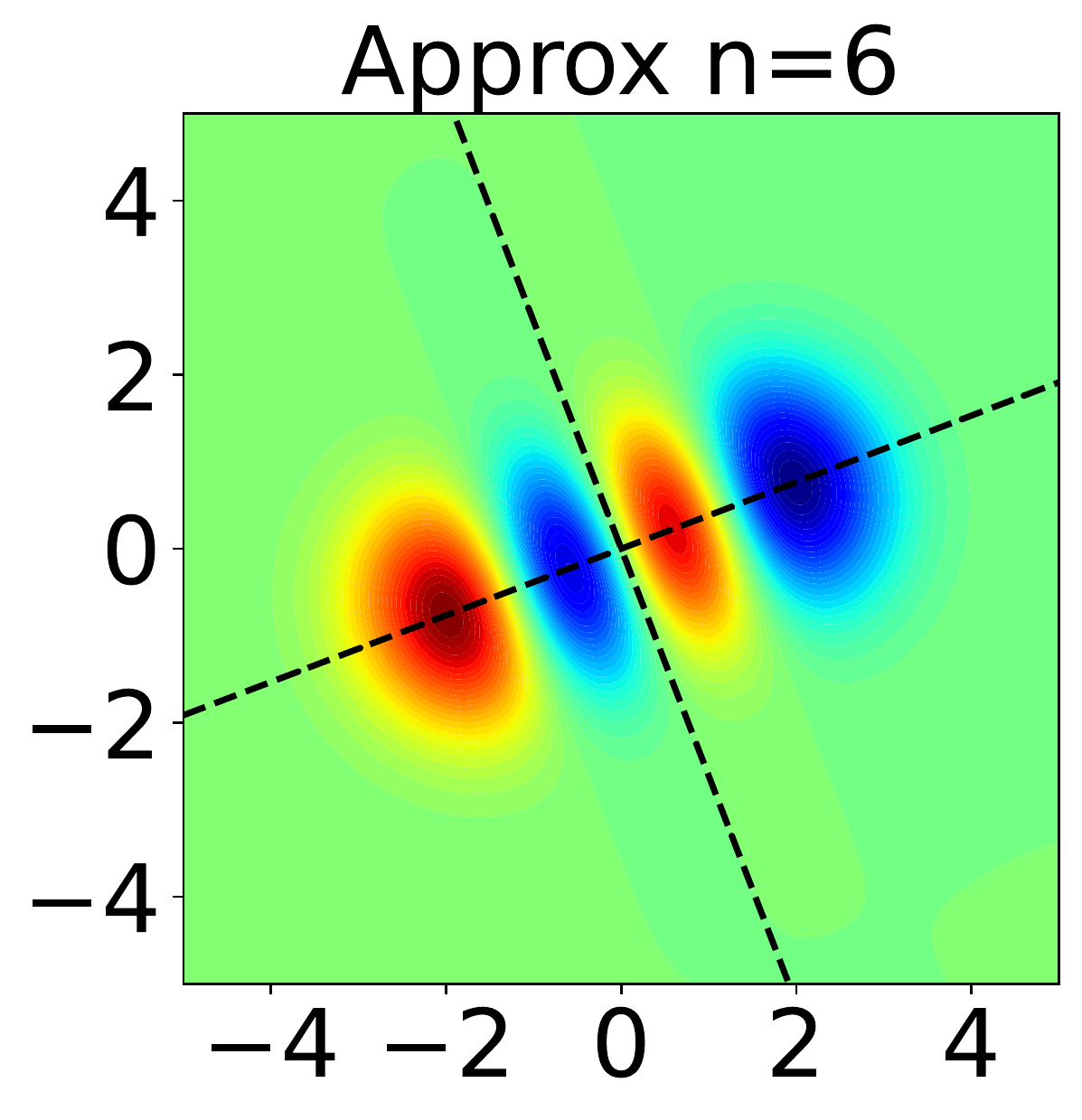}
\includegraphics[width=2.5cm,height=2.5cm]{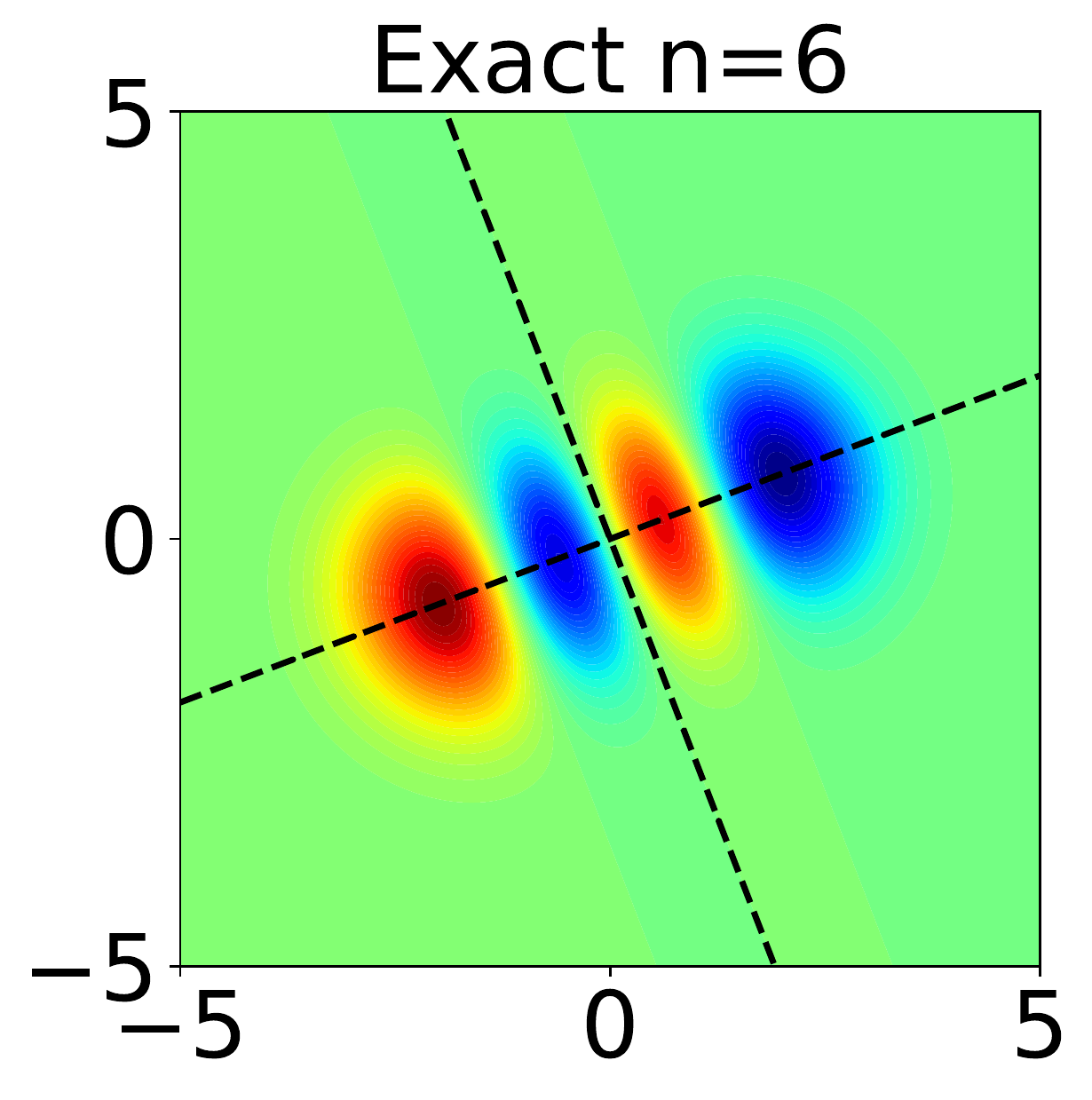}
\includegraphics[width=2.5cm,height=2.5cm]{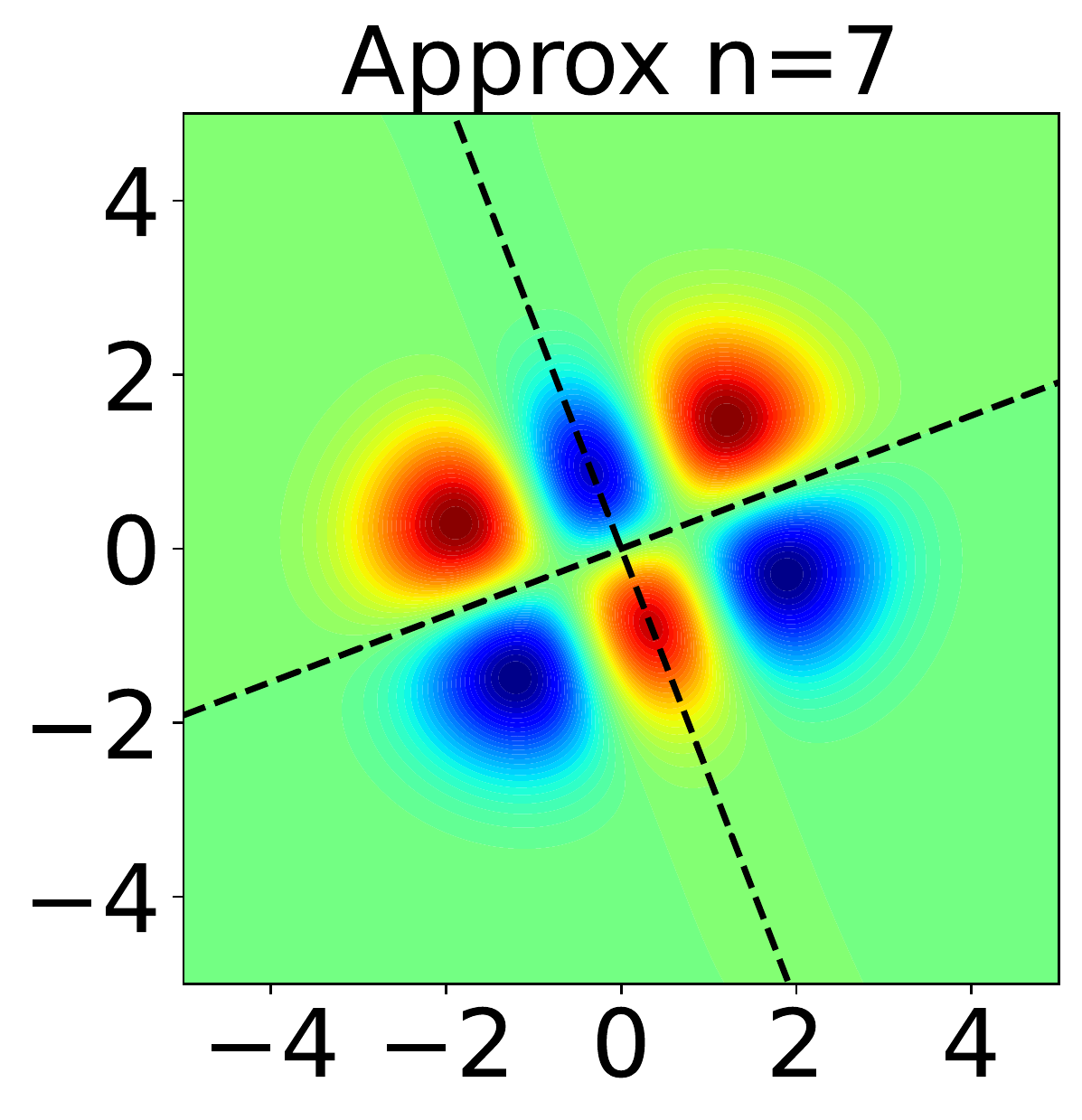}
\includegraphics[width=2.5cm,height=2.5cm]{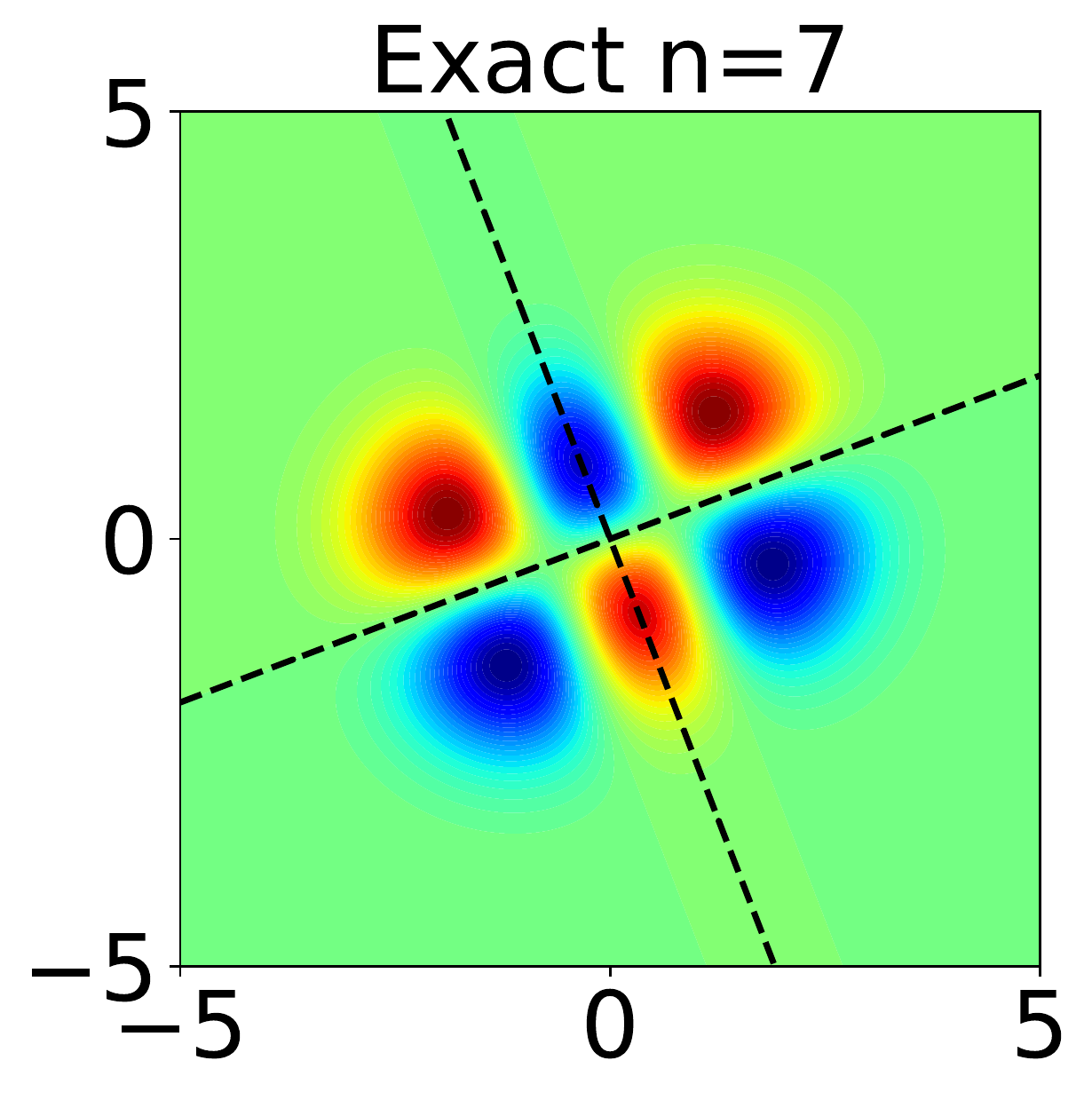}
\includegraphics[width=2.5cm,height=2.5cm]{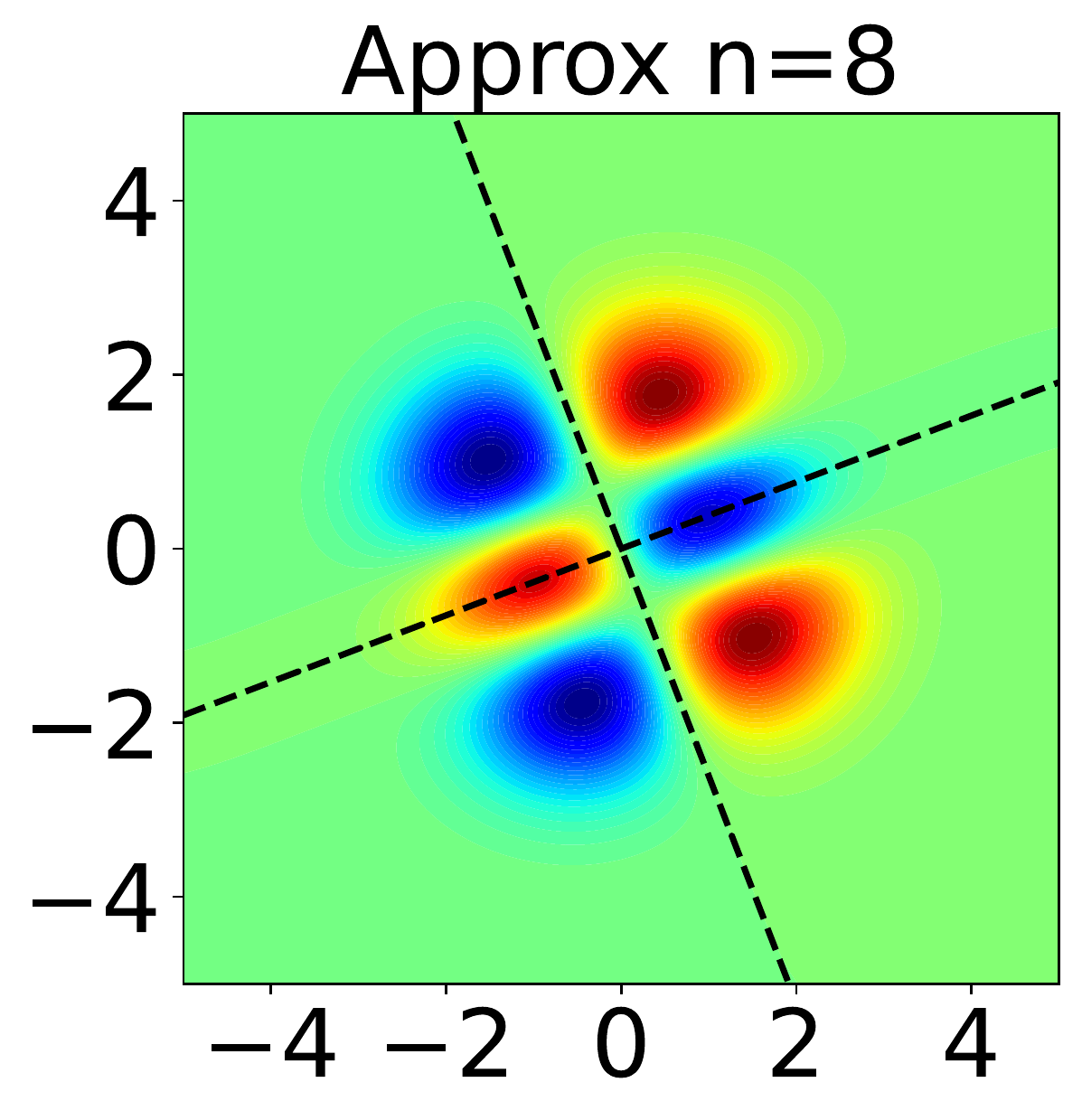}
\includegraphics[width=2.5cm,height=2.5cm]{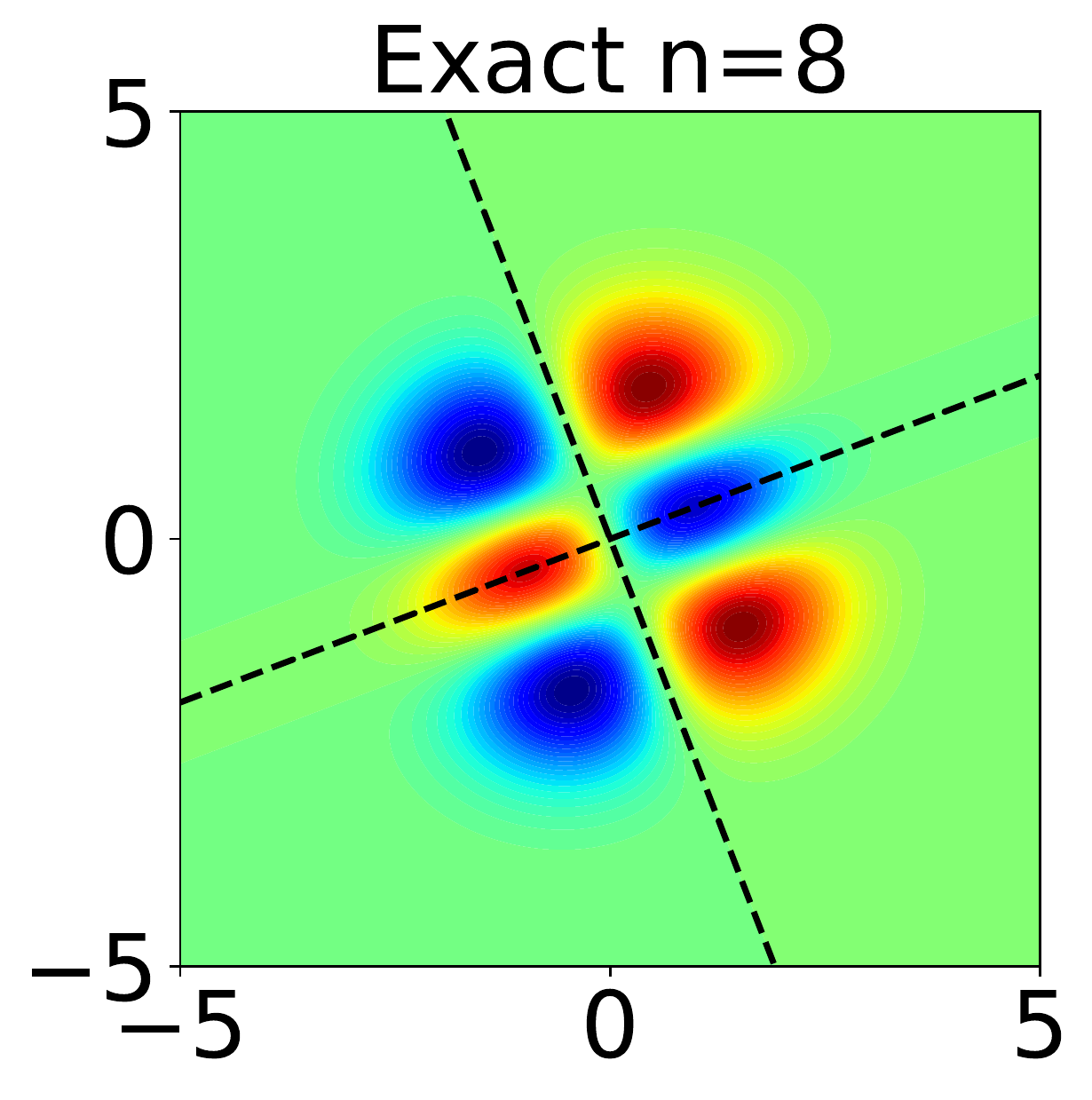}\\
\includegraphics[width=2.5cm,height=2.5cm]{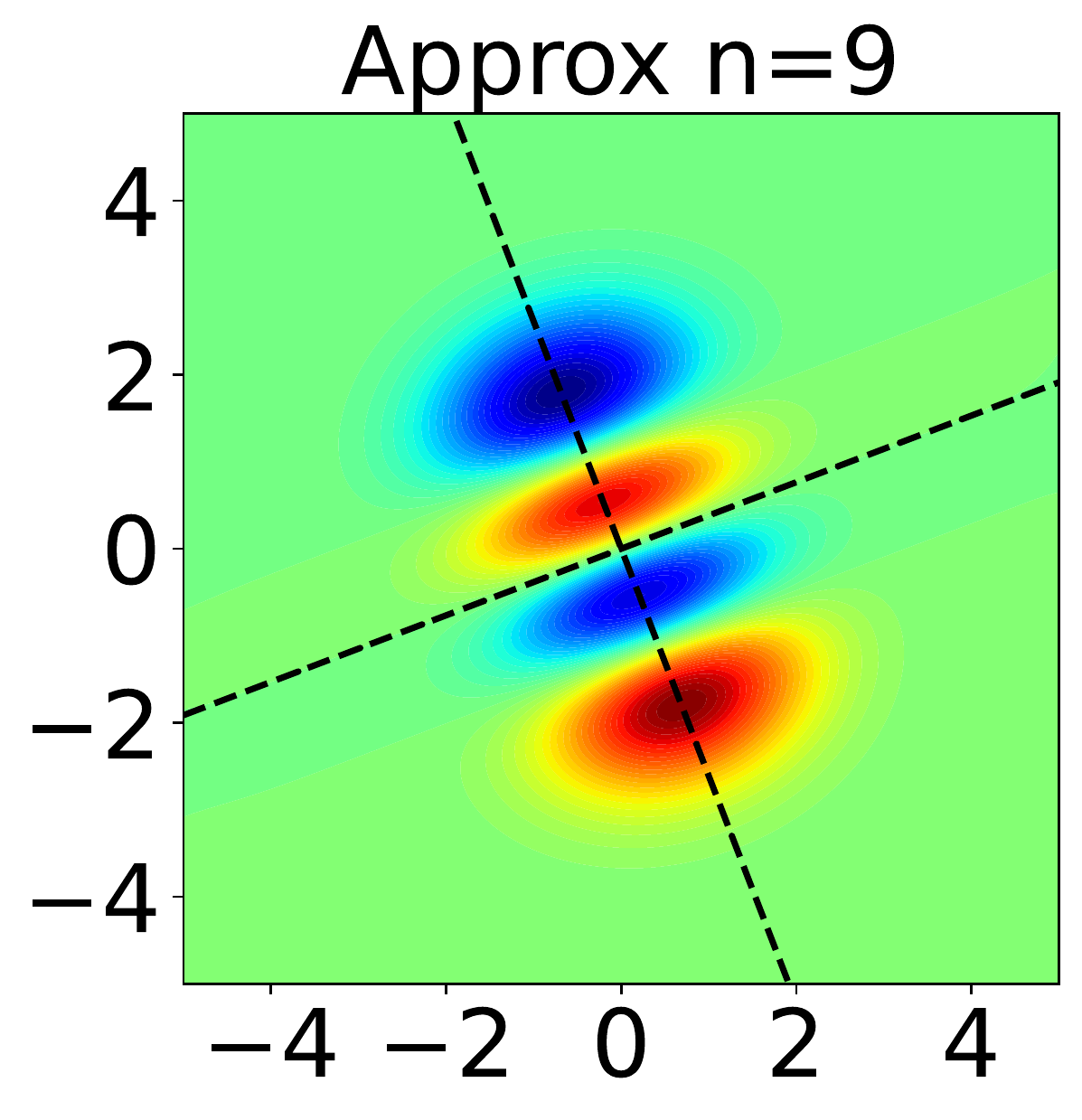}
\includegraphics[width=2.5cm,height=2.5cm]{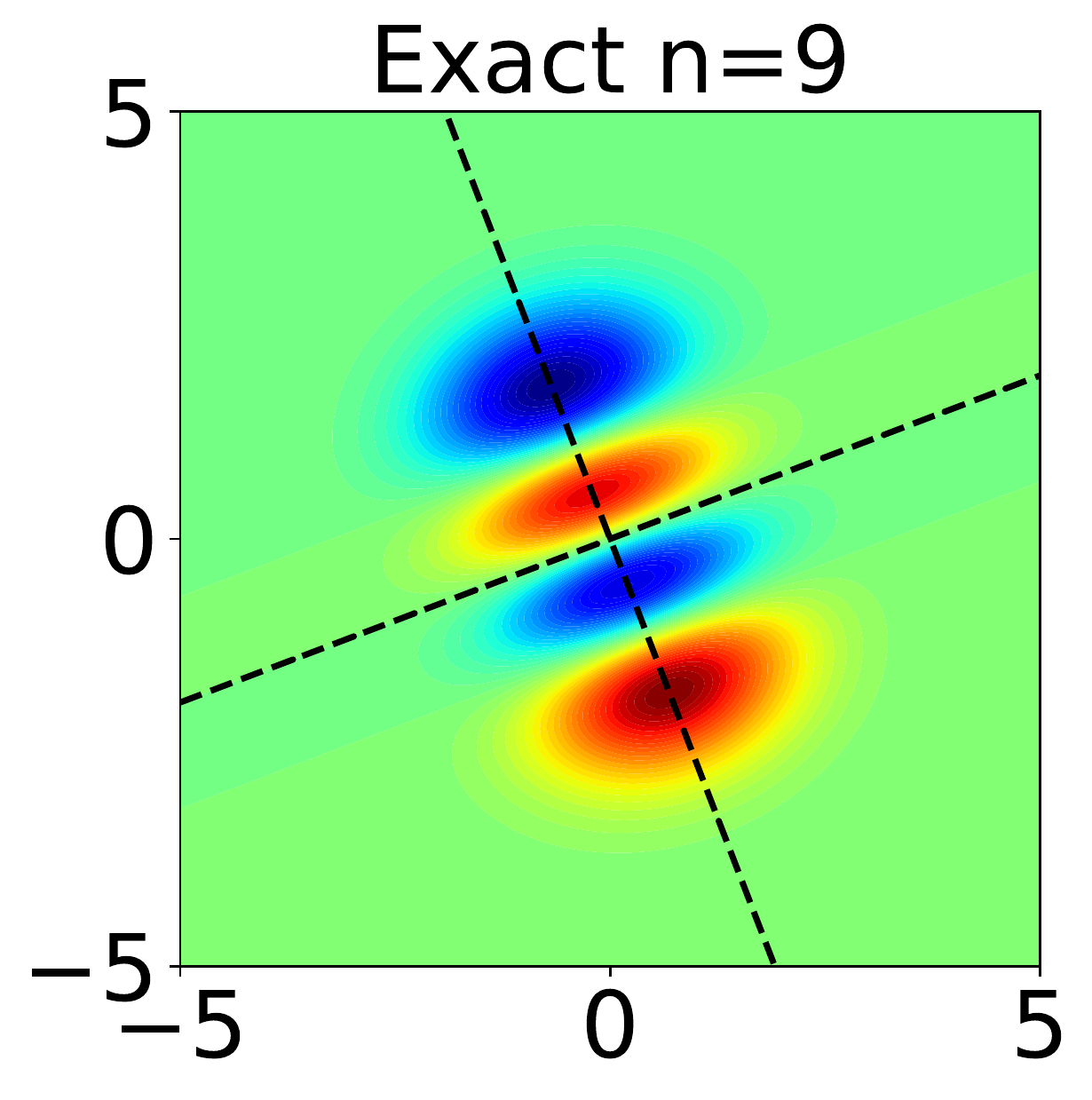}
\includegraphics[width=2.5cm,height=2.5cm]{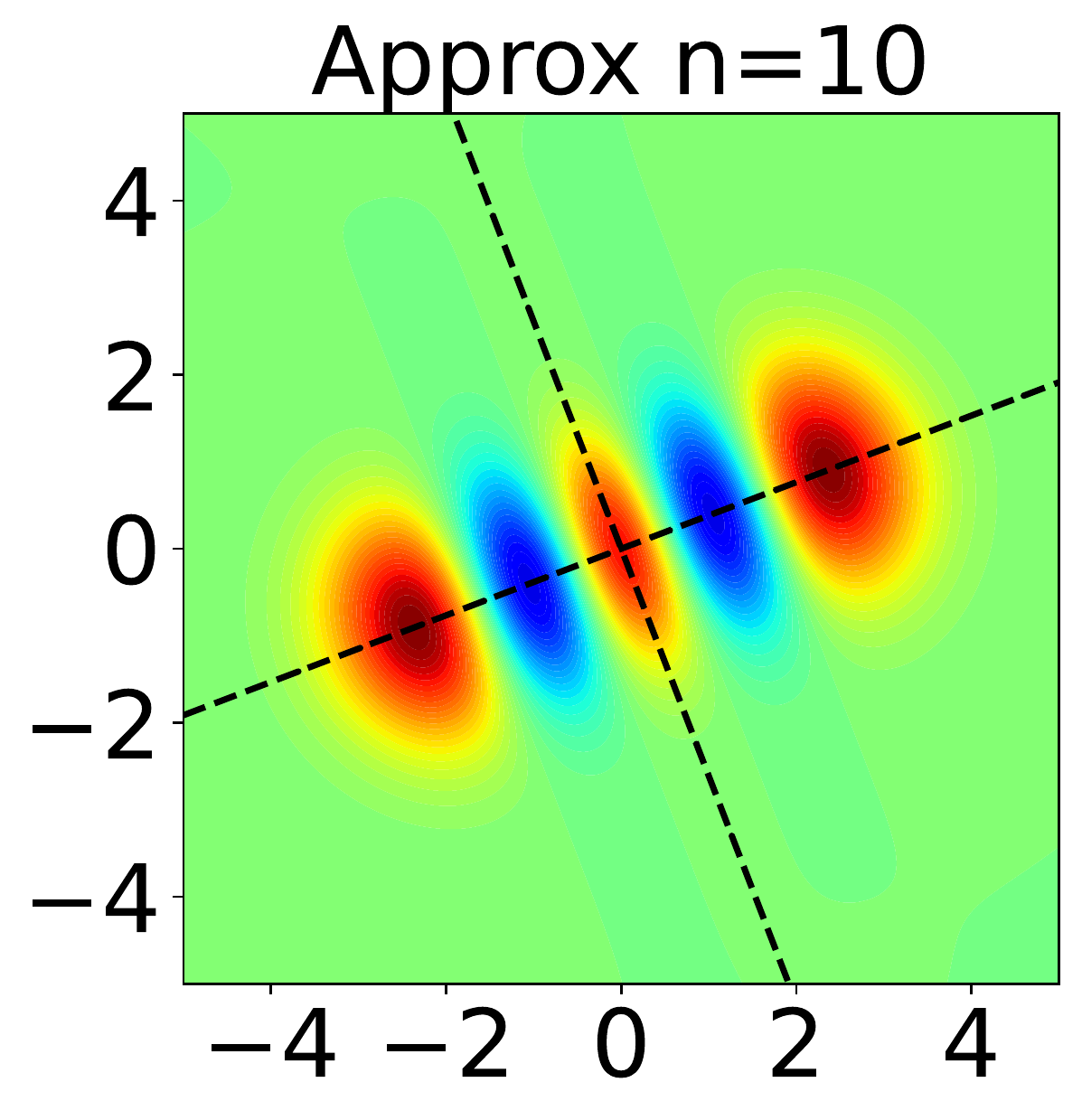}
\includegraphics[width=2.5cm,height=2.5cm]{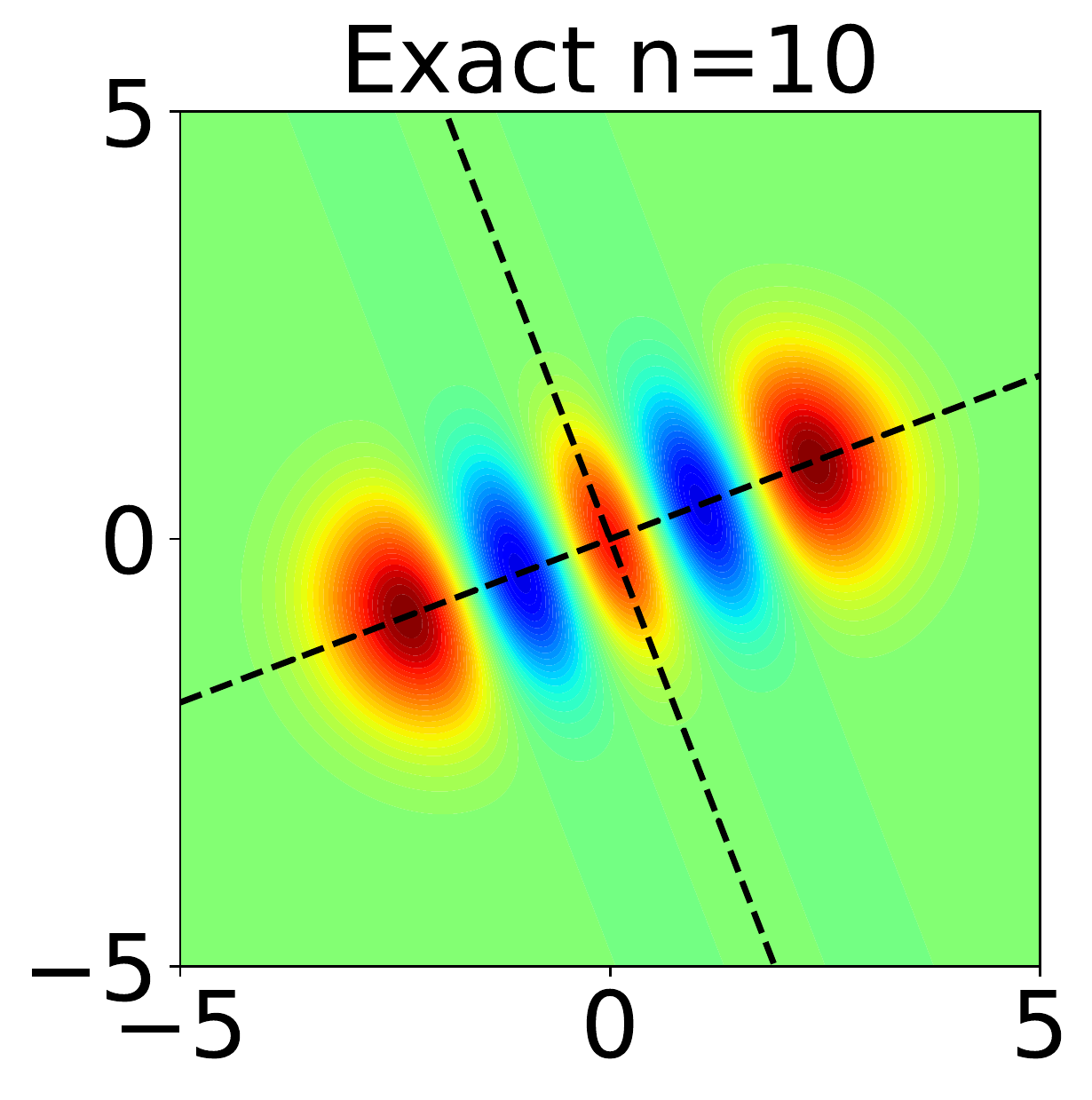}
\includegraphics[width=2.5cm,height=2.5cm]{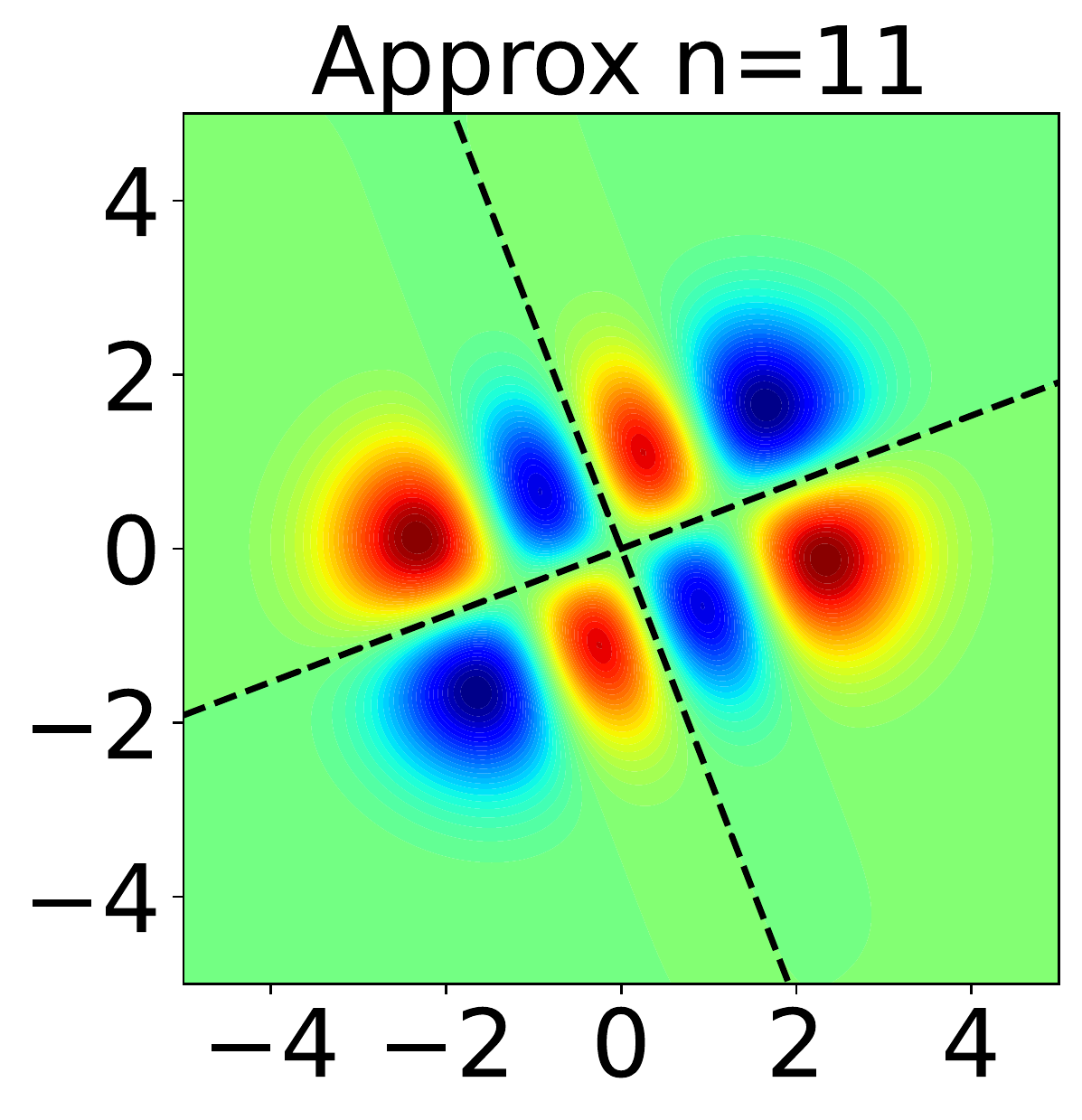}
\includegraphics[width=2.5cm,height=2.5cm]{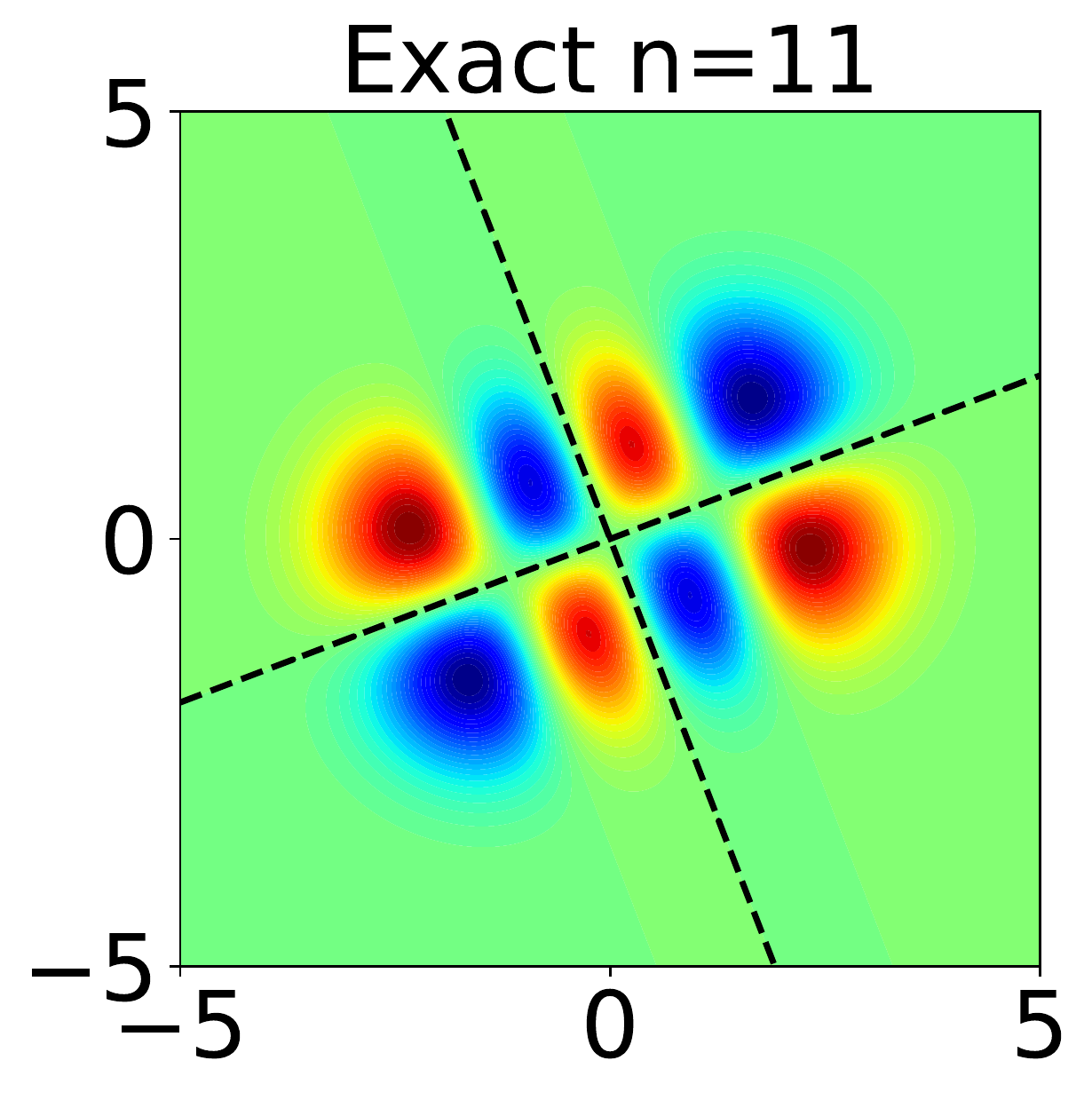}\\
\includegraphics[width=2.5cm,height=2.5cm]{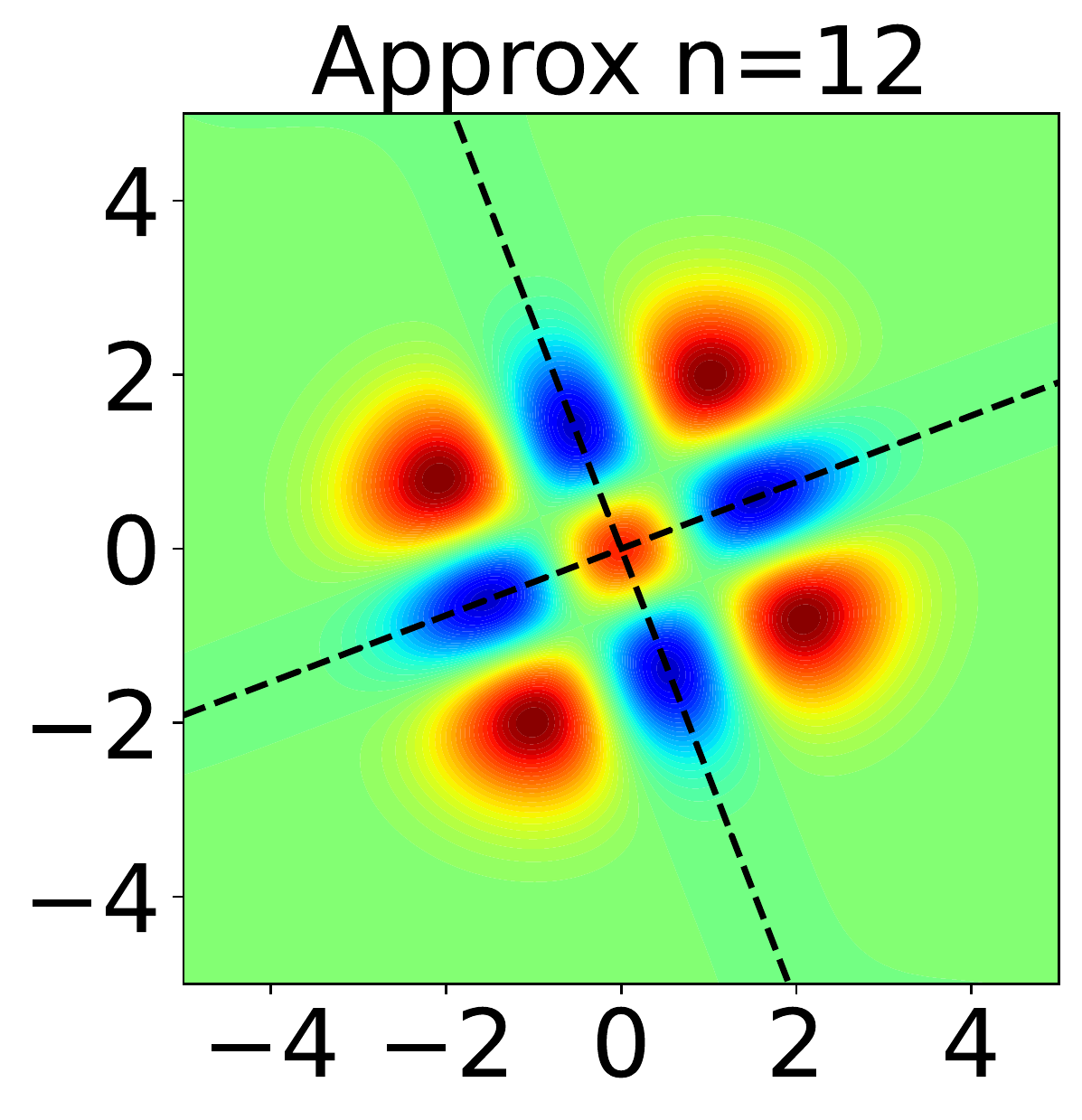}
\includegraphics[width=2.5cm,height=2.5cm]{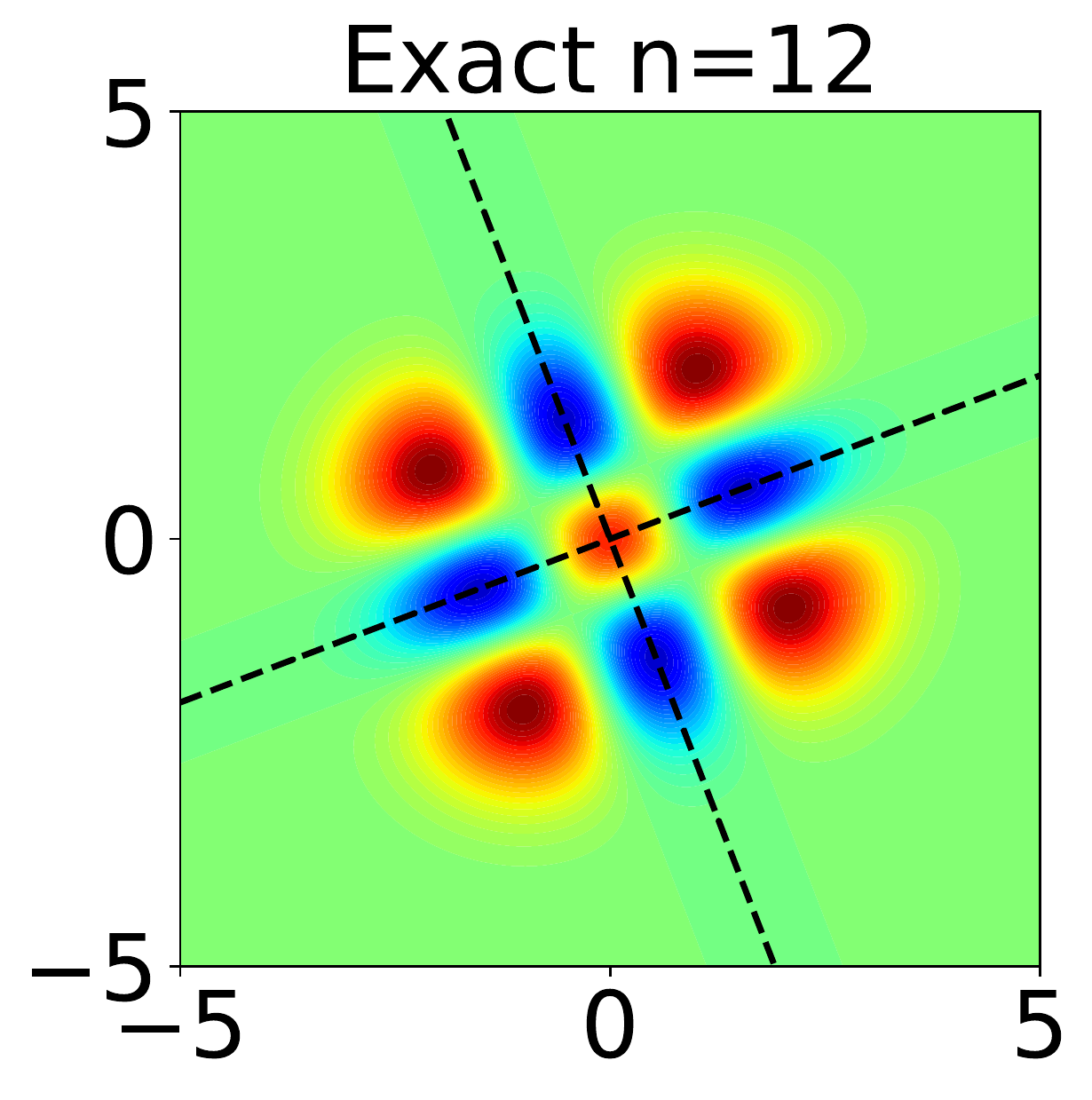}
\includegraphics[width=2.5cm,height=2.5cm]{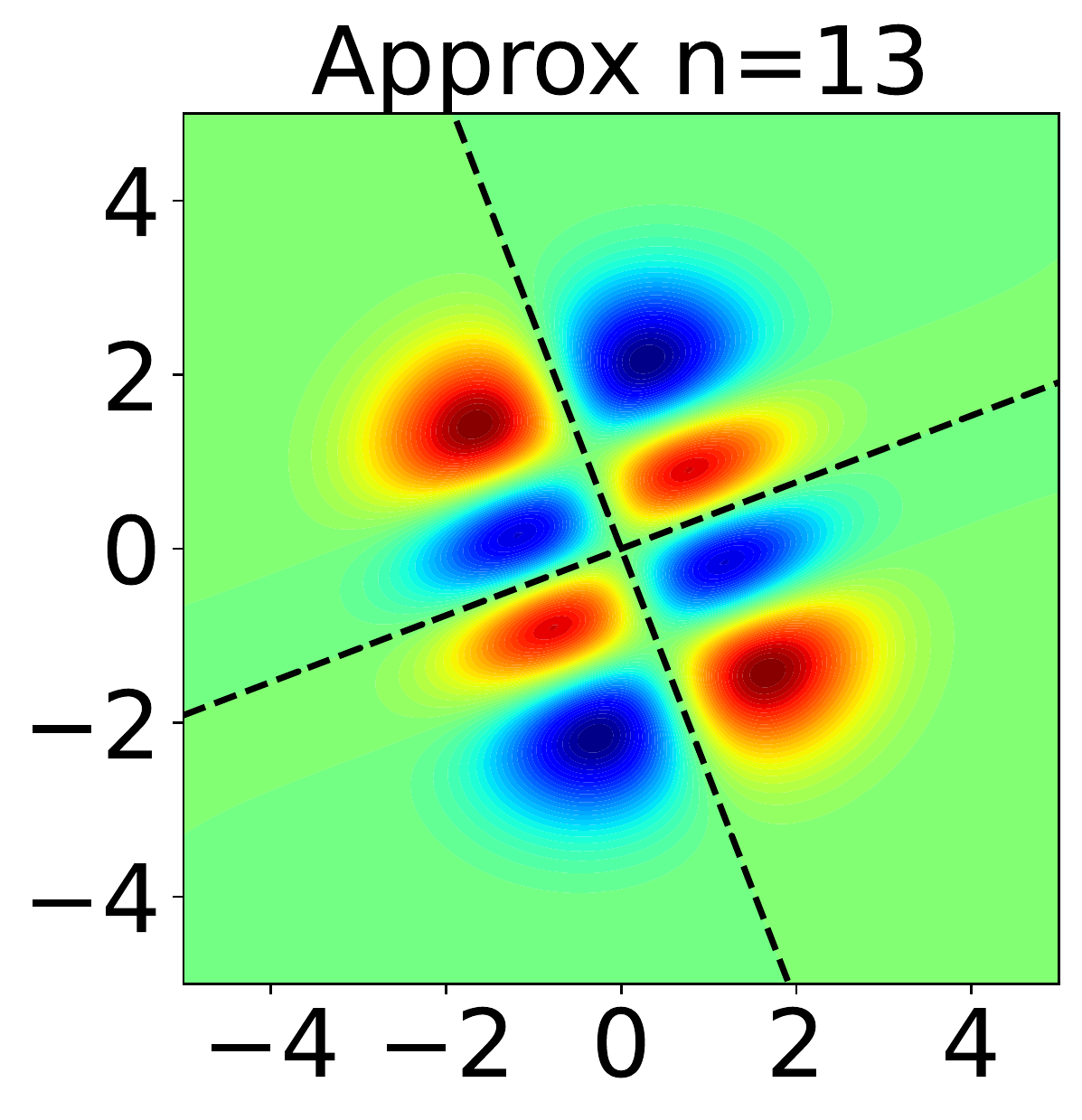}
\includegraphics[width=2.5cm,height=2.5cm]{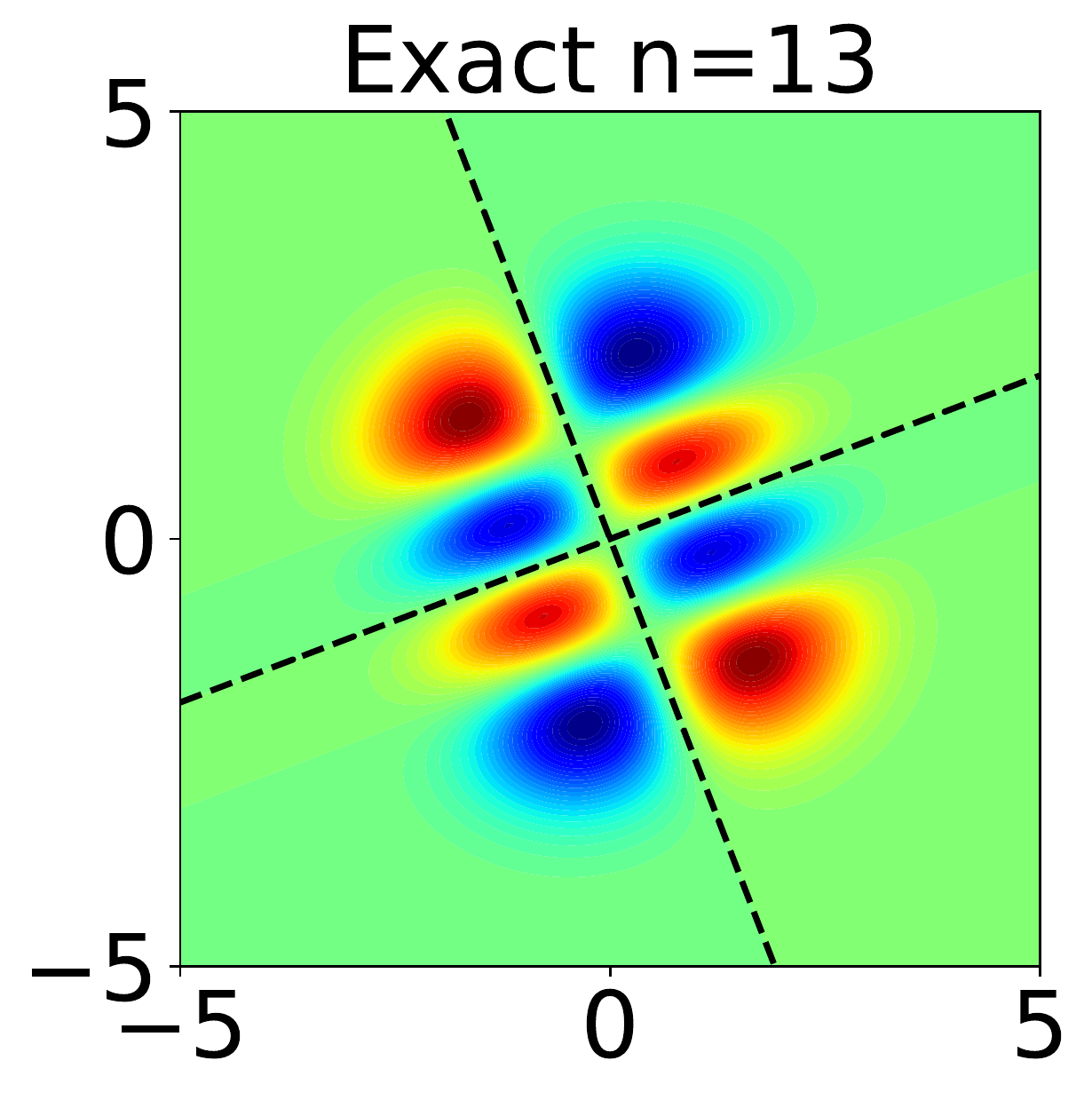}
\includegraphics[width=2.5cm,height=2.5cm]{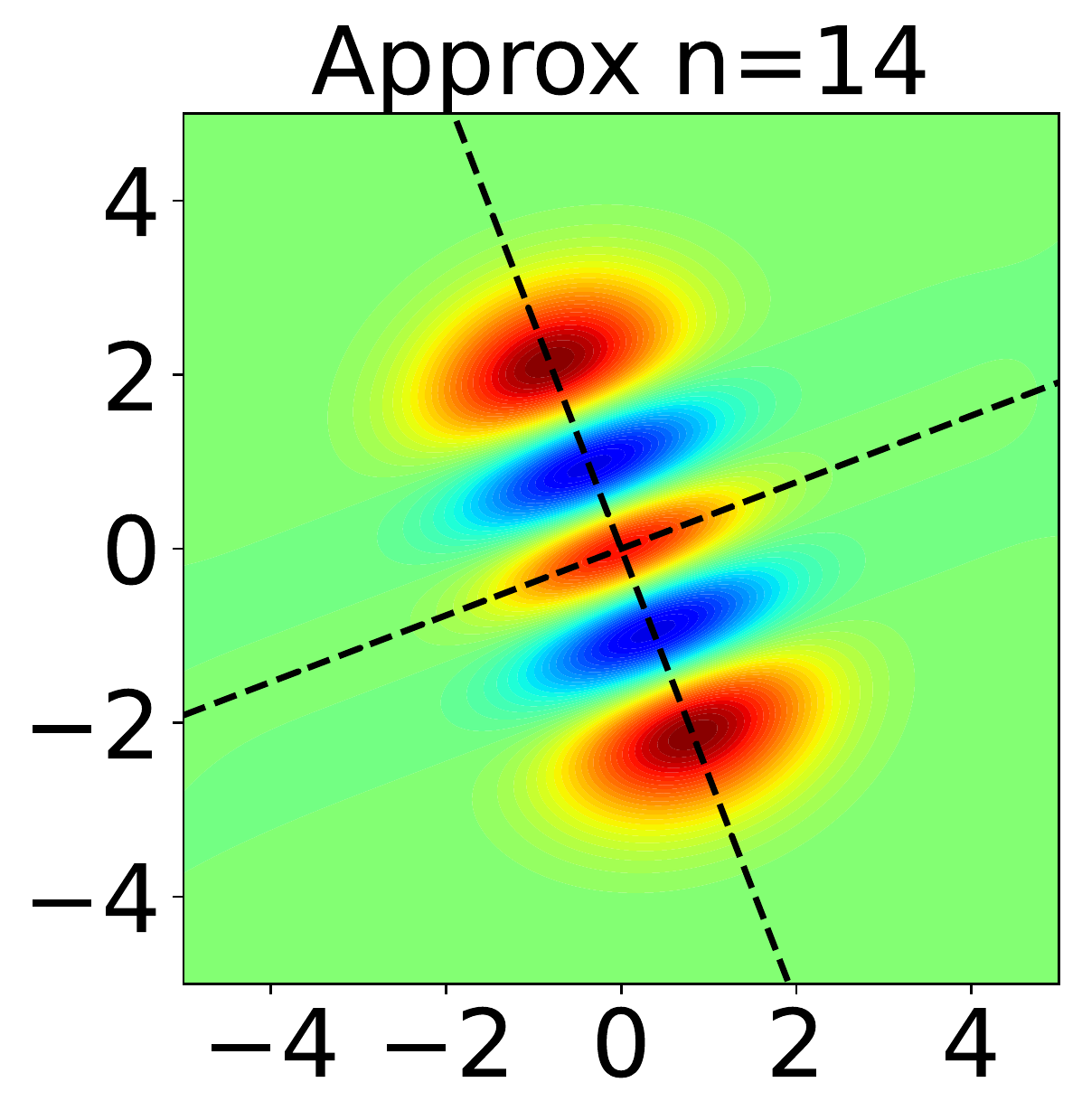}
\includegraphics[width=2.5cm,height=2.5cm]{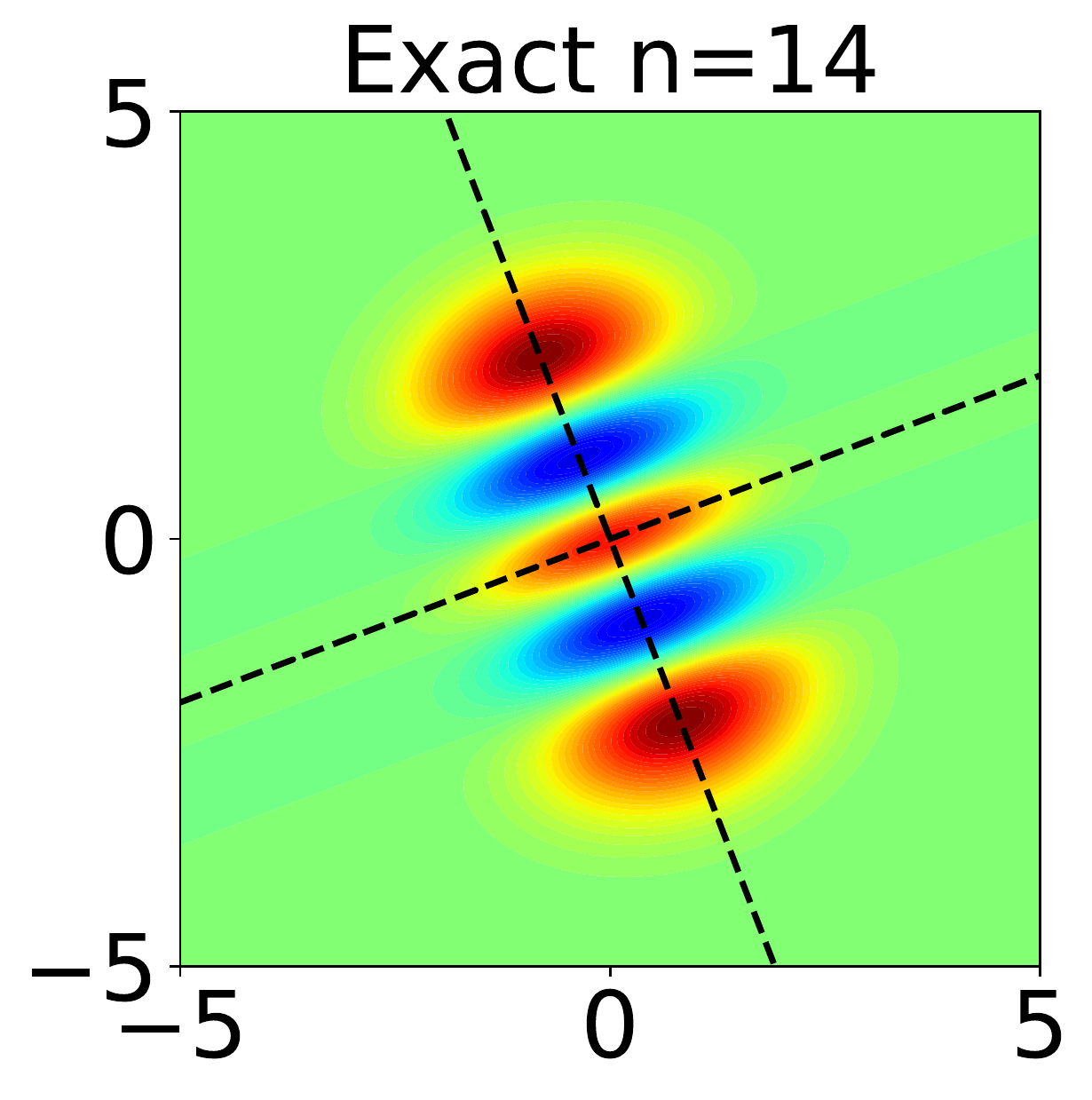}\\
\includegraphics[width=2.5cm,height=2.5cm]{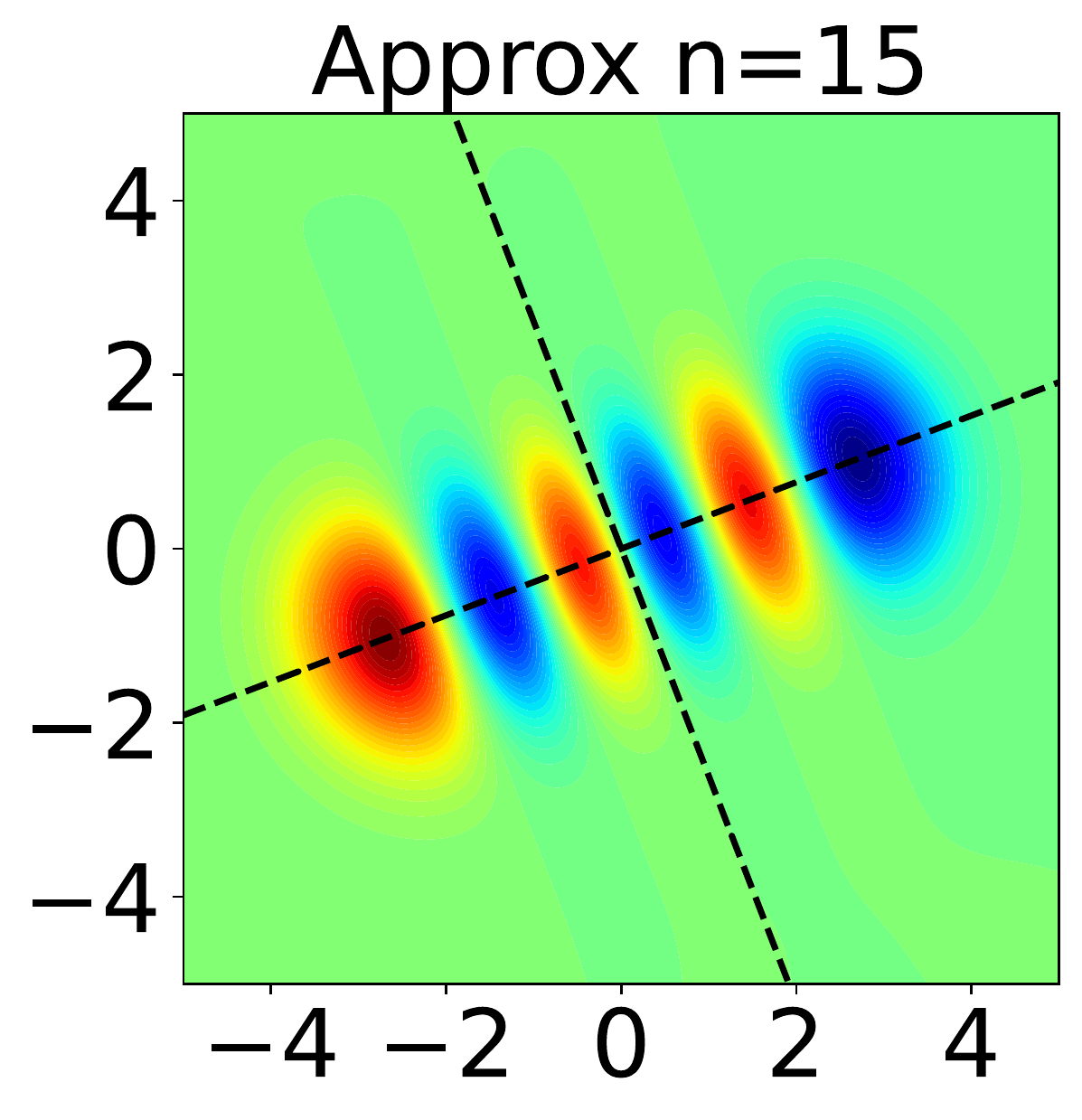}
\includegraphics[width=2.5cm,height=2.5cm]{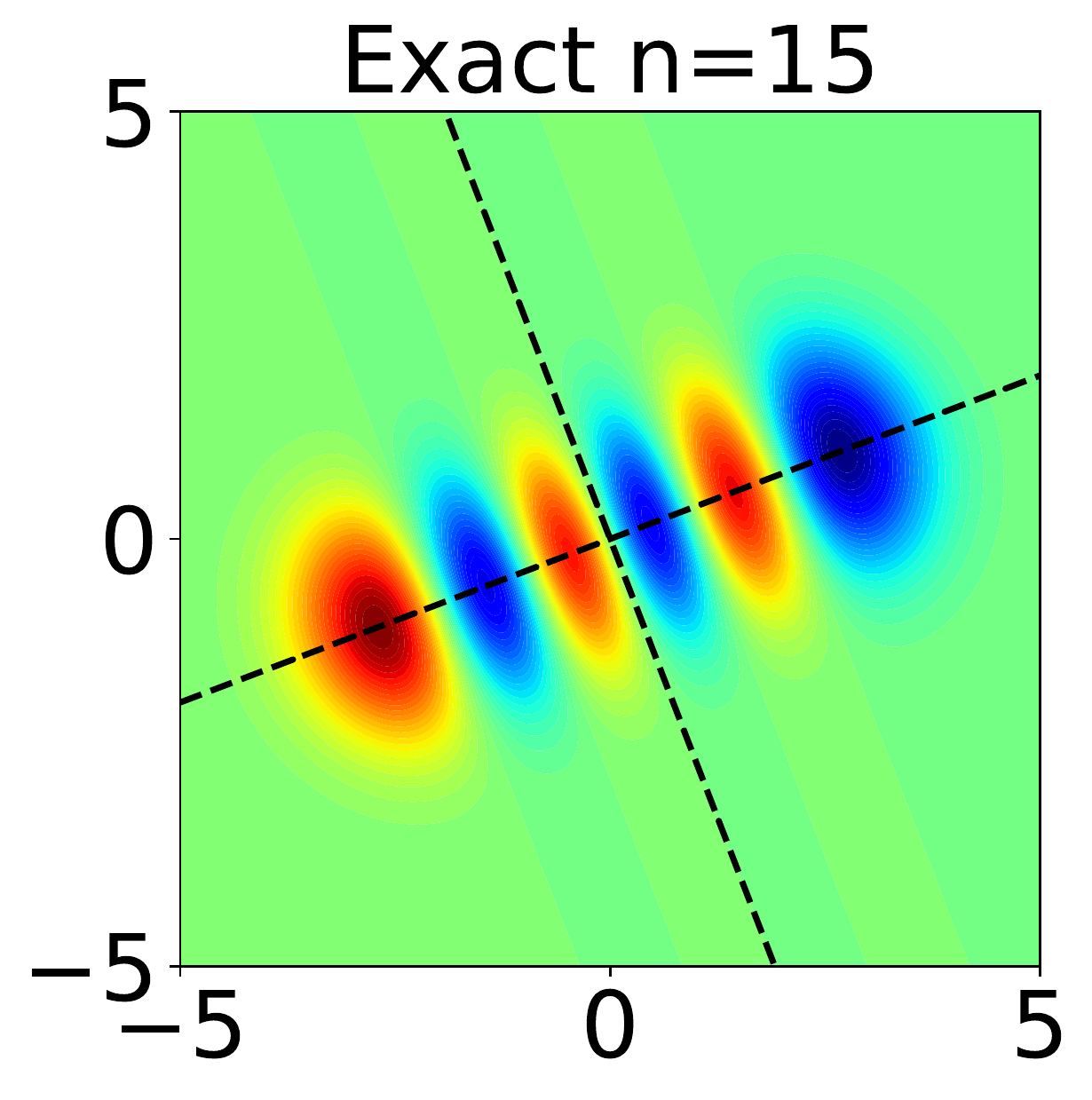}\\
\includegraphics[width=6cm,height=0.5cm]{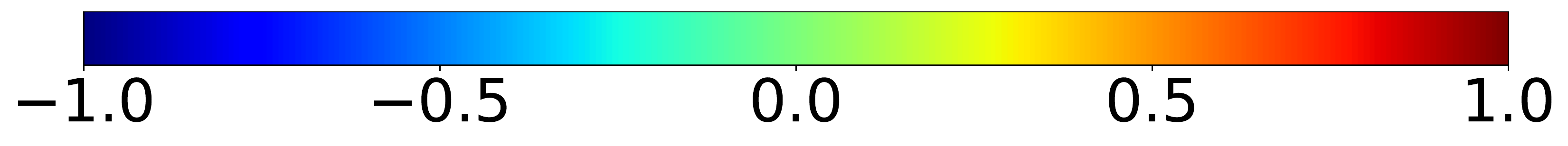}
\caption{The contour plots of the first 16 eigenfunctions for two-dimensional coupled harmonic
oscillator example in coordinate $(x_1,x_2)$. The two dashed lines are
$y_1=0$ and $y_2=0$, respectively.}\label{fig_2dCHO}
\end{figure}

The corresponding approximate wavefunctions obtained by TNN-based machine learning method and the exact wavefunctions
are shown in Figure \ref{fig_2dCHO}, which implies that the approximate wavefunctions have very good accuracy
even for the coupled harmonic oscillator.

\subsubsection{Five-dimensional coupled harmonic oscillator}
Then, we investigate the performance of the proposed method for the case of five-dimensional coupled harmonic oscillator.
Here, the Hamiltonian operator is defined as (\ref{eq_HO_x}) with the matrix $A$ being replaced with the following matrix
\begin{eqnarray}\label{eq_A}
A=
\begin{bmatrix}
  1.05886042&  0.01365034&  0.09163945&  0.11975290&  0.05625013\\
  0.01365034&  1.09613742&  0.10887930&  0.07448974&  0.07407652\\
  0.09163945&  0.10887930&  1.00935913&  0.05588543&  0.08968956\\
  0.11975290&  0.07448974&  0.05588543&  1.17627129&  0.06049045\\
  0.05625013&  0.07407652&  0.08968956&  0.06049045&  0.94969417
\end{bmatrix}.
\end{eqnarray}
Then the exact energy is
\begin{eqnarray}
E_{n_1,n_2,\cdots,n_5}=\sum_{i=1}^5\Big(\frac{1}{2}+n_i\Big)\mu_i^{1/2},
\end{eqnarray}
where $\mu_1=0.88021303$, $\mu_2=0.90973982$, $\mu_3=1.02312382$, $\mu_4=1.10243017$,
$\mu_5=1.37481559$, each of the five quantum numbers $n_i$, $i=1,2,\cdots,5$ takes the values $0,1,2,\cdots$.

%We use 16 TNNs to learn the 16 lowest energy states, each subnetwork of a single TNN with depth 2 and width 50,
%and $p=10$.
%The Adam optimizer is employed with a learning rate 0.001 and epochs of 100000.
%We truncate the computational domain from $\mathbb R^2$ to $[-5,5]^2$, use 100 equal
%subintervals and 16 Gauss points quadrature scheme in each subinterval.
We use $16$ TNNs to learn the lowest $16$ energy states.
A larger TNN structure than the 2-dimensional example is used, the rank is chosen to be $p=50$,
the subnetwork is built with depth $3$ and width $100$.
The Adam optimizer is employed with a learning rate 0.001 and epochs of 500000.
Then the final result is given by the subsequent 10000 steps LBFGS.
The same 99 points Hermite-Gauss quadrature scheme as that in the last two examples is sufficient to this example.

Table \ref{table_5dCHO} shows the corresponding numerical results, where we can find the
proposed numerical method can obtain obviously better accuracy than that in \cite{LiZhaiChen}.
\begin{table}[!htb]
\caption{Errors of five-dimensional coupled harmonic oscillator problem for the 16 lowest energy states.}\label{table_5dCHO}
\begin{center}
\begin{tabular}{ccccc}
\hline
$n$&  $(n_1,n_2,n_3,n_4,n_5)$&   Exact $E_n$&  Approx $E_n$&  ${\rm err}_E$\\
\hline
0&   (0,0,0,0,0)&   2.562993697776131&   2.562993699775476&   7.801e-10\\
1&   (1,0,0,0,0)&   3.501190387362160&   3.501190399601748&   3.496e-09\\
2&   (0,1,0,0,0)&   3.516796517763949&   3.516796531293733&   3.847e-09\\
3&   (0,0,1,0,0)&   3.574489532179441&   3.574489543933587&   3.288e-09\\
4&   (0,0,0,1,0)&   3.612960443213281&   3.612960454781468&   3.202e-09\\
5&   (0,0,0,0,1)&   3.735519003914085&   3.735519020687214&   4.490e-09\\
6&   (2,0,0,0,0)&   4.439387076948189&   4.439387114214263&   8.394e-09\\
7&   (1,1,0,0,0)&   4.454993207349979&   4.454993243128091&   8.031e-09\\
8&   (0,2,0,0,0)&   4.470599337751768&   4.470599382176197&   9.937e-09\\
9&   (1,0,1,0,0)&   4.512686221765470&   4.512686265259870&   9.638e-09\\
10&  (0,1,1,0,0)&   4.528292352167259&   4.528292394259157&   9.295e-09\\
11&  (1,0,0,1,0)&   4.551157132799310&   4.551157174501202&   9.163e-09\\
12&  (0,1,0,1,0)&   4.566763263201100&   4.566763304302355&   9.000e-09\\
13&  (0,0,2,0,0)&   4.585985366582751&   4.585985402475187&   7.827e-09\\
14&  (0,0,1,1,0)&   4.624456277616591&   4.624456313384630&   7.735e-09\\
15&  (0,0,0,2,0)&   4.662927188650432&   4.662927231227110&   9.131e-09\\
\hline
\end{tabular}
\end{center}
\end{table}

\begin{figure}[htb!]
\centering
\includegraphics[width=3.75cm,height=3.75cm]{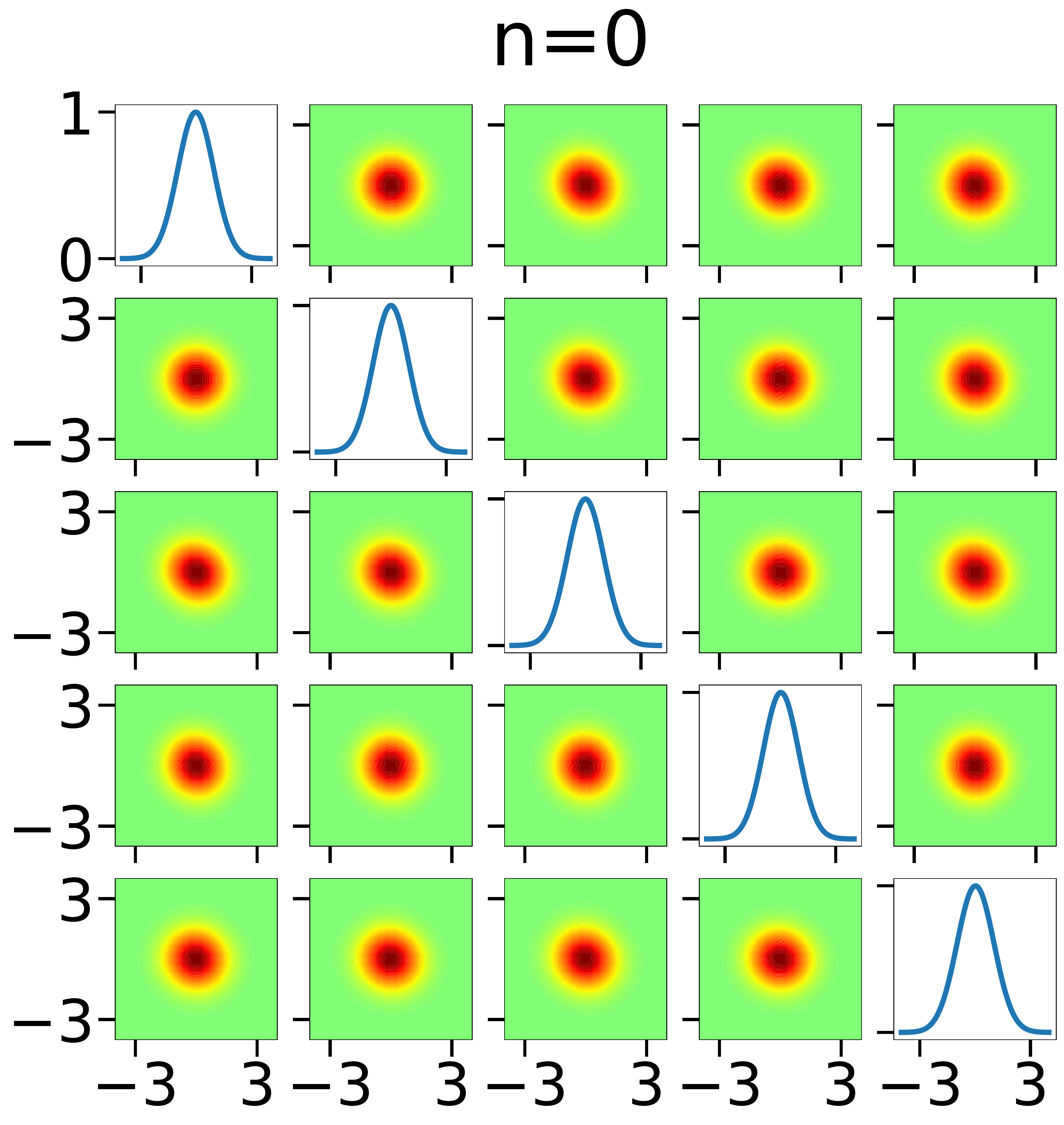}
\includegraphics[width=3.75cm,height=3.75cm]{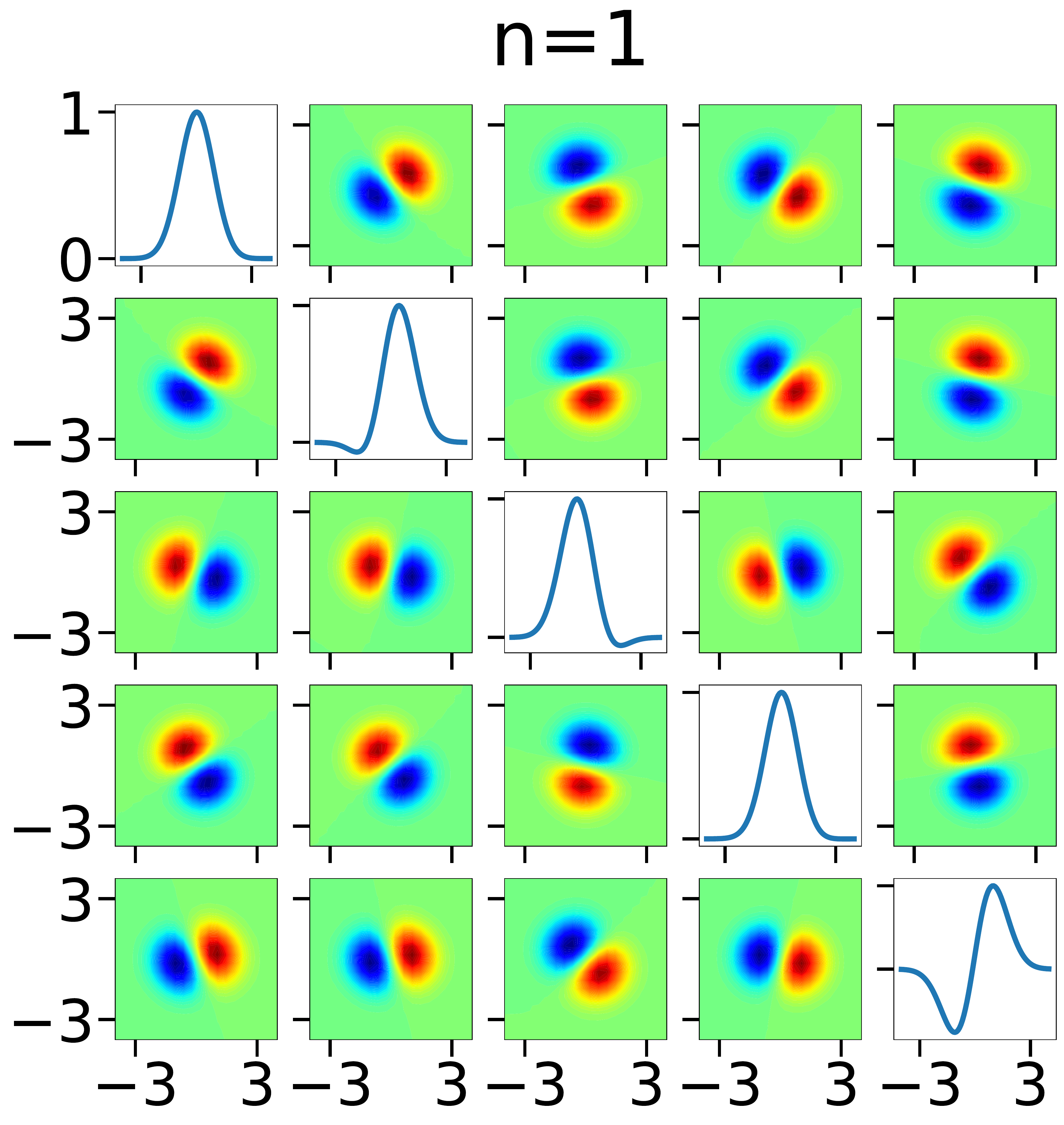}
\includegraphics[width=3.75cm,height=3.75cm]{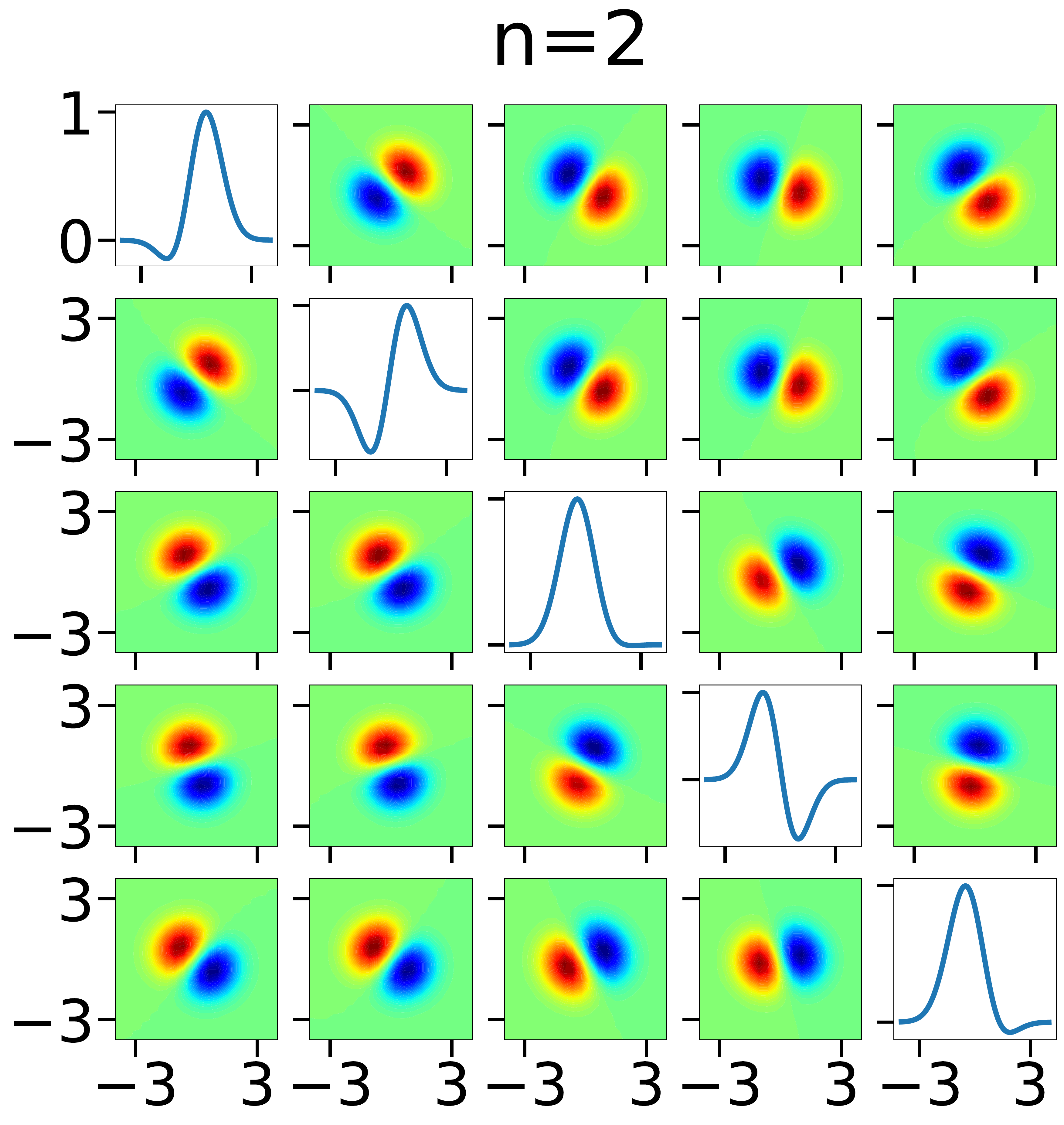}
\includegraphics[width=3.75cm,height=3.75cm]{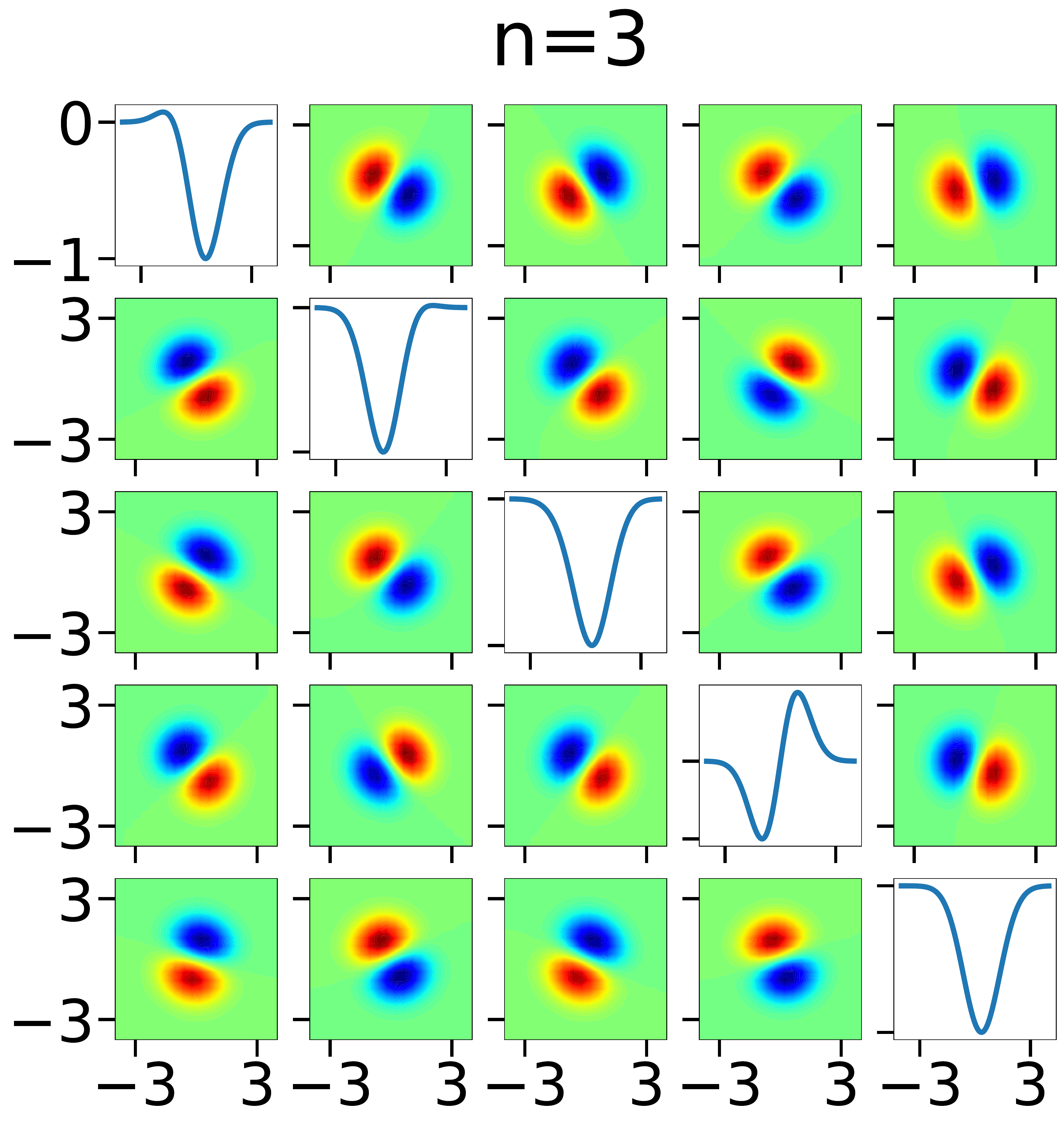}\\
\includegraphics[width=3.75cm,height=3.75cm]{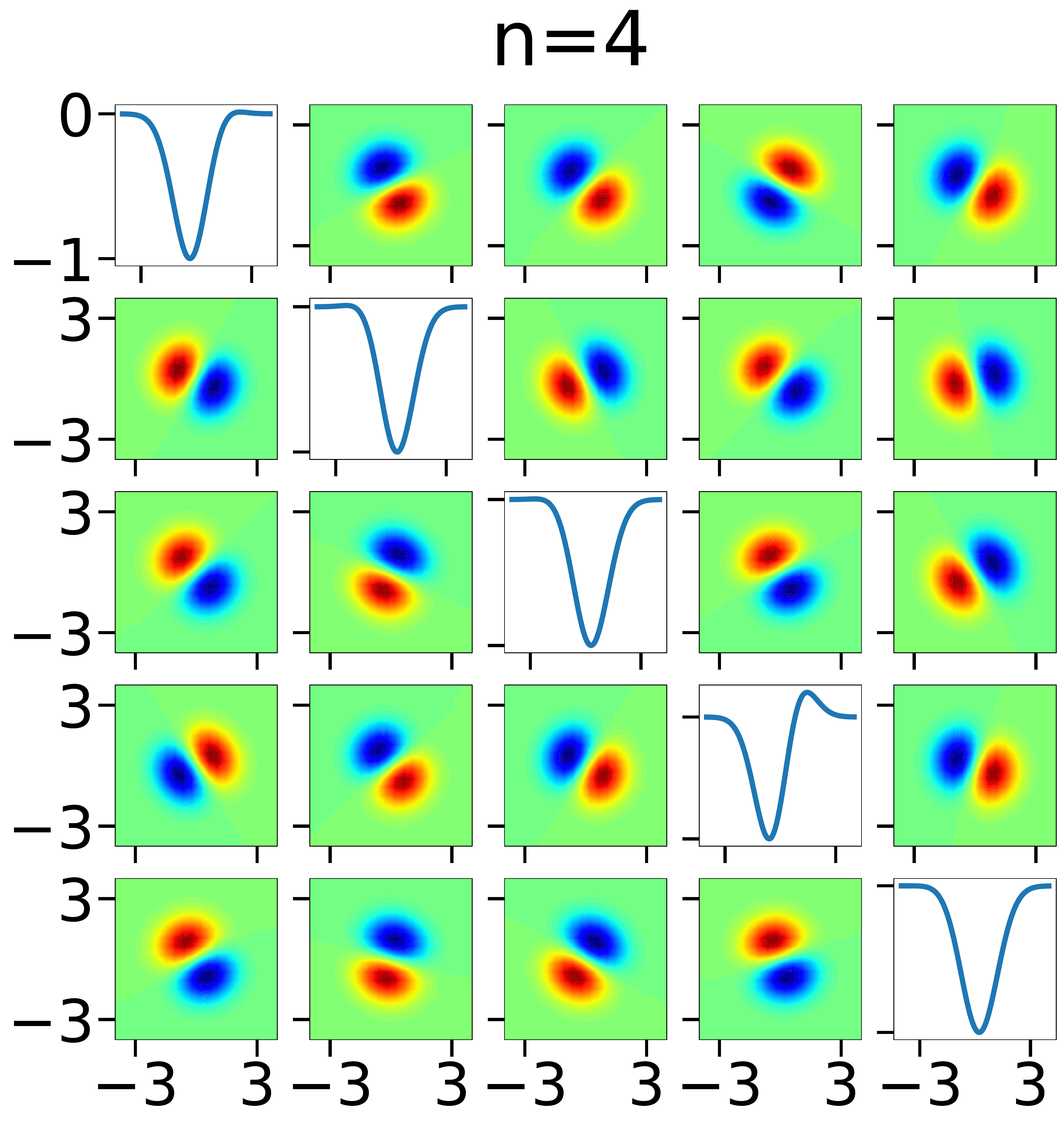}
\includegraphics[width=3.75cm,height=3.75cm]{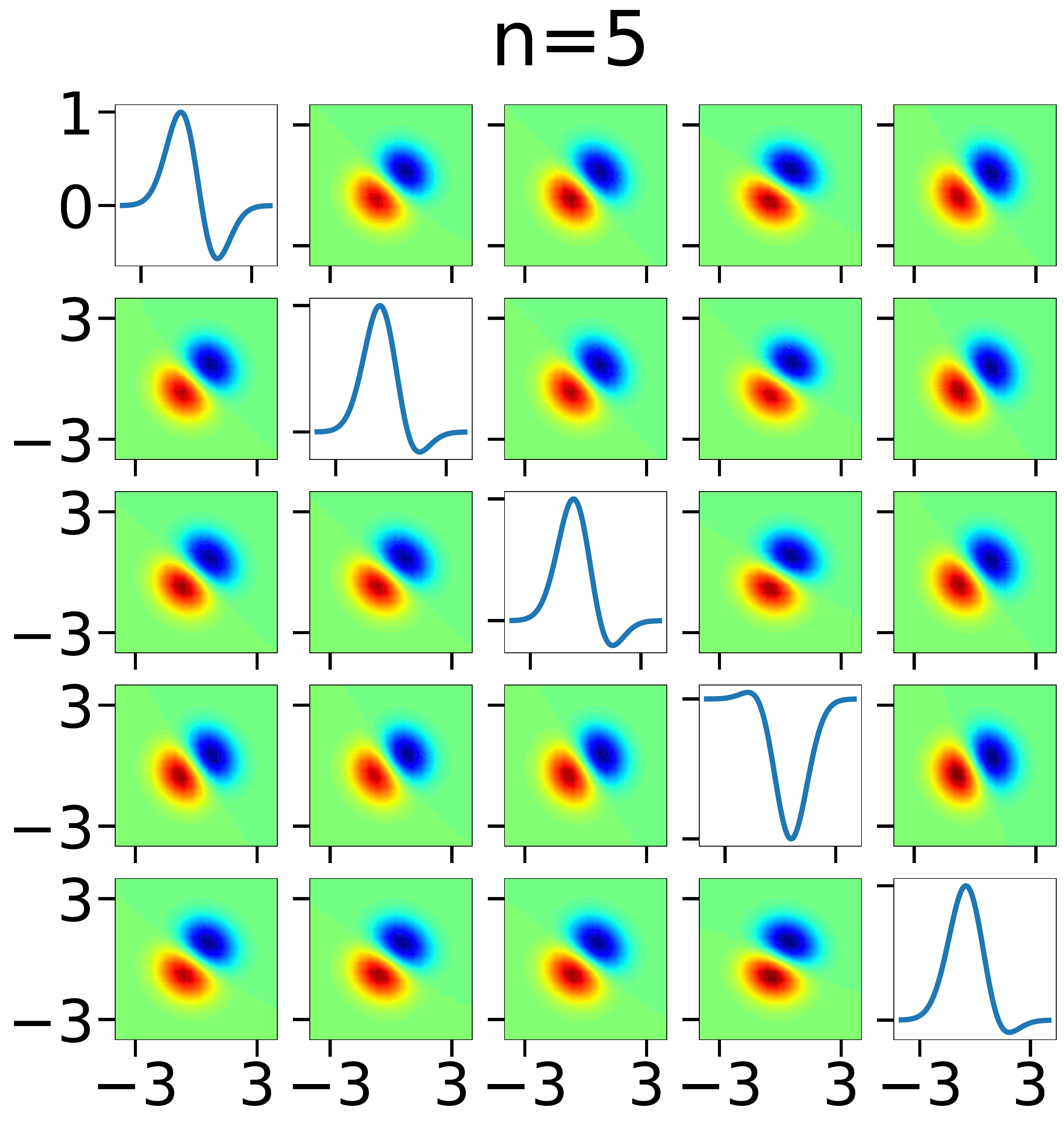}
\includegraphics[width=3.75cm,height=3.75cm]{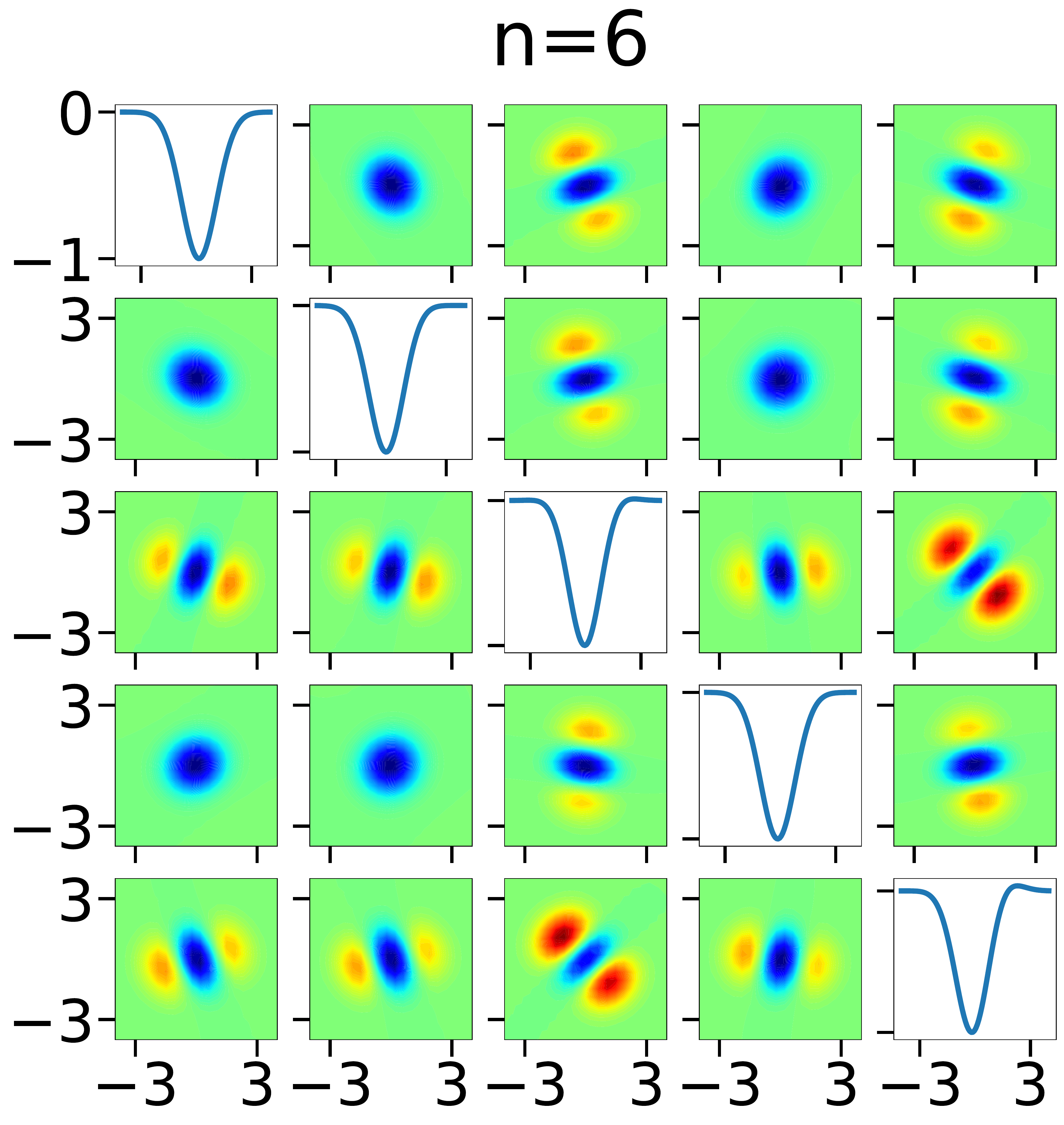}
\includegraphics[width=3.75cm,height=3.75cm]{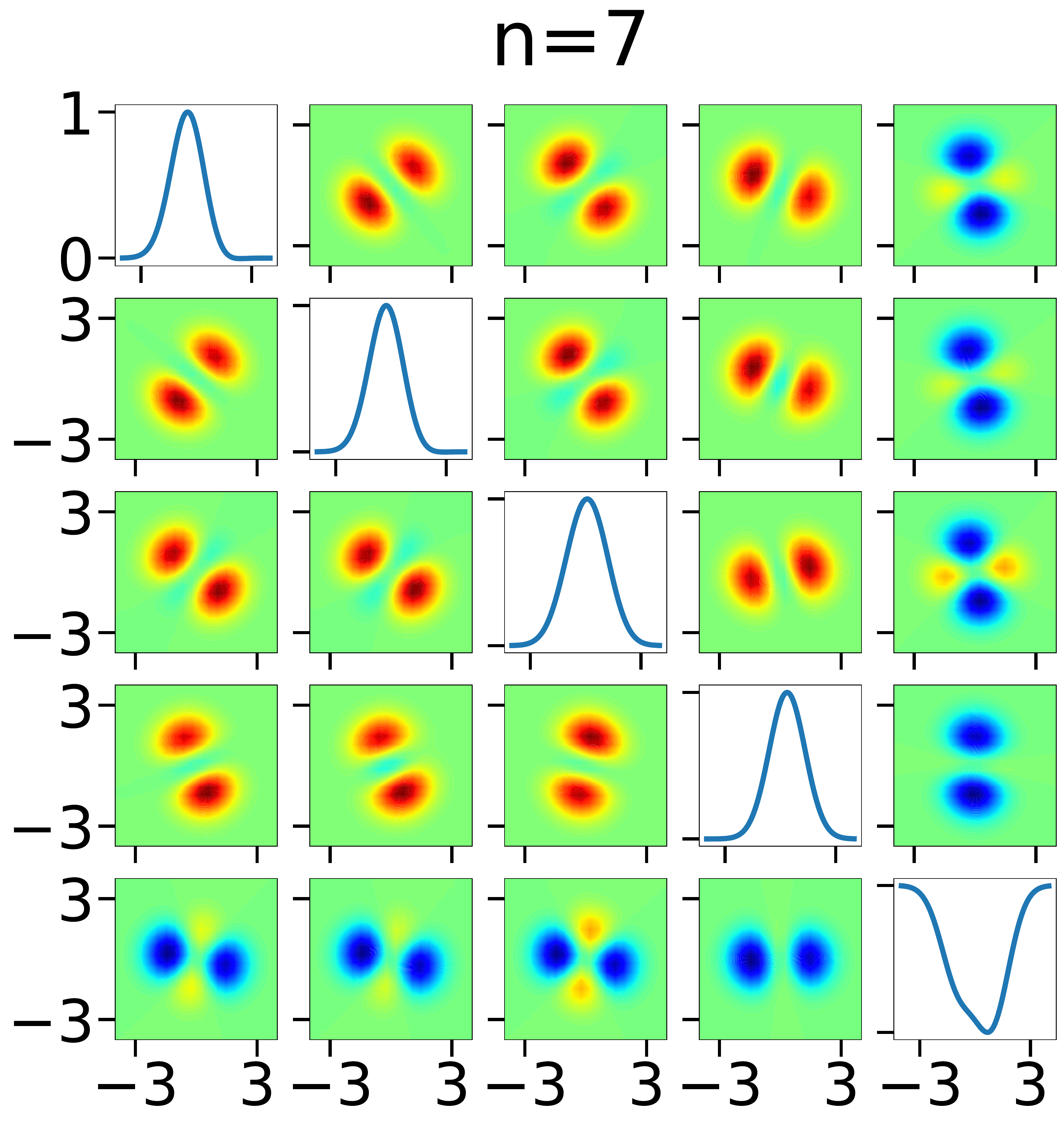}\\
\includegraphics[width=3.75cm,height=3.75cm]{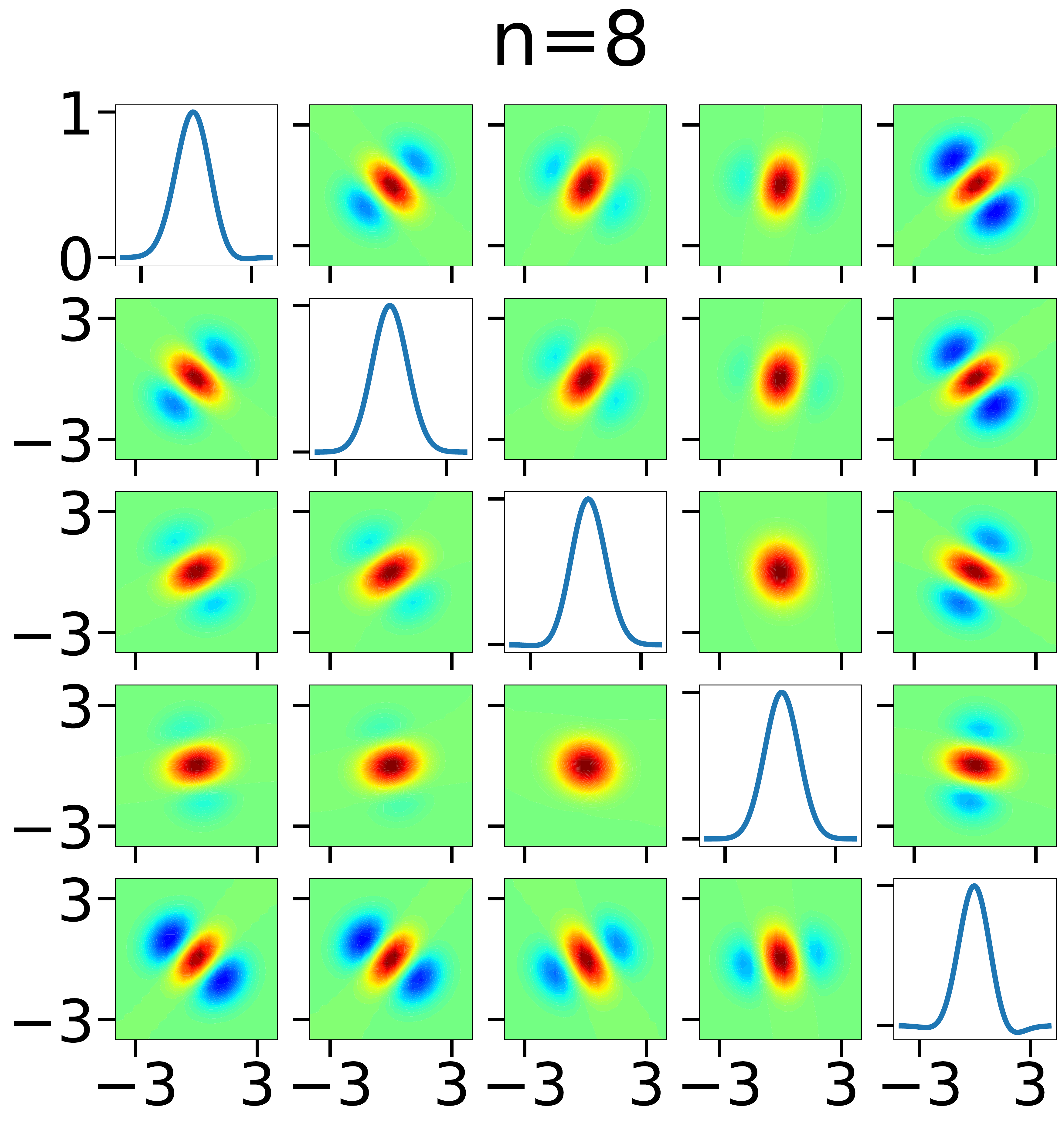}
\includegraphics[width=3.75cm,height=3.75cm]{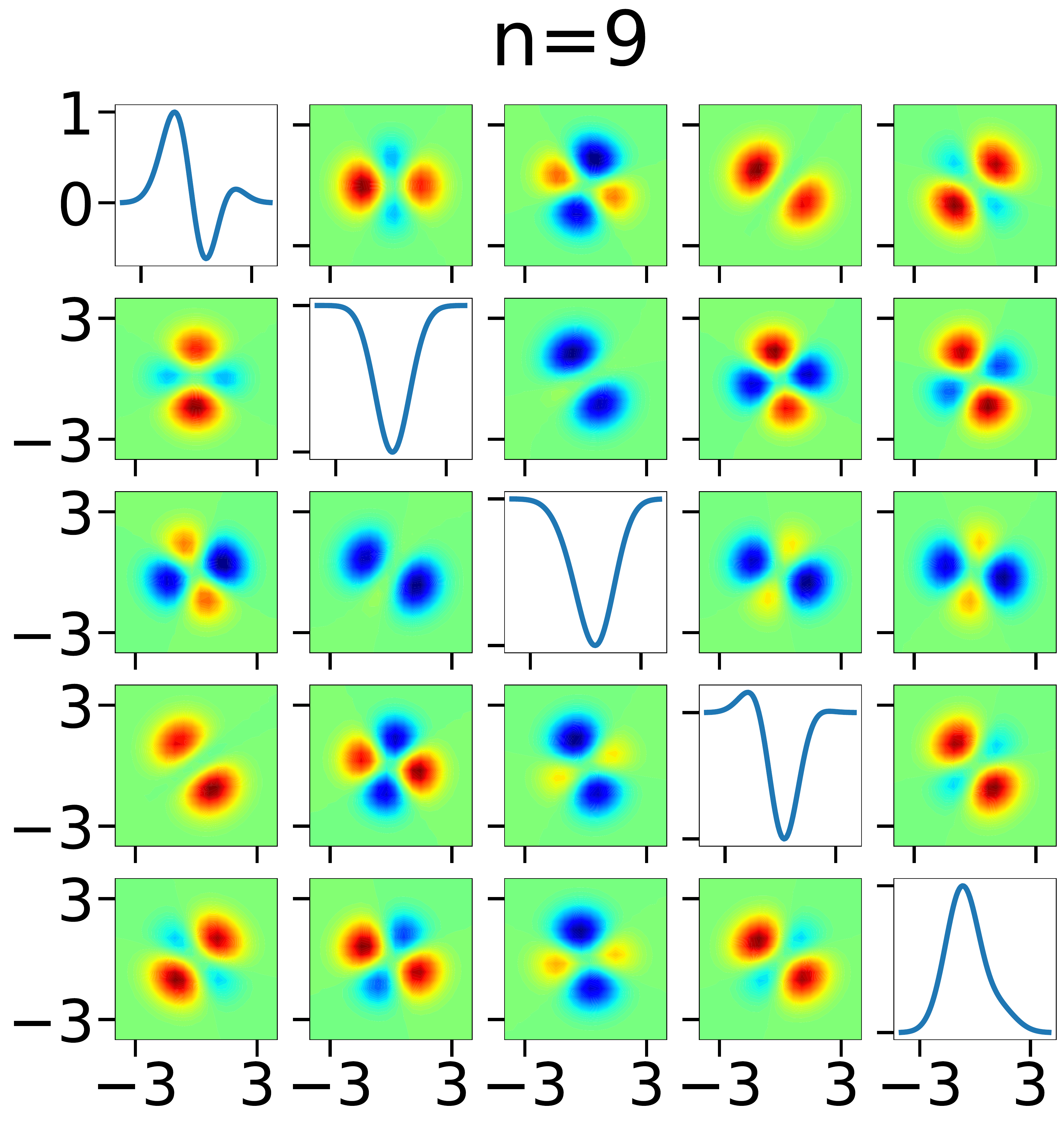}
\includegraphics[width=3.75cm,height=3.75cm]{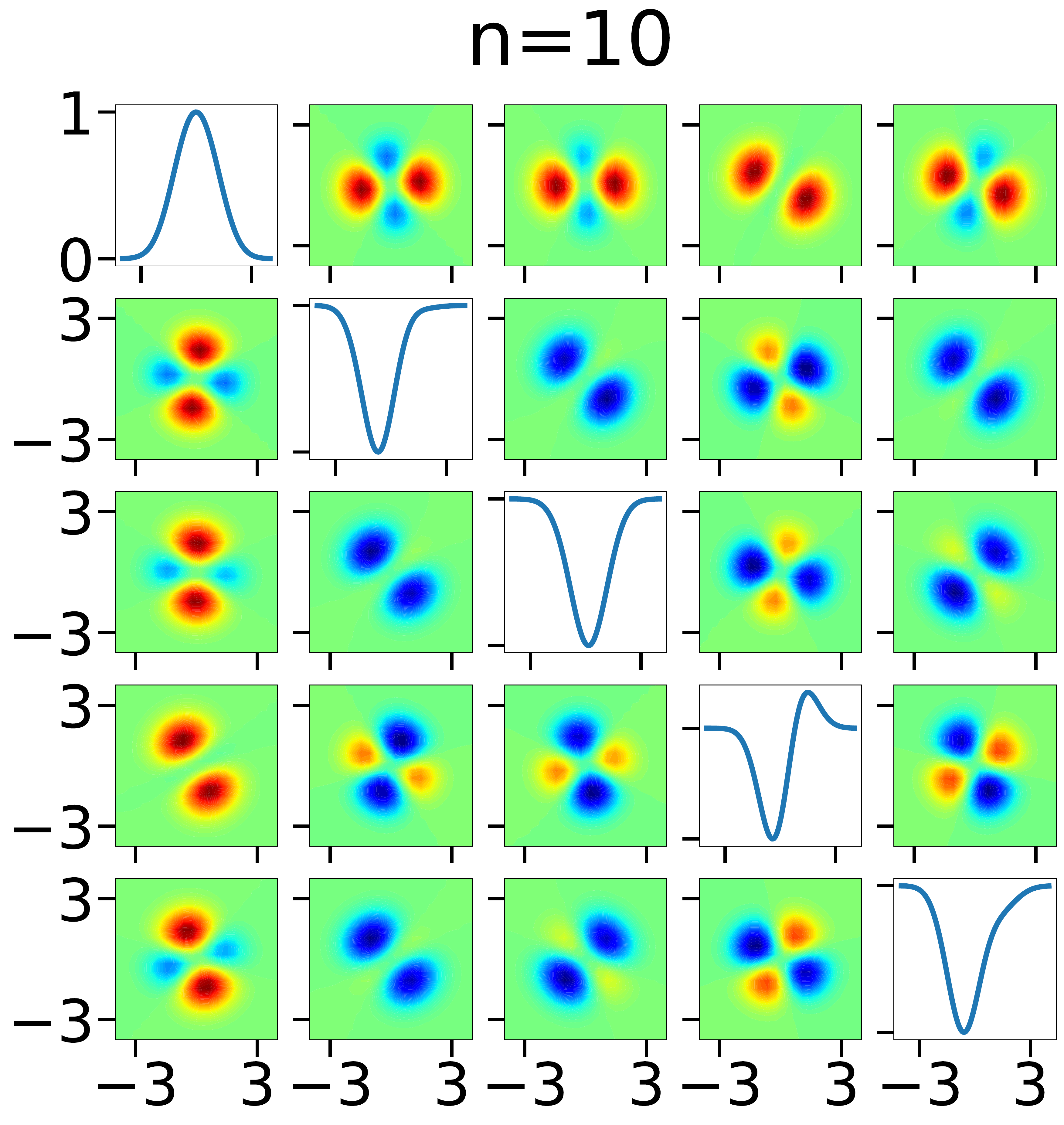}
\includegraphics[width=3.75cm,height=3.75cm]{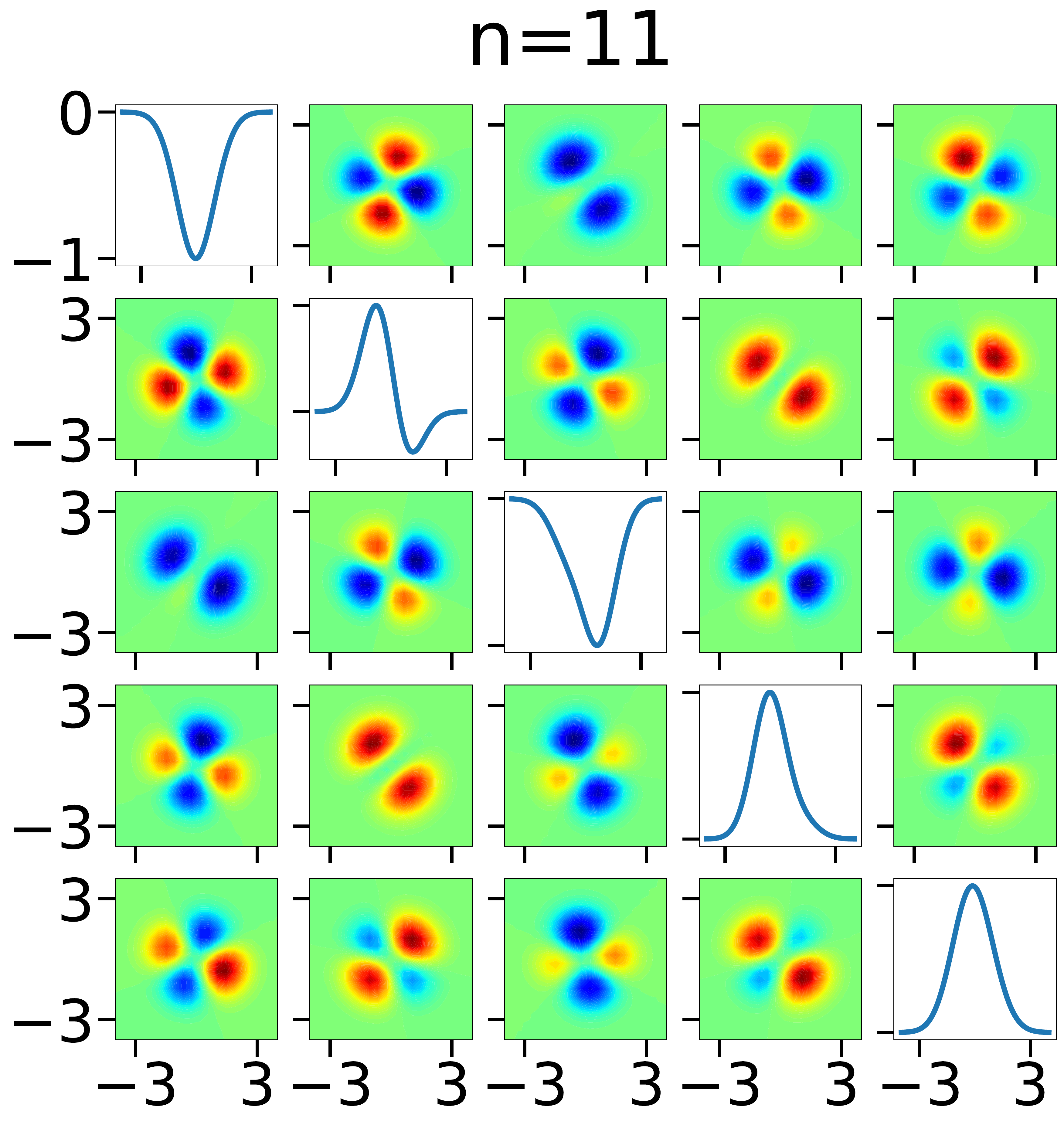}\\
\includegraphics[width=3.75cm,height=3.75cm]{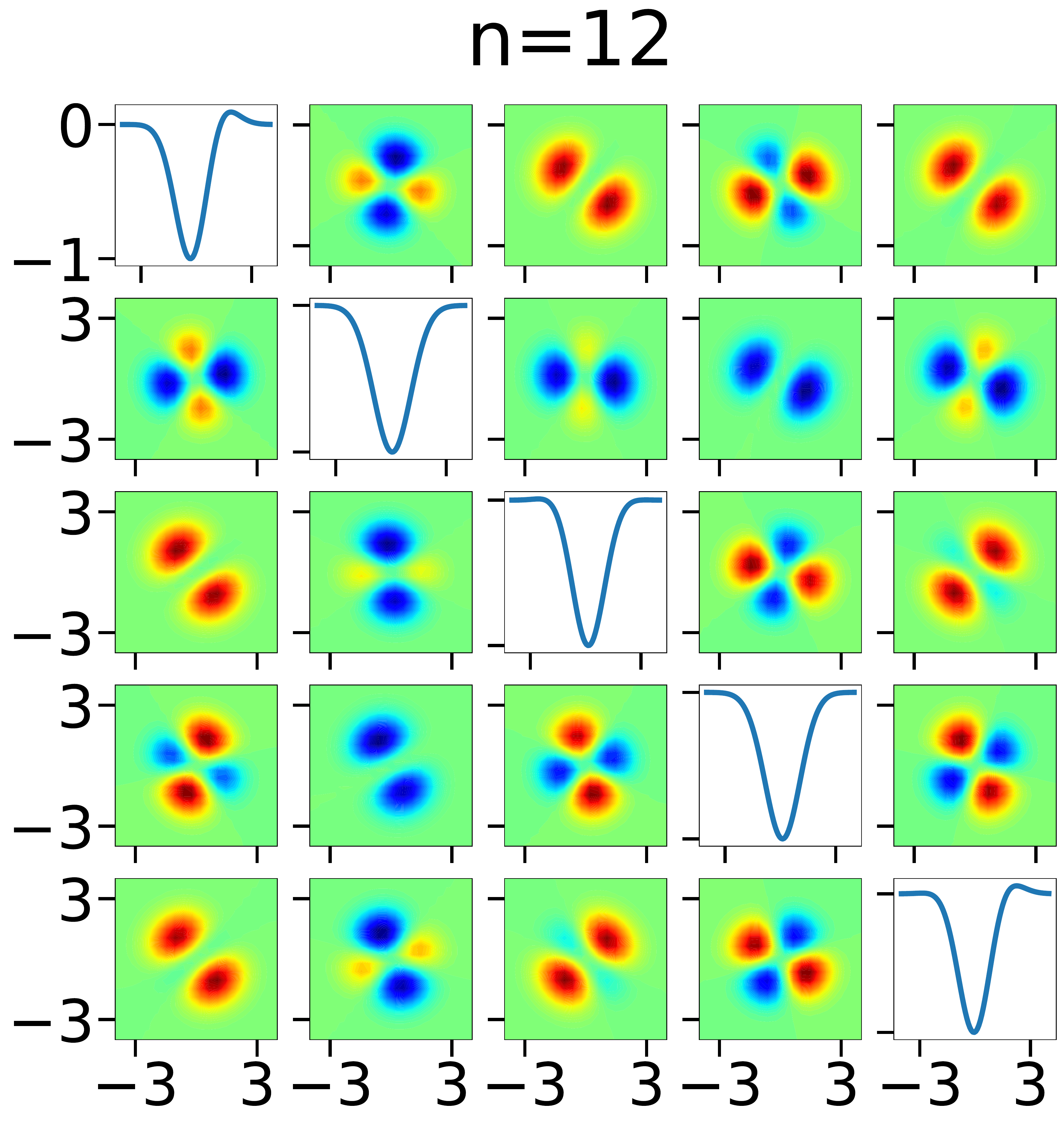}
\includegraphics[width=3.75cm,height=3.75cm]{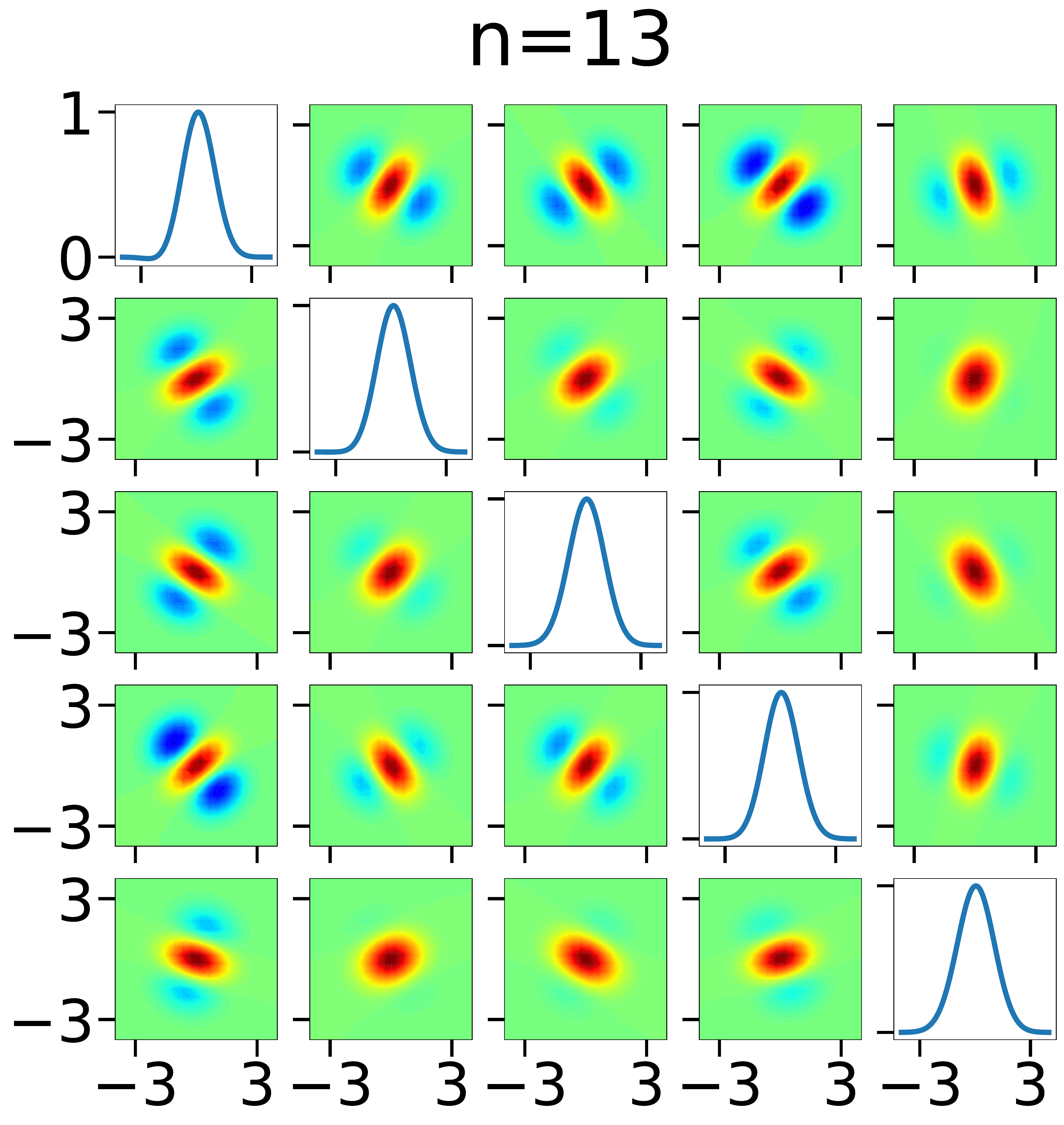}
\includegraphics[width=3.75cm,height=3.75cm]{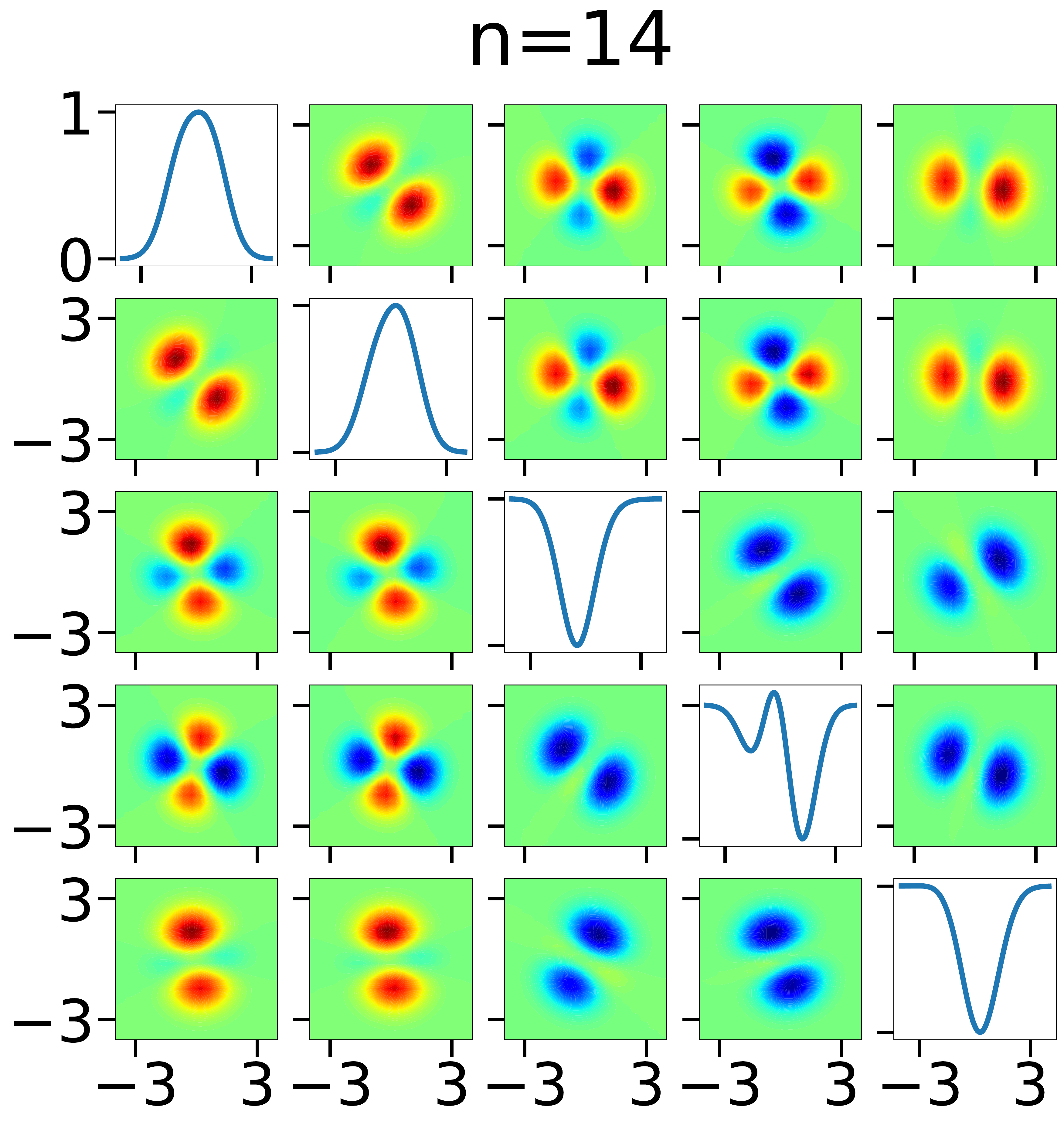}
\includegraphics[width=3.75cm,height=3.75cm]{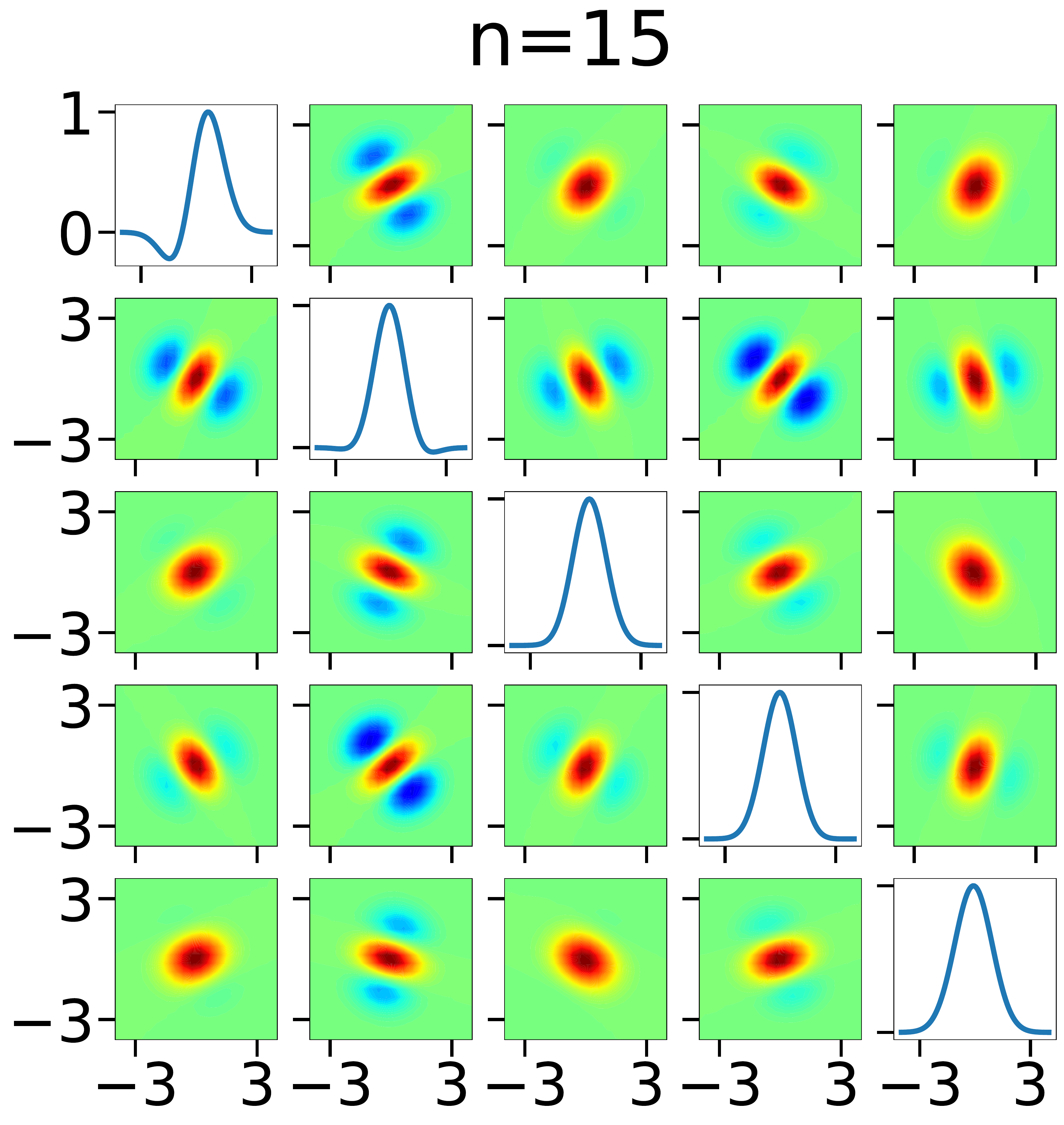}\\
\includegraphics[width=6cm,height=0.5cm]{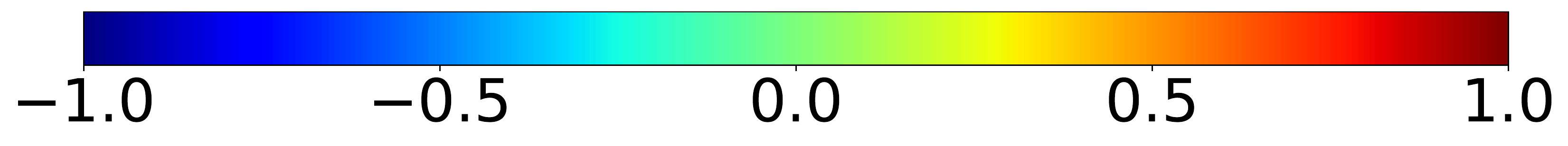}
\caption{The contour plots of the first 16 eigenfunctions for five-dimensional coupled harmonic
oscillator example.}\label{fig_5dCHO}
\end{figure}

Figure \ref{fig_5dCHO} shows the corresponding approximate wavefunctions obtained by
TNN-based machine learning method. From Figure \ref{fig_5dCHO},
the proposed numerical method here also has good accuracy for the five-dimensional eigenvalue problem of
the coupled harmonic oscillator.

\subsection{Energy states of hydrogen atom}
In this section, we study energy states of hydrogen atom.
The wave function $\Psi(x,y,z)$ of the hydrogen atom satisfies the following Schr\"{o}dinger equation
\begin{eqnarray}
-\frac{1}{2}\Delta\Psi-\frac{\Psi}{|\mathbf r|}=E\Psi,
\end{eqnarray}
where $|\mathbf r|=(x^2+y^2+z^2)^{1/2}$.
The exact energy of the hydrogen atom are $E_{n}=-\frac{1}{2n^2}$ and there are $n^2$ states consist with energy $E_n$.

In order to compute the singular integrals of the Column potential terms $1/|\mathbf r|$,
we adopt spherical coordinates $(r, \theta, \varphi)$ with density $r^2\sin\theta$.
Then the wave function $\Psi(\mathbf r)$ should be written as $\Psi(r,\theta, \varphi)$.
The Laplace $\Delta$ has following expression
\begin{eqnarray}
\Delta\Psi&=&\frac{\partial^2\Psi}{\partial r^2}+\frac{2}{r}\frac{\partial\Psi}{\partial r}+\frac{1}{r^2}\frac{\partial^2\Psi}{\partial\theta^2}
+\frac{\cos\theta}{r^2\sin\theta}\frac{\partial\Psi}{\partial\theta}
+\frac{1}{r^2\sin^2\theta}\frac{\partial^2\Psi}{\partial\varphi^2}\nonumber\\
&=&\frac{1}{r^2}\frac{\partial}{\partial r}\left(r^2\frac{\partial\Psi}{\partial r}\right)+\frac{1}{r^2\sin\theta}\frac{\partial}{\partial\theta}
\left(\sin\theta\frac{\partial\Psi}{\partial\theta}\right)
+\frac{1}{r^2\sin^2\theta}\frac{\partial^2\Psi}{\partial\varphi^2}.
\end{eqnarray}
The $99$ points Laguerre-Gauss quadrature is used in the direction $r$ and $16$ points
Legendre-Gauss quadrature with subintervals length $\frac{\pi}{64}$ in directions $\theta,\varphi$.
The TNN structure is defined as follows
\begin{eqnarray}\label{eq_TNN_H}
\Psi(r,\theta,\varphi)=\sum_{j=1}^pc_j\phi_{r,j}(\beta r)e^{-\frac{\beta r}{2}}\cdot\phi_{\theta,j}(\theta)\cdot\big(\phi_{\varphi,j}(\varphi)\sin(\varphi/2)+\gamma_j\big),
\end{eqnarray}
where $\phi_r=(\phi_{r,1},\cdots,\phi_{r,p})$, $\phi_\theta=(\phi_{\theta,1},\cdots,\phi_{\theta,p})$ and $\phi_\varphi=(\phi_{\varphi,1},\cdots,\phi_{\varphi,p})$ are three FNNs with depth $3$ and width $50$, and $p=20$.
The activation function is selected as $\sin(x)$.
The trainable parameter $\gamma_j$ is introduced to satisfy periodic boundary conditions
$\Psi(r,\theta,0)=\Psi(r,\theta,2\pi)$.

In implement, we use 15 TNNs to learn the lowest 15 energy states, each TNN is defined as (\ref{eq_TNN_H}).
The Adam optimizer is employed with a learning rate 0.0003 and epochs of 100000 and then the L-BFGS
in the subsequent 10000 steps to produce the final result.
Table \ref{table_H} shows the final energy approximations and corresponding errors.

\begin{table}[!htb]
\caption{Errors of energy states of hydrogen atom for the 15 lowest energy states.}\label{table_H}
\begin{center}
\begin{tabular}{ccccc}
\hline
$n$&  State&   Exact $E_n$&  Approx $E_n$&  ${\rm err}_E$\\
\hline
1&   $1s$&   $1/2$ &   -0.499999764803373&   -4.704e-07\\
2&   $2s$&   $1/8$ &   -0.124999995127689&   -3.898e-08\\
2&   $2p$&   $1/8$ &   -0.124999991163279&   -7.069e-08\\
2&   $2p$&   $1/8$ &   -0.124999974204279&   -2.064e-07\\
2&   $2p$&   $1/8$ &   -0.124999911448570&   -7.084e-07\\
3&   $3s$&   $1/18$&   -0.055555553061913&   -4.489e-08\\
3&   $3p$&   $1/18$&   -0.055555552597916&   -5.324e-08\\
3&   $3p$&   $1/18$&   -0.055555552115414&   -6.192e-08\\
3&   $3p$&   $1/18$&   -0.055555545116880&   -1.879e-07\\
3&   $3d$&   $1/18$&   -0.055555438704375&   -2.103e-06\\
3&   $3d$&   $1/18$&   -0.055555383940678&   -3.089e-06\\
3&   $3d$&   $1/18$&   -0.055555370476165&   -3.331e-06\\
3&   $3d$&   $1/18$&   -0.055555025332032&   -9.544e-06\\
3&   $3d$&   $1/18$&   -0.055554992074146&   -1.014e-05\\
4&   $4s$&   $1/32$&   -0.031249964498854&   -1.136e-06\\

\hline
\end{tabular}
\end{center}
\end{table}

\section{Conclusions}
In this paper, based on the deep Ritz method, we design a type of TNN-based
machine learning method to compute the leading
multi-eigenpairs of high dimensional eigenvalue problems. The most important advantage of
TNN is that the high dimensional integrations of TNN functions
can be calculated with high accuracy and efficiency.  Based on the high accuracy and efficiency of
the high dimensional integration, we can build the corresponding machine learning method for solving
high dimensional problems with the high accuracy.  The presented numerical examples show that the
proposed machine learning method in this paper can obtain obviously better accuracy than the
Monte-Carlo-based machine learning methods.

In our numerical implementation, we also find that the accuracy and stability of the machine learning
process should be paid more attention. These are necessary to get the final high accuracy for solving
high dimensional problems by using machine learning methods.

Actually, the proposed TNN and the corresponding machine learning method can be extended to other
high dimensional problems such as Schr\"{o}dinger equations, Boltzmann equations, Fokker-Planck equations,
stochastic equations, multiscale problems and so on.
This means TNN-based machine learning method can bring more practical applications in
physics, chemistry, biology, material science, engineering and so on. These will be our future work.

\end{document}